\documentclass{commatDV}

\usepackage{DLde}
\usepackage{graphicx}


\def\Om{\Omega}
\def\om{\omega}
\def\ph{\varphi}
\def\rh{\rho}
\def\ps{\psi}
\def\Ph{\Phi}
\def\al{\alpha}
\def\be{\beta}
\def\Ga{\Gamma}
\def\ga{\gamma}
\def\de{\delta}
\def\ep{\varepsilon}
\def\et{\eta}
\def\la{\lambda}
\def\La{\Lambda}
\def\si{\sigma}
\def\Si{\Sigma}
\def\ze{\zeta}
\def\ta{\tau}
\def\ka{\varkappa}
\def\ch{\chi}
\def\Th{\Theta}
\def\th{\theta}
\def\De{\Delta}
\renewcommand\de{\delta}
\def\Ad{\operatorname{Ad}}
\def\dif{\bar\partial}
\def\ad{\operatorname{ad}}
\def\sgn{\operatorname{sgn}}
\def\id{\operatorname{id}}
\def\SL{\operatorname{SL}}
\def\Hom{\operatorname{Hom}}
\def\GL{\operatorname{GL}}
\def\Gr{\operatorname{Gr}}
\def\gr{\operatorname{gr}}
\def\tr{\operatorname{tr}}
\def\x{\times}
\def\SU{\operatorname{SU}}

\renewcommand\S{\operatorname S}
\def\Ker{\operatorname{Ker}}
\def\pd#1{\nfrac{\partial}{\partial{#1}}}
\def\bpd{\bar\partial}
\def\hpd{\partial}
\def\M{(M,\mathcal O)}
\renewcommand\Im{\operatorname{Im}}
\def\Int{\operatorname{Int}}

\def\Bih{\operatorname{Bih}}
\renewcommand\i{^{-1}}
\def\row#1#2#3{#1_{#2},\ldots,#1_{#3}}
\def\Aut{\operatorname{Aut}}
\def\ev{\operatorname{ev}}

\def\rd{\operatorname{rd}}
\def\red{\operatorname{red}}
\renewcommand\S{\operatorname{S}}
\renewcommand\H{\operatorname{H}}
\def\I{\operatorname{I}}
\def\odd{\operatorname{odd}}
\def\diag{\operatorname{diag}}
\def\Alt{\operatorname{Alt}}
\def\Der{\mathfrak{der}}

\setcounter{tocdepth}{2}

\makeatletter
\renewcommand {\ssbegin}[2][*]
  {\refstepcounter{subsection}
\if*#1
\addcontentsline{toc}{subsection}{\protect\numberline{\thesubsection}#2}%
\else
\addcontentsline{toc}{subsection}{\protect\numberline{\thesubsection}#2\ (#1)}%
\fi
  \def \secno {\gdef \secno {}{\ssecfont \thesubsection.\hskip 2ex}%
  }%
\if*#1
   \begin{#2}
\else
\begin{#2}[#1]
\fi}

\renewcommand {\sssbegin}[1]
  {\refstepcounter{subsubsection}
\addcontentsline{toc}{subsubsection}{\protect\numberline{\thesubsubsection}#1}%
  \def \secno {\gdef \secno {}{\ssecfont \thesubsubsection.\hskip 2ex}%
  }%
   \begin{#1}}

\makeatother

\title{Non-split supermanifolds associated with the cotangent bundle}

\author[Arkady Onishchik]{\fbox{Arkady Onishchik}}

\affiliation{N/A} 



\abstract{Here, I study the problem of classification of non-split supermanifolds
having as retract the split supermanifold $(M,\Omega)$, where $\Omega$ is the
sheaf of holomorphic forms on a given complex manifold $M$ of dimension
$> 1$. I propose a
general construction associating with any $d$-closed $(1,1)$-form $\omega$ on $M$
a supermanifold with retract $(M,\Omega)$ which is non-split whenever the
Dolbeault class of $\omega$ is non-zero. In particular, this gives a non-empty
family of non-split supermanifolds for any flag manifold $M\ne
\mathbb{CP}^1$. In the case where $M$ is
an irreducible compact Hermitian symmetric space, I get a~complete
classification of non-split supermanifolds with retract $(M,\Omega)$. For each
of these supermanifolds, the 0- and 1-cohomology with values in the tangent
sheaf are calculated. As an example, I study the $\Pi$-symmetric
super-Grassmannians introduced by Yu. Manin.
}

\VOLUME{30}
\NUMBER{3}
\firstpage{49}
\DOI{https://doi.org/10.46298/cm.9613}

\begin{paper}    

\tableofcontents


\subsection*{Editor's note} The work on the problem of deformations of analytical supermanifolds, especially the study of non-splitness, was dormant for ca 40 years after the discovery of the first non-split example by Green. Apart from Onishchik and his students nobody studied this problem. The problem  drew new attention of both mathematicians and theoretical physicists after Donagi and Witten showed that the moduli space of super Riemann surfaces is not split and how this fact affects working with modules of string theory, see \cite{DW}.

This work by Onishchik was preprinted in 1997 as Pr\'epubl. Univ. Poitiers D\'epart de Math. N.~109. It  was found very helpful several times since, see \cite{V6} -- \cite{V2} and~ \cite{BV}. 

A.L.Onishchik used to tell me that he understood the meaning of the ``non-split" supermanifold having learned the papers by Green \cite{9} and Vaintrob \cite{33}, \cite{34}. Here, I updated the references; my additions are marked with a~ $\ast$. I also considerably edited English, but each time I read this text after a break I see something to be corrected, so I am afraid there still is something left. 
I also changed the outdated or \textit{ad hoc} notation of several supergroups and superalgebras using the currently used notation and inserted a couple of clarifying parenthetical remarks (marked by \textit{D.L.}). 

For basics of supermanifold theory, I recommend \cite{Del} and \cite{Lsos} which still contain many results, notions and open problems not covered in other sources; see also comments in \cite[Section~4.8]{Mo}. I also recommend the wonderful introduction into the theory of schemes and ringed spaces \cite{MaAG}, and the definition and calculation of curvature tensors of almost complex supermanifolds, and real-complex supermanifolds endowed with a non-integrable distribution (see \cite{BGLS}), examples of such supermanifolds are all superstrings usually considered in the works of physicists and most of the super Grassmannians. 

I am thankful to E.~Vishnyakova for her help in editing this manuscript. 

In what follows, ``I" means ``Onishchik''.~\textit{D.Leites.}

\section{Introduction}\label{S0}

One of the most important features of the theory of complex analytic
supermanifolds is the existence of non-split supermanifolds. The simplest
example is the superquadric $\mathcal{Q}^{1|2}$ in the
projective superplane $\mathcal{CP}^{2|2}$, see Example \ref{E1.10} below; it is of dimension $1|2$ and has as
its base the projective line $\mathbb{CP}^1$. This superquadric
belongs to one of four series of homogeneous complex supermanifolds
constructed by Yu. Manin \cite{20} --- the flag supermanifolds; as a~rule, they are non-split.

With any supermanifold a~split one, called its \textit{retract}, is associated. In
this paper, I study non-split supermanifolds with retract $(M,\Om)$, where
$\Om$ is the sheaf of holomorphic forms on a~complex manifold $M$. I~
present a~construction assigning to any $d$-closed (1,1)-form $\om$ on $M$
a supermanifold with retract $(M,\Om)$; this supermanifold is non-split whenever $\om$
has a~non-zero Dolbeault cohomology class. In particular, for any
compact K\"ahler manifold $M$, we obtain a~family of supermanifolds with
retract $(M,\Om)$ parametrized by $H^{1,1}(M,\mathbb C)$, all members of which are non-split,
except the one corresponding to 0. This family is non-empty,
e.g., when $M$ is a~flag manifold.

The next problem is the classification of \textit{all} non-split supermanifolds
with retract $(M,\Om)$, where $M$ is a~flag manifold. I solve it in the case
where $M$ is an irreducible Hermitian symmetric space. In this case, the
family mentioned above contains precisely one non-split supermanifold. I
prove that this is the only non-split supermanifold with retract $(M,\Om)$ if
one excludes the case of the Grassmannians $M = \Gr^{n}_{s}$, where $2\le s\le n-2$,
while the non-split supermanifolds for such a~Grassmannian form an
1-parameter family. The proof is based on certain general results
concerning classification of supermanifolds with a~given retract. We also
calculate the 0- and 1-cohomology of the tangent sheaf for all the
supermanifolds associated with the cotangent bundle over a~compact
irreducible Hermitian symmetric space.

The well known examples of supermanifolds studied here are the
supermanifolds of $\Pi$-symmetric flags that form one of Manin's series
mentioned above. If $M$ is a~symmetric space, then $M = \Gr^{n}_{s}$, and we have
the $\Pi$-symmetric super-Grassmannians $\Pi\Gr^{n|n}_{s|s}$. As a~corollary,
we calculate the 0- and 1-cohomology of their tangent sheaf. Note that
this special question initiated the study exposed in this paper. For the three
other series of super-Grassmannians this calculation was performed in \cite{22}, \cite{26},
\cite{27}.

\section{Generalities about superalgebras and supermanifolds}\label{S1}

\subsection{Algebraic background}\label{ss1.1}

To fix the notation, we give here some definitions. 

Let $\mathbb Z_2$ denote not the ring of 2-adic integers but $\mathbb Z/2\mathbb Z = \{\bar 0,\bar 1\}$. A \textit{vector superspace} is any $\mathbb Z_2$-graded vector space $V$.
In this paper, the ground field is the field of complex numbers $\mathbb C$,
By definition, we have
$$
V = V_{\bar 0}\oplus V_{\bar 1},
$$
where $V_{\bar 0},\; V_{\bar 1}$ are vector subspaces called the \textit{even
part} and the \textit{odd part} of $V$, respectively. The \textit{non-zero} elements of these
subspaces are said to be \textit{even} or \textit{odd}, respectively, and we define the {\it
parity function} by setting for \textit{non-zero} elements
$$
p(x) = \begin{cases} \bar 0\text{ if $ x\in V_{\bar 0}$}\\ \od\text{ if $x\in V_{\bar 1}$}.
\end{cases}
$$

We write $\dim V = n|m$, where $\dim V_{\bar 0} = n,\; \dim V_{\bar 1} = m$;
this
is the \textit{superdimension} of a~vector superspace $V$. A standard example of a
vector superspace of dimension $n|m$ is $\mathbb C^{n|m} = \mathbb C^n\oplus \Pi
\mathbb C^m$, where $\Pi$ is the change of parity functor.

A \textit{superalgebra} is a~$\mathbb Z_2$-graded algebra over $\mathbb C$, i.e.,  a
vector superspace $A = A_{\bar 0}\oplus A_{\bar 1}$ endowed with a
multiplication $(a,b)\mapsto ab$ satisfying the following condition
$$
A_iA_j\subset A_{i+j}\text{~~for any~~} i,j\in\mathbb Z_2.
$$
A \textit{morphism} $f: A\tto B$ of superalgebras is, by definition, a~\textit{parity preserving}
homomorphism of algebras, i.e.,  satisfying $f(A_i)
\subset B_i$ for any $i\in\mathbb Z_2$. If $A$ and
$B$ are superalgebras with units, we also assume that $f(1) = 1$.

Let $A = \bigoplus_{i\in\mathbb Z} A_i$ be a~graded (i.e.,  a~$\mathbb Z$-graded)
algebra. This means that
$$
A_iA_j\subset A_{i+j}\text{~~for any~~} i,j\in\mathbb Z.
$$
Setting
$$
A_{\bar 0} = \bigoplus_{i\in\mathbb Z} A_{2i},\;
A_{\bar 1} = \bigoplus_{i\in\mathbb Z} A_{2i+1},
$$
we, clearly, endow $A$ with a~parity ($\mathbb Z_2$-grading), turning it into a
superalgebra. In this case, we say that the $\mathbb Z$-grading and the
$\mathbb Z_2$-grading in $A$ are \textit{compatible}.

\ssbegin[The exterior a.k.a. Grassmann algebra]{Example}\label{E1.2}
Let $E$ denote a~complex vector space of dimension $m$ and let $\bigwedge
E$ be the exterior (or Grassmann) algebra over $E$. Then, we have the natural
$\mathbb Z$-grading
$$
\bigwedge E = \bigoplus_{p=0}^n\bigwedge^p E
$$
making $\bigwedge E$ a~graded algebra. Using the above procedure, we can
regard $\bigwedge E$ as a~superalgebra. Note that setting $p(x) = \od$ for any $x\in E$ we endow $\bigwedge E$ with a~(natural) parity.
Any basis $\xi_1,\ldots,\xi_m$ of $E$ is a~set of (odd) generators of
$\bigwedge E$, and we often write
$$
\bigwedge E = \bigwedge_{\mathbb C}(\xi_1,\ldots,\xi_m).
$$
\end{Example}

Many formulas of Linear Algebra are superized by means of the Sign Rule \lq\lq \textbf{if something of parity $a$ is moved past something of parity $b$, the sign $(-1)^{ab}$ accrues; formulas defined on homogeneous elements are extended to all elements via linearity}". Here are examples.

A superalgebra $A$ is called \textit{supercommutative} if
$$
ba = (-1)^{p(a)p(b)}ab
$$
for any even or odd $a,b\in A$. The associativity of $A$ is meant in the
usual sense. Clearly, $\bigwedge E$ from Example \ref{E1.2} is an associative
supercommutative superalgebra with unit.

Let $V$ and $W$ be two vector superspaces. Then, the tensor product $U = V\otimes
W$ becomes a~superspace if we endow it with the following $\mathbb Z_2$-grading:
$$
U_{\bar 0} = V_{\bar 0}\otimes W_{\bar 0}\oplus V_{\bar 1}\otimes W_{\bar 1},
\quad
U_{\bar 1} = V_{\bar 0}\otimes W_{\bar 1}\oplus V_{\bar 1}\otimes W_{\bar 0},
$$
For any two superalgebras $A$ and $B$, let us endow the superspace $A\otimes
B$ with the multiplication
$$
(a_1\otimes b_1)(a_2\otimes b_2) = (-1)^{p(b_1)p(a_2)}(a_1 a_2)\otimes
(b_1 b_2),\quad a_i\in A,\; b_i\in B.
$$
Then, $A\otimes B$ is a~superalgebra (the \textit{tensor product} of $A$ and
$B$). The tensor product of two associative (supercommutative) superalgebras is
associative (respectively, supercommutative).

Let $V$ and $W$ be two vector superspaces. Then, the vector space $L(V,W)$ of all
linear mappings $V\tto W$ is endowed with the following $\mathbb Z_2$-grading:
$$
L(V,W)_k = \{f\in L(V,W) \mid f(V_i)\subset W_{i+k},\; i\in\mathbb Z_2\},\; k\in
\mathbb Z_2.
$$
Thus, a~non-zero $f\in L(V,W)$ is even (odd) if it preserves (respectively, changes)
the parity. For example, any morphism of superalgebras is, by definition,
even.

Regarding $\mathbb C$ as $\mathbb C^{1|0}$, we get a~natural $\mathbb Z_2$-grading in
the dual vector space $V^* = L(V,\mathbb C)$ of a~superspace $V$. Clearly,
\[
\text{$(V^*)_{\bar 0} = \{f\in V^* \mid f(V_{\bar 1}) = 0\}$ and
$(V^*)_{\bar 1} = \{f\in V^* \mid f(V_{\bar 0}) = 0\}$ }
\]
are identified with
$(V_{\bar 0})^*$ and $(V_{\bar 1})^*$, respectively. For another vector
superspace $W$, the superspace $V^*\otimes W$ is identified with $L(V,W)$,
as usual.

\ssbegin[Associative superalgebra of supermatrices]{Example}\label{E1.3}
Let $V$ be a~vector superspace. Then, $L(V) = L(V,V)$, the associative algebra  with unit 
of all linear transformations of $V$ is a~superalgebra if we endow
it with the above $\mathbb Z_2$-grading.

A corresponding example can be constructed by means of matrices. Let
$\mathbb M_{n+m}(\mathbb C)$
denote the (associative, with unit) algebra of $(n+m)\x(n+m)$ matrices over
$\mathbb C$. Endow this algebra with a~$\mathbb Z_2$-grading in the following way.
Write a~matrix $X\in \mathbb M_{n+m}(\mathbb C)$ in the form
$$
X =\begin{pmatrix}
X_{00} & X_{01} \\
X_{10} & X_{11}
\end{pmatrix},
$$
where $X_{00}$ fills the first $n$ rows and $n$ columns. Then, $p(X) = \bar 0$
whenever $X_{01} = 0,\; X_{10} = 0$, and $p(X) = \od$ whenever $X_{00} = 0,\;
X_{11} = 0$. Clearly, this $\mathbb Z_2$-grading endows $\mathbb M_{n+m}(\mathbb C)$ with a
superalgebra structure. We denote it by $\mathbb M_{n|m}(\mathbb C)$.

Let $V$ be a
vector superspace of dimension $n|m$. Choosing in $V$ a~basis 
\[
\text{$\row e1n,
\row f1m$ with $e_i\in V_{\bar 0},\; f_j\in V_{\bar 1}$,}
\]
 we get a
natural isomorphism of superalgebras $L(V)\simeq\mathbb M_{n|m}(\mathbb C)$.
\end{Example}

Let $\mathfrak g$ be a~superalgebra and let us agree to denote the multiplication
in $\mathfrak g$ by $[-,-]$ and to call it the \textit{bracket}. We say that
$\mathfrak g$ is a~\textit{Lie superalgebra} if the following conditions are
satisfied for any $x, y, z\in \mathfrak{g}$:
$$
\begin{aligned}
{}[y,x] &= -(-1)^{p(x)p(y)}[x,y],\\
{}[x,[y,z]] &= [[x,y],z] + (-1)^{p(x)p(y)}[y,[x,z]].
\end{aligned}
$$
A graded algebra $\mathfrak g$ is called a~\textit{graded Lie superalgebra} if
$\mathfrak g$ is a~Lie
superalgebra being provided with a $\mathbb Z_2$-grading (parity). Observe that these two gradings are not necessarily compatible.

\ssbegin[Lie superalgebra of supermatrices]{Example} \label{E1.4}
As in the classical case, there is a~functor $\mathfrak L$ converting associative
superalgebras into Lie ones. If $A$ is an associative superalgebra, then
$\mathfrak L(A)$ is the vector superspace $A$ endowed with the bracket
$$
[a,b] = ab - (-1)^{p(a)p(b)}ba)\text{~~for any $a,b\in A$}.
$$

Let us note some special cases.

If $A$ is supercommutative, then $\mathfrak L(A)$ is a~Lie superalgebra with zero
bracket called \textit{commutative} Lie superalgebra.

For any vector superspace $V$, the superalgebra $L(V)$ from Example \ref{E1.2}
gives the general linear Lie superalgebra $\mathfrak{gl}(V) = \mathfrak L(L(V))$.

By the same example, we get the Lie superalgebra $\mathfrak{gl}_{n|m}(\mathbb C) =
\mathfrak L(\mathbb M_{n|m}(\mathbb C))$.
\end{Example}

\ssbegin[Lie superalgebra of superderivations]{Example}\label{E1.5}
Let $A$ be an arbitrary superalgebra. A 
linear transformation
$v\in\mathfrak{gl}(A)$ is called a~\textit{derivation} of $A$ if
$$
v(ab) = v(a)b + (-1)^{p(v)p(a)}av(b)\text{~~for any $a,b\in A$}.
$$
Denote
$$
\mathfrak{der}\,A := (\mathfrak{der}\,A)_{\bar 0}\oplus (\mathfrak{der}\,A)_{\bar 1},
$$
where $(\mathfrak{der}\,A)_i\subset\mathfrak{gl}(A)_i,\; i\in\mathbb Z_2$, is the
vector space of even or odd derivations of $A$. One checks easily that
$\mathfrak{der}\,A$ is a~subalgebra of $\mathfrak{gl}(A)$ and, hence, a~Lie
superalgebra (the \textit{superalgebra of derivations} of $A$).
\end{Example}

Let $\mathfrak g$ be a~Lie superalgebra. For any $x\in L$, define the {\it
adjoint operator} $\ad_x\in L(\mathfrak g)$ by setting
$$
\ad_x\,(y) := [x,y]\text{~~for any~~}y\in\mathfrak g.
$$

A straightforward verification shows that $\ad_x\in\mathfrak{der}\,\mathfrak g$ for
any $x\in\mathfrak g$ and  
${\ad: \mathfrak g\tto\mathfrak{der}\,\mathfrak g}$
is a~homomorphism of Lie superalgebras.

Clearly, $\mathfrak g_{\bar 0}$ is a~usual Lie algebra. If $x\in\mathfrak g_{\bar 0}$,
then $\ad_x$ is an even derivation. Restricting it onto $\mathfrak g_p$, where $p =
\bar 0,\bar 1$,
we get two linear representations $\ad_p$ of $\mathfrak g_{\bar 0}$ in $\mathfrak g_p$, where $p = \bar 0,\bar 1$.

Similarly, let $\mathfrak g = \bigoplus_{p\in\mathbb Z}\mathfrak g_p$ be a~graded Lie
superalgebra. Then, $\mathfrak g_0$ is a~Lie algebra, and the restriction $\ad|_{
\mathfrak g_0}$ determines a~representation $\ad_p$ of $\mathfrak g_0$ in $\mathfrak g_p$
for any $p\in\mathbb Z$.

If $A$ is an associative and supercommutative superalgebra, then $\Der\,A$ is naturally
an $A$-module according to the rule
$$
(au)(b) = au(b)\text{~~for any~~}u\in\Der\,A,\;a,b\in A.
$$

The same definitions apply to sheaves of graded algebras on a~topological
space $M$. If, in particular, $\mathcal A$ is a~sheaf of associative supercommutative
graded algebras, then the sheaf
$\mathcal Der\mathcal A$ of derivations of $\mathcal A$ is defined, which is a~sheaf of
Lie superalgebras and a~sheaf of $\mathcal A$-modules on $M$.

Let us consider the case where $A$ is a~graded algebra and regard it as a
superalgebra with respect to the compatible $\mathbb Z_2$-grading. Then, 
$\Der\,A$ has a~natural structure of the graded Lie superalgebra.
More precisely, $\Der\,A = \sum_{p\in\mathbb Z}(\Der\,A)_p$, where
$$
(\Der\,A)_p = \{\de\in\Der\,A \mid \de(A_q)\subset A_{q+p}\ \ \text{ for any }
q\in\mathbb Z\}.
$$
One can always define the \textit{grading derivation} $\ep\in (\Der\,A)_0$
given by the formula
\begin{equation}
\ep(f) = pf\text{ for any$f\in A_p$ and $p\ge 0$}.\label{(1.1)}
\end{equation}
One easily checks that
\begin{equation}
[\ep,v] = pv\text{ for any } v\in(\Der\,A)_p,\; p\in\mathbb Z.\label{(1.2)}
\end{equation}

\ssbegin[Vectorial Lie superalgebras]{Example}\label{E1.6}
Consider a~complex vector space $E$ of dimension $m$ and its corresponding the Grassmann algebra $A = \bigwedge E$ (see Example \ref{E1.2}).
Denote $W(E) = \Der\,A$. These Lie superalgebras constitute one of the
``Cartan type" series of simple (for $m\ge 2$) finite-dimensional vectorial Lie
superalgebras (see \cite{15}).

We need the well known description of derivations from $W(E)$ in terms of
multilinear
forms. Any $u\in W(E)_p$ is determined by its restriction onto $E = A_1$,
which can be an arbitrary linear mapping $E\tto A_{p+1} = \bigwedge^{p+1} E$.
Thus, $W(E)_p$ is isomorphic, as a~vector space, to $L(E,\bigwedge^{p+1}E)$,
which can be identified with $\bigwedge^{p+1}E\otimes E^*$. Let us denote by
$i(\ph)\in W(E)_p$ the derivation corresponding to a~linear mapping
$\ph\in\bigwedge^{p+1}E\otimes E^*$.

The elements of the latter vector space can be regarded as vector-valued
$(p+1)$-forms on $E^*$, i.e.,  as anti-symmetric
$(p+1)$-linear mappings $(E^*)^{p+1}\tto E^*$. Regarding $A$ as the set of all
anti-symmetric multilinear forms on $E^*$, we have
$$
\begin{aligned}
i(\ph)&(a)(x_1,\ldots,x_{p+q}) \\ &=\nfrac1{(p+1)!(q-1)!}
\sum_{\al\in S_{p+q}}(\sgn\al)
a(\ph (x_{\al_1},\ldots,x_{\al_{p+1}}),x_{\al_{p+2}},\ldots,x_{\al_{p+q}})\text{~~for the $x_k\in E^*$.}
\end{aligned}
$$

Denote by $\xi_j$ for $j = 1,\ldots,m$ a~basis of $E$, and by $\xi_j^*$ for $j =
1,\ldots,m$ the dual basis of $E^*$.
Clearly, the derivations $\pd{\xi_j} = i(\xi_j^*)\in W(E)_{-1}$ for $ j =
1,\ldots, m$, constitute a~basis of the $A$-module $W(E)$. It follows that
the derivations
$$
\xi_{i_1}\ldots\xi_{i_{p+1}}\pd{\xi_j}\text{~ for $i_1<\ldots<i_{p+1}$ and $j = 1,\ldots,
m,$}
$$
constitute a~basis of $W(E)_p$ over $\mathbb C$. In particular, we see that
$W(E)_p$ is non-zero only for $-1\le p\le m$.

We also write 
$$
i(\ph)(a) = a\barwedge\ph\text{~~ for any $a\in A$ and $\ph\in A\otimes E^*$}.
$$
A similar operation can be defined for two vector-valued forms of arbitrary
degrees. Namely, let
$\ph\in A_p\otimes E^*$ and $\psi\in A_q\otimes E^*$ be given. Regarding these
tensors as $E^*$-valued $p$- and $q$-forms on $E^*$, we define the form
$\ph\barwedge\psi\in A_{p+q-1}\otimes E^*$ by the formula
\begin{equation}
\begin{aligned}
(\ph&\barwedge\psi)(x_1,\ldots,x_{p+q-1}) \\
&=\nfrac1{(p-1)!q!}\sum_{\al\in S_{p+q-1}}(\sgn\al)
\ph(\psi (x_{\al_1},\ldots,x_{\al_q}),x_{\al_{q+1}},\ldots,x_{\al_{p+q-1}})
\end{aligned}
\label{(1.3)}
\end{equation}
for any $x_k\in E^*$. This operation can be used to express the bracket in
$W(E)$. More precisely, define the bilinear operation $\{-,-\}$ on
$A\otimes E^*$ by setting
\begin{equation}
\{\ph,\psi\} = \psi\barwedge\ph - (-1)^{(p-1)(q-1)}\ph\barwedge\psi\label{(1.4)}
\end{equation}
for any $\ph\in A_p\otimes E^*$ and $\psi\in A_q\otimes E^*$. Then, 
$$
i(\{\ph,\psi\}) = [i(\ph),i(\psi)].
$$

In what follows, we will use the linear mapping $j: \bigwedge^pE\to
L(E,\bigwedge^{p+1}E)$ given by the formula
\begin{equation}
j(\ps)(u) = \ps u\text{~~for any~~}u\in E.\label{(1.5)}
\end{equation}
It is injective whenever $p < m$. Clearly,
\begin{equation}
i(j(\ps)) = \ps\ep,\ \ \ps\in\bigwedge^pE.\label{(1.6)}
\end{equation}
Regarding $L(E,\bigwedge^{p+1}E)$ as $\bigwedge^{p+1}E\otimes E^*$, we easily
see that
\begin{equation}
j(\ps) = \sum_{k=1}^m(\ps\xi_k)\otimes \xi_k^*.\label{(1.7)}
\end{equation}
Finally, regarding elements of $L(E,\bigwedge^{p+1}E)$ as vector-valued
$(p+1)$-forms on $E^*$, we obtain
\begin{equation}
j(\ps)(x_1,\ldots, x_{p+1}) = p!\sum_{k=1}^m(-1)^{k-1}\ps(x_1,\ldots,
\hat x_k,\ldots, x_{p+1})x_k,\text{~~where~~}x_l\in E^*.\label{(1.8)}
\end{equation}

On the other hand, there is the contraction mapping $c: \bigwedge^{p+1}E
\otimes E^*\tto\bigwedge^pE$ given by the formula
$$
c(\ph)(x_1,\ldots,x_p) = \sum_{k=1}^m\ph(\xi_k^*,x_1,\ldots,x_p)(\xi_k),
\text{~~where~~} x_l\in E^*.
$$
An easy calculation shows that $cj = p!(m-p)\id$. It follows that
\begin{equation}
\bigwedge^{p+1}E\otimes E^* = \Im j\oplus\Ker c\text{~~whenever $p < m$.}\label{(1.9)}
\end{equation}
\end{Example}

\subsection{Complex supermanifolds}\label{ss1.7}
The word ``supermanifold" will mean the same as in \cite{2}, \cite{2b}, \cite{18}, but
the complex-analytic version of the theory will be considered (see \cite{20}).
Let us begin with a~more general notion of the ringed space.

A $\mathbb Z_2$-\textit{graded ringed space} is a~pair $\M$, where $M$ is a
topological space and $\mathcal O$ is a~sheaf of associative unital supercommutative
superalgebras on $M$. A \textit{morphism} between two $\mathbb Z_2$-graded ringed spaces $(M,\mathcal O_M) \to
(N,\mathcal O_N)$ is a~pair $(f, f^*)$, where $f: M\to N$ is a~continuous mapping and $f^*: \mathcal O_N \to \mathcal O_M$ a~morphism of sheaves of superalgebras. In particular, if $F = (f,f^*)$ is an automorphism of a~ringed space $\M$, then we can consider the mapping $f_* = (f^*)\i$ instead of $f^*$; this is an automorphism of the sheaf $\mathcal O$ over $M$. The automorphisms of $\M$ form the group $\Aut\M$.

\ssbegin[Complex-analytic supermanifolds]{Example}\label{E1.8}
On the space $\mathbb C^n$, consider the sheaf
$$
\mathcal F_{n|m} :=
\bigwedge(\row\xi 1m)\otimes\mathcal F_n = \bigwedge_{\mathcal F_n}(\row\xi 1m),
$$
where $\mathcal F_n$ is the sheaf of germs of holomorphic functions on
$\mathbb C^n$. Here we assume that the functions from $\mathcal F_n$ are
even, while $\xi_j$ are odd. A \textit{superdomain} in $\cC^{n|m}$ is, by
definition, a~$\mathbb Z_2$-graded ringed space of the form $(U,\mathcal F_{n|m})$,
where $U$ is an open subset of $\mathbb C^n$.
\end{Example}

A \textit{complex-analytic supermanifold} of dimension $n|m$ is a
$\mathbb Z_2$-graded ringed space that is locally isomorphic to a
superdomain in $\cC^{n|m}$. Thus, if $\M$ is a~supermanifold, then for
any point $x_0\in M$ there exist a~neighborhood $U$ of $x_0$ in $M$
and an isomorphism of the ringed space $(U,\mathcal O|_U)$ onto a~superdomain
$(\widetilde U,\mathcal F_{n|m})$ in $\cC^{n|m}$ called a~\textit{chart} on
$U$. Let $x_1,\dots ,x_n$ denote the standard coordinates in $\mathbb C^n$.
Identifying $(U,\mathcal O)$ with the superdomain by means of the chart,
we get the elements $x_i$ for $i = 1,\dots ,n$, and $\xi_j$ for $j = 1,\dots ,m$ of
$\mathcal O(U)$ called the \textit{local coordinates} on $U$.

Let $U$ (resp. $V$) be two open subsets of $M$ admitting two charts
with local coordinates for $i = 1,\dots ,n$, and $\xi_j$ for $j = 1,\dots ,m$
(resp. $y_i$ for $i = 1,\ldots,n$ and $\et_j$ for $j = 1,\ldots,m$). Then, in $U\cap V$
we can write
\begin{equation}
\begin{aligned}
y_i &= \ph_i(x_1,\ldots,x_n,\xi_1,\ldots,\xi_m),\text{~~where~~}  i = 1,\ldots,n;\\
\et_j &= \ps_j(x_1,\ldots,x_n,\xi_1,\ldots,\xi_m), \text{~~where~~} j = 1,\ldots,m,
\end{aligned}\label{(1.10)}
\end{equation}
where $\ph_i,\;\ps_j$ are, respectively, even and odd sections of
$\mathcal F_{n|m}$ called the \textit{transition functions}. Similarly, there are
transition functions realizing the inverse coordinate transformation.
An \textit{atlas} of $\M$ is a~cover of $M$ by open subsets that admit certain
charts; any atlas determines a~supermanifold up to isomorphism.

Let $\M$ be a~supermanifold and
\begin{equation}
\mathcal J = (\mathcal O_{\bar 1}) = \mathcal O_{\bar 1} + (\mathcal O_{\bar 1})^2 \label{(1.11)}
\end{equation}
the subsheaf of ideals of $\mathcal O$ generated by the subsheaf
$\mathcal O_{\bar 1}$ of odd elements. Manin (see \cite{20}) denoted 
\[
\text{$\mathcal O_{\rd} := \mathcal O/\mathcal J$.
and $M_{\rd} := (M,\mathcal O_{\rd})$}.
\] 
So, $M_{\rd}$ is a~usual complex analytic
manifold of dimension $n$ called the \textit{odd reduction} of $\M$, and we
have a~morphism 
\[
\red = (\id,p_0): (M,\mathcal O_{\rd})\tto (M,\mathcal O),
\]
 where
$p_0: \mathcal O\tto\mathcal O_{\rd}$ is the canonical projection. This morphism
takes the odd local coordinates $\xi_j$ to 0
and the even ones $x_i$ to certain local coordinates $X_1,\dots ,X_n$ on
$M_{\rd}$. Clearly, any chart of $\M$ determines a~chart on $M$, and, on the
intersection of two charts, the transition functions transforming $X_i =
p_0(x_i)$ into $Y_i = p_0(y_i)$ (see eq.~ \eqref{(1.10)}) have the form
$$
Y_i = \ph_i(X_1,\ldots,X_n,0,\ldots,0).
$$

Though one should distinguish between the coordinates $x_i$ of $\M$ and the
coordinates $X_i$ on $M$, they often are denoted in the same way. In what
follows, we
usually denote the sheaf $\mathcal O_{\rd}$ by $\mathcal F$. The complex manifold
$M_{\rd} = (M,\mathcal F)$ will usually be denoted just by $M$.

Any morphism of supermanifolds 
\[
F = (f,f^*): (M,\mathcal O_M)\tto (N,\mathcal O_N)
\]
induces a~morphism
of manifolds $M\to N$. This just means that the mapping ${f: M\to N}$
is holomorphic. As a~consequence, we get a~canonical homomorphism of groups from
$\Aut\M$ to $\operatorname{Bih}\,M$, the group of all biholomorphic transformations of $M$.

Any superdomain is, clearly, a~supermanifold. More complicated examples will
be given below.

\ssbegin[Supermanifold $(M,\Om)$]{Example}\label{E1.9}
Let $M$ be a~complex manifold of dimension $n$ and $\Om = \bigoplus_{p=0}^n
\Om^p$ be the sheaf of holomorphic exterior forms on $M$. Then, $(M,\Om)$ is
a supermanifold of dimension $n|n$. Indeed, let $U$ be an open subset of
$M$, where a~chart with local coordinates $\row x1n$ is defined. Clearly,
the sheaf $\Om|_U$ can be identified with $\bigwedge_{\mathcal F_n}(\row{dx}1m)$.
Denoting $\xi_j := dx_j$, we see that the $x_i,\,\xi_j$ are local coordinates for
$(M,\Om)$. If $V$ is another open subset with local coordinates $y_i$ and
$\et_j := dy_j$, then the transition functions in $U\cap V$ have the form
$$
\begin{aligned}
y_i &= \ph_i(\row y1n),\ \ i = 1,\ldots,n,\\
\et_j &= \sum_{k=1}^n\nfrac{\partial y_j}{\partial x_k}\xi_k,\ \ j = 1,\ldots,
n,
\end{aligned}
$$
where $\ph_i$ are the usual transition functions for $M$.
\end{Example}

The simplest class of supermanifolds are the so-called split ones.
Let $(M,\mathcal F)$ be a~complex manifold and $\mathcal E$ a~locally free
analytic sheaf on it. Defining $\mathcal O = \bigwedge_{\mathcal F}\mathcal E$, we get
a supermanifold $\M$. A supermanifold is called \textit{split}
if it is isomorphic to a~supermanifold of this form.

The structure
sheaf ${\mathcal O}$ of a~split supermanifold admits the ${\mathbb Z}$-grading
${\mathcal O} = \bigoplus_{p\ge 0} {\mathcal O}_{p}$, where
$$
{\mathcal O}_{p}\simeq \bigwedge^{p}_{\mathcal F} {\mathcal E};
$$
this ${\mathbb Z}$-grading on it is compatible with the ${\mathbb Z_{2}}$-grading.
In what
follows, we often omit the subscript ${\mathcal F}$ while denoting the exterior
powers, the tensor products etc. of the sheaves of ${\mathcal F}$-modules.

Let $U$ be a~coordinate neighborhood in $M$, over which the sheaf $\mathcal E$ is
free or, which is the same, the corresponding vector bundle $\mathbf E$
is trivial. Then, we can choose special local coordinates of $\M$ in
$U$; these are $x_i,\;\xi_j$, where $\row x1n$ are local coordinates of
$M$, while $\row\xi 1m$ is a~basis of the free $\mathcal F_U$-module
$\Ga(U,\mathcal E)$. These local coordinates will be called \textit{splitting} ones.
The transition functions between two systems of splitting coordinates (see eq.~
\eqref{(1.10)}) have the following special form:
$$
\begin{aligned}
y_i &= \ph_i(\row x1n),\quad i = 1,\ldots,n,\\
\et_j &= \sum_{j=1}^m\ps_{jk}(\row x1n)\xi_k,\quad j = 1,\ldots,m,
\end{aligned}
$$
where $\ph_i,\;\ps_{jk}$ are holomorphic functions in $x_i$, and the matrices
$(\ps_{jk})$ are the transition functions of the vector bundle $\mathbf E$.

A classical example of a~split supermanifold is $(M,\Om)$ (see Example \ref{E1.9}).
The sheaf $\Om$ corresponds to the cotangent bundle $\mathbf E = \mathbf T(M)^*$
over $M$. As splitting coordinates one can choose the $x_i$ and $\xi_j := dx_j$,
where $x_i$ are local holomorphic coordinates on $M$.

Another important example is that of the \textit{complex projective superspace}.
                                                     
\ssbegin[Projective superspace $\mathcal{CP}^{n|m}$]{Example}\label{E1.10}
Formally, a~``point" of the projective superspace $\mathcal{CP}^{n|m}$ is
determined by a~row of ``homogeneous coordinates"
\[
    (z_0 \colon \dotsc \colon z_n \colon \ze_1 \colon \dotsc \colon \ze_m),
\]
where $p(z_i) = \bar 0$ and $p(\ze_j) = \od$ and $(\row z0n)\ne 0$,
which is defined up to multiplication by a~non-zero complex number. As $M$,
we take the usual projective space $\mathcal{CP}^{n}$; its points are given by
the homogeneous coordinates $(z_0:\ldots :z_n)$. As usual, consider the
cover of $M$ by the affine open sets $U_k = \{z_k\ne 0\}$ for any $k = 0,\ldots,n$.
In $U_k$, we can uniquely write the coordinate row $(z,\ze)$ as
$$
(x^{(k)}_1,\ldots,x^{(k)}_k,1,x^{(k)}_{k+1},\ldots,x^{(k)}_n,
\xi^{(k)}_1,\ldots,\xi^{(k)}_m),
$$
where $x^{(k)}_i,\;\xi^{(k)}_j$ are, by definition, the local coordinates of
the supermanifold
$\mathcal{CP}^{n|m}$ in $U_k$, expressed through homogeneous coordinates by
$$
\begin{aligned}
x^{(k)}_i &= \begin{cases}\nfrac{z_{i-1}}{z_k} &\text{for $1\le i\le k$}\\
\nfrac{z_i}{z_k} &\text{for $k+1\le i\le n$},\end{cases}\\
\xi^{(k)}_j &= \nfrac{\ze_j}{z_k},\; 1\le j\le m.
\end{aligned}
$$
One can easily write down the transition functions, showing that $\mathcal{CP}^{n|m}$
is a~split supermanifold. The sheaf $\mathcal E$ is $\mathcal F(-1)^m$, where
$\mathcal F(-1)$ is the invertible sheaf determined by a~hyperplane of
$\mathbb{CP}^n$.
\end{Example}

\subsection{Subsupermanifold, retract}\label{ssRetract}
Let $(M,\mathcal O_M)$ be a~supermanifold and $\mathcal I$ be a~$\mathbb Z_2$-graded subsheaf of
ideals of $\mathcal O_M$.  Setting
$$
N = \{x\in M \mid \pi(\ph)(x) = 0\ \ \text{~~for all~~}\ph\in(\mathcal O_M)_x\},\ \
\mathcal O_N = (\mathcal O_M/\mathcal I)|_N,
$$
we get the $\mathbb Z_2$-graded ringed space $(N,\mathcal O_N)$. If this space is a
supermanifold, then it is called a~\textit{submanifold} of $(M,\mathcal O_M)$. If
the sheaf of ideals $\mathcal I$ is generated, over an open set $U\subset M$, by
its
homogeneous sections $\ph_1,\ldots,\ph_s$, then it is usual to say that the
submanifold is determined in $U$ by the equations $\ph_i = 0$ for $i = 1,\ldots,
s$.

For example, the subsheaf $\mathcal J$, given by the formula eq.~\eqref{(1.11)}, determines the reduction
$M_{\rd}$ of $(M,\mathcal O_M)$, which is thus a~submanifold of $(M,\mathcal O_M)$.
For other examples, see Subsections~\ref{E1.12}, \ref{E1.13}.

There is a~construction that to any supermanifold
$\M$ assigns a~split one. Consider the filtration
\begin{equation}
{\mathcal O} = {\mathcal J}^{0} \supset {\mathcal J}^{1} \supset {\mathcal J}^{2}
\supset \ldots\label{(1.12)}
\end{equation}
of $\mathcal O$ by the powers of the subsheaf of ideals ${\mathcal J}$ given by the formula
\eqref{(1.11)}. The associated graded sheaf
$$
\gr{\mathcal O}=\bigoplus_{p\ge 0}\gr^p\mathcal O,
$$
where $\gr^{p}{\mathcal O} = {\mathcal J}^{p}/{\mathcal J}^{p+1}$,
gives rise to the split supermanifold $(M,\gr{\mathcal O})$.

Indeed, $\gr{\mathcal O} \simeq \bigwedge_{{\mathcal F}}{\mathcal E}$,
where $\mathcal F = \gr^{0}{\mathcal O} = \mathcal O_{\rd}$ and ${\mathcal E} =
\gr^{1}{\mathcal O}$ is a~locally free sheaf of $\mathcal F$-modules. Clearly, $\M$
and $(M,\gr\mathcal O)$ have the same dimension.
The supermanifold $(M,\gr\mathcal O)$ is called the \textit{retract} of the supermanifold $\M$. 
           
A
supermanifold is split if and only if it is isomorphic to its retract. If the
$x_i$ and $\xi_j$ are arbitrary local coordinates of $\M$ in a~neighborhood $U
\subset M$, then $X_i = x_i +\mathcal J^2$ and $ \Xi_j = \xi_j +\mathcal J$ are splitting
local coordinates of $(M,\gr\mathcal O)$ in $U$, and one gets the transition
functions between these splitting coordinates, if one takes the terms of
degree 0 (respectively 1) in $\xi_j$ in the transition functions $\ph_i$
(respectively $\ps_j$) for $\M$ (see eq.~\eqref{(1.10)}).

Thus, we see that with any supermanifold $\M$ two objects of the classical
complex analytic geometry are associated: the complex manifold $(M,\mathcal F)$
and the holomorphic vector bundle $\mathbf E$ over $(M,\mathcal F)$ corresponding
to the sheaf $\mathcal E$. It turns out that $\M$ is not, in general, determined
by these two objects up to an isomorphism, since there exist non-split
supermanifolds. For examples, see below.

To settle, if a~given supermanifold $\M$ is split, one can consider the
following exact sequences of sheaves over $M$:
\begin{equation}
\begin{aligned}
0\tto\mathcal J&\tto\mathcal O\overset {p_0}\rightarrow\mathcal F\tto 0,\\
0\tto\mathcal J^2&\tto\mathcal J\overset {p_1}\rightarrow\mathcal E\tto 0.
\end{aligned}\label{(1.13)}
\end{equation}
If the supermanifold $\M$ is split, then both these exact sequences
are split, i.e.,  there exist homomorphisms $q_0:\mathcal F\tto\mathcal O$ and $q_1:
\mathcal E\tto\mathcal J$ such that $p_iq_i = \id,\;i = 0,1$. The obstructions
to splitness lie in certain sheaf cohomology of $M$ (see \cite[Ch.4, Sect. 2]{20}, \cite[Ch.4, Sects.  6, 7]{2}, \cite[Ch.3, Sects. 6, 7]{2b}).

\subsection{Super-Grassmannians
}\label{ss1.11-}
In this subsection, we will briefly consider certain examples of complex
supermanifolds
introduced by Yu. Manin in \cite{20}. Actually, four series of compact
complex supermanifolds corresponding to the following four series of
classical complex Lie superalgebras, were constructed:
\begin{enumerate}
\item[$(1)$]
$\mathfrak{gl}_{n|m}(\mathbb C)$ --- the general linear Lie superalgebra of the
vector superspace $\mathbb C^{n|m}$;
\item[$(2)$]
$\mathfrak{osp}_{n|m}(\mathbb C)$ --- the orthosymplectic Lie superalgebra that
annihilates a~non-degene\-rate even symmetric bilinear form in $\mathbb C^{n|m},
\;m$ being even;
\item[$(3)$]
$\mathfrak{pe}_{n|n}(\mathbb C)$ --- the linear Lie superalgebra that
annihilates a~non-degenerate odd symmetric bilinear form in
$\mathbb C^{n|n}$ (Manin denoted $\mathfrak{pe}_{n|n}(\mathbb C)$ by $\pi\fsp_n(\mathbb C)$ in \cite{20}, see also \cite{V4}, but A.~Weil's suggestion to call the odd non-degenerate bilinear form, and the Lie superalgebra/supergroup it preserves, \textit{periplectic} took over and is now universally accepted together with the name \textit{queer} for the following purely super analog of $\fgl_n$, \textit{D.L.});
\item[$(4)$]
$\mathfrak{q}_{n}(\mathbb C)$ --- the linear Lie superalgebra that commutes with an
odd involution in $\mathbb C^{n|n}$.
\end{enumerate}

These supermanifolds are called the \textit{flag supermanifolds} in case
(1), the \textit{supermanifolds of isotropic flags} in cases (2) and (3),
and
the \textit{supermanifolds of $\Pi$-symmetric flags} in case (4). We will
call them the \textit{classical flag supermanifolds}. They are,
in most cases, non-split. 

Here we consider the classical flag supermanifolds
under assumption that
the flags have the minimal possible length; these are so-called {\it
super-Grassmannians}. The super-Grassmannians are basic in Manin's
constructions, because the flag supermanifolds are defined inductively as
relative super-Grassmannians over the flag supermanifolds of lesser length.

As in Example \ref{E1.8}, we denote by $e_1,\ldots,e_n,f_1,\ldots,f_m$ the standard
basis of $\mathbb C^{n|m}$.

\ssbegin[The super-Grassmannian]{Example}\label{E1.12} The super-Grassmannian $\Gr^{n|m}_{k|l}$ of
$(k|l)$-dimensional subspaces in $\cC^{n|m}$ is a~natural generalization
of the projective superspace $\mathcal{CP}^{n|m} = \Gr^{n+1|m}_{1|0}$. Its
structure is determined by the $(n+m)\times (k+l)$ coordinate matrix
$$
Z = \begin{pmatrix} Z_{00} & Z_{01}\\ Z_{10} & Z_{11}\end{pmatrix},
$$
where $Z_{00}$ and $Z_{11}$ are $n\x k$- and $m\x l$-matrices, respectively,
whose entries are even homogeneous coordinates, while $Z_{01}$ and
$Z_{10}$ are $n\x l$- and $m\x k$-matrices, respectively, whose entries are
odd ones. It is supposed that $Z_{00}$ and $Z_{11}$ are complex matrices of
ranks $k$ and $l$, respectively, so that each of them determines a~point of
the complex Grassmannian $\Gr^{n}_{k}$ or $\Gr^{m}_{l}$, respectively. 

Thus, we
get an element $x_0$ of the manifold $M = \Gr^{n}_{k}\times\Gr^{m}_{l}$; this
manifold is
the reduction of the super-Grassmannian. The matrix $Z$ is to be regarded up
to the following equivalence: 
\[
\text{$Z\sim Z'$ if $Z' = ZQ$, where $Q$ is an
invertible $k\x l$-matrix.}
\]
 If we fix an invertible $k\x l$-submatrix of $Z$,
then the remaining entries of $Z$ give us the even and the odd local
coordinates in a~neighborhood of $x_0$. Using the equivalence, we can assume
that the fixed submatrix is the unit matrix $I_{k|l} = \begin{pmatrix} I_k & 0\\
0 & I_l\end{pmatrix}$.

For example, choose 
\[
x_0 = \langle e_{n-k+1},\ldots,e_n,f_1,\ldots,f_l
\rangle = (\langle e_{n-k+1},\ldots,e_n\rangle,\langle f_1,\ldots,f_l
\rangle). 
\]
Then, the coordinate matrix can be written in the form
\begin{equation}
Z = \begin{pmatrix} X & \Xi\\ I_k & 0\\ 0 & I_l\\ \H &
Y\end{pmatrix},\label{(1.14)}
\end{equation}
where 
\[
\begin{array}{l}
X = (x_{ij}),\; Y = (y_{\al p}),\; \Xi = (\xi_{ip}),\; \H =
(\et_{\al i}),\\
i = 1,\ldots,n-k,\; j = 1,\ldots,k,\; p = 1,\ldots,l,\; \al =
l+1,\ldots, n. 
\end{array}
\]
Here $x_{ij}$ and $y_{\al p}$ are even local coordinates
satisfying $x_{ij}(x_0) = y_{\al p}(x_0) = 0$, while $\xi_{ip}$ and
$\et_{\al i}$ are odd ones. In particular, we have
\[
\dim\Gr^{n|m}_{k|l} = n(n-k)+m(m-l)\,\mid\,n(m-l)+m(n-k).
\]

In the case where $0 < k < n$ and $0 < l < m$, the supermanifold
$\Gr^{n|m}_{k|l}$ is non-split. The simplest non-split
super-Grassmannian is $\Gr^{2|2}_{1|1}$ of dimension $2|2$.
\end{Example}

\ssbegin[The isotropic super-Grassmannian. Superquadric]{Example}\label{E1.13} Let an even non-degenerate symmetric bilinear form $b$ be given in $\cC^{n|m}$. Then, it is possible to define the subsupermanifold $\I\Gr^{n|m}_{k|l}$ of $\Gr^{n|m}_{k|l}$, consisting of subspaces that are (totally) isotropic with respect to $b$; this $\I\Gr^{n|m}_{k|l}$ is called {\it isotropic super-Grassmannian}).

If $n$ is odd, then we get the simplest non-split supermanifolds 
for 
\[
\text{$n = 3$, $k = 1$, $m = 2s\ge 2$, $l = 0$. }
\]
The supermanifold
$\cQ^{1,m} = \I\Gr^{3|m}_{1|0}$ is called the
\textit{superquadric} in the projective superplane $\mathcal{CP}^{2|m}$. In
homogeneous coordinates (see Example \ref{E1.9}), we can express the superquadric by
the equation
$$
z_0^2 - z_1z_2 + \sum_{i=1}^s\ze_i\ze_{s+i} = 0.
$$
The local coordinates on the superquadric are given by the formula
$$
\begin{aligned}
x &= \nfrac{z_0}{z_1},\; \xi_j = \nfrac{\ze_j}{z_1}\text{ ~for $z_1\ne 0$};\\
y &= \nfrac{z_0}{z_2},\; \et_j = \nfrac{\ze_j}{z_2}\text{ ~for $z_2\ne 0$},
\end{aligned}
$$
and the transition functions have the form
$$
\begin{aligned}
y &= x^{-1}(1 + x^{-2}\sum_{i=1}^s\xi_i\xi_{s+i})^{-1},\\
\et_j &= x^{-2}(1 + x^{-2}\sum_{i=1}^s\xi_i\xi_{s+i})^{-1}\xi_j,\text{~~where~~} j =
1,\ldots,2s.
\end{aligned}
$$

Historically, this was (for $m = 2$) one of the first examples of non-split
supermanifolds (see \cite{9}, \cite{2}, \cite{20}).
\end{Example}

\ssbegin[The odd isotropic super-Grassmannian]{Example}\label{E1.14}
Quite similarly, an isot\-ro\-pic super-Grassmannian, associated with an odd non-degenerate anti-symmetric (or symmetric) bilinear form $b$ is defined. In this case, $n = m$, and we denote by $\I_{\odd}\Gr^{n|n}_{k|l}$ the corresponding submanifold of $\Gr^{n|n}_{k|l}$.
\end{Example}

\ssbegin[The $\Pi$-symmetric super-Grassmannian]{Example}\label{E1.15}
Suppose that $m = n$ and that an odd involutive linear transformation
$\Pi$ of the vector superspace $\mathbb C^{n|n}$ is given. Then, we can define
the submanifold $\Pi\Gr^{n|n}_{s|s}$ of $\Gr^{n|n}_{s|s}$ that
consists of $\Pi$-invariant subspaces of dimension $s|s$ (the $\Pi$-{\it
symmetric super-Grassmannian}). This super-Grassmannian is one of the main
objects of our study, and therefore it will be considered in more details in
Section~5. We only mention here that the retract of $\Pi\Gr^{n|n}_{s|s}$ is the
supermanifold $(\Gr^{n}_{s},\Om)$ of Example \ref{E1.9}.
\end{Example}

\section{Tangent sheaf and vector fields}\label{S2}

\subsection{Tangent space and tangent sheaf}\label{ss2.1}
We retain the notation of Subsection~\ref{ss1.7}. Let $\M$ be a~complex
supermanifold. Fix a~point $x\in M$. Using local coordinates $x_1, \dotsc, x_n, \xi_1, \dotsc, \xi_m$ in a~neighborhood of $x$, we can identify the superalgebra $\mathcal O_x$ with $\mathbb C\{x_1,\ldots,x_n\} \otimes\bigwedge_{\mathbb C}(\xi_1, \ldots, \xi_m)$ for any $x \in M$. Notice that this is a~local superalgebra whose unique maximal ideal is $\fm_x := (x_1,\ldots,x_n,\xi_1,\ldots,\xi_m)$. The vector superspace $T_x\M = (\fm_x/\fm_x^2)^*$ is called the {\it tangent space} to $\M$ at the point $x$. Since $\mathcal F = \mathcal O/\mathcal J$, we have the exact sequence
$$
0\tto\mathcal J_x\tto \fm_x\tto n_x\tto 0,
$$
where $\fn_x$ is the maximal ideal of $\mathcal F_x$. This implies the
following exact sequence:
$$
0\tto\mathcal J_x/\fm_x\mathcal J_x\tto \fm_x/(\fm_x)^2\to
\mathfrak \fn_x/(\mathfrak \fn_x)^2\tto 0.
$$
The fiber at $x$ of the vector bundle $\mathbf E$ corresponding to $\M$ is
$\mathcal J_x/\fm_x\mathcal J_x = E_x$. Since $T_x(M) =
(\fn_x/\fn_x^2)^*$ is the tangent vector space to $M$ at $x$, we get
the exact sequence
$$
0\tto T_x(M)\tto T_x\M\tto E_x^*\tto 0.
$$
This gives the canonical identifications
$$
T_x\M_{\bar 0} = T_x(M),\ \ \  T_x\M_{\bar 1} = E_x^*.
$$

The \textit{tangent sheaf} of a~supermanifold $\M$ is by definition the sheaf
$\mathcal T = \mathcal Der\,\mathcal O$ of derivations of the structure sheaf $\mathcal O$.
Its stalk at $x \in M$ is the Lie superalgebra $\mathfrak{der}_{\mathbb C}\mathcal O_x$
of derivations of the superalgebra $\mathcal O_x$. Its sections
are caled \textit{holomorphic vector fields} on $\M$. The vector superspace
$\mathfrak v\M = \Ga(M,\mathcal T)$ of all holomorphic vector fields is
finite-dimensional whenever $M$ is compact. We regard it as a~complex Lie
superalgebra with the bracket
\begin{equation}
[X,Y] = XY - (-1)^{p(X)p(Y)}YX.\label{(2.1)}
\end{equation}

Fix a~point $x \in M$. Any $\delta \in \mathfrak {der}\,\mathcal O_x$ is such that
$\delta (\mathfrak m^2_x) \subset \mathfrak m_x$, and hence defines a~linear mapping
$\widetilde\de: \mathfrak m_x/\mathfrak m^2_x \tto \mathcal O_x/\mathfrak m_x
= \mathbb C$ which is an element of $T_x\M$. This permits us
to define an even linear mapping $\ev_x: \mathfrak v\M\tto T_x\M$ by
$$
\ev_x(v) = \widetilde v_x.
$$
We note that, in contrast with the non-super case, a~vector field $v$
is not, in general, uniquely determined by its values $\widetilde v_x$
at all $x \in M $.

Endow the tangent sheaf $\mathcal T$ with the following filtration:
\begin{equation}
\mathcal T = \mathcal T_{(-1)} \supset \mathcal T_{(0)} \supset \ldots
\supset \mathcal T_{(m)} \supset \mathcal T_{(m+1)} = 0,\label{(2.2)}
\end{equation}
where
$$
\mathcal T_{(p)} = \{\delta \in \mathcal T \mid \delta (\mathcal O)
  \subset \mathcal J^p, \delta (\mathcal J)
  \subset \mathcal J^{p+1}\}\text{ for any }p \ge 0.
$$
Thus, we have obtained a~filtered sheaf of Lie superalgebras. Let $\gr\mathcal T$ denote
the corresponding graded sheaf of algebras. Any $v\in\mathcal T_{(p)}$ maps
$\mathcal J^q$ to $\mathcal J^{q+p}$, inducing a~derivation from
$(\mathcal T_{\gr})_p$, where $\mathcal T_{\gr} = \mathcal Der\,\gr\mathcal O$. As a~result,
we get a~homomorphism $\si_p: \mathcal T_{(p)}\to(\mathcal T_{\gr})_p$. It is easy to check
(see \cite{21}) that the following assertion is true.

\ssbegin[An exact sequence]{Proposition}\label{P2.1}
The following sequences of sheaves are exact:
$$
0\tto\mathcal T_{(p+1)}\tto\mathcal T_{(p)}\overset{\si_p}\to(\mathcal T_{\gr})_p\tto 0,\text{~~where~~}
p\ge -1.
$$

The homomorphisms $\si_p$, where $p\ge -1$, determine an isomorphism of the
sheaves of graded algebras $\gr\mathcal T\tto\mathcal T_{\gr} = \mathcal Der\,\gr\mathcal O$.

\end{Proposition}

In what follows, we will use the cohomology groups $H^p(M,\mathcal T)$ with
values in the tangent sheaf.  Recall that they are finite-dimensional vector
spaces if $M$ is compact. We have $H^0(M,\mathcal T) = \mathfrak v\M$. Since $\mathcal T$
is a~sheaf of Lie superalgebras, we can define the corresponding operation
in
\[
H^{\bcdot}(M,\mathcal T) = \bigoplus_{p\ge0} H^p(M,\mathcal T), 
\]
giving a~graded algebra.
This operation will be denoted by $[ - , -]$; it coincides on $H^0(M,\mathcal T)$
with the bracket defined above. The filtration \eqref{(2.2)} gives rise to a~natural
filtration in $H^{\bcdot}(M,\mathcal T)$, so we get a~filtered algebra.

\subsection{The tangent sheaf of the split supermanifold
}\label{ss2.3}
Here we make some remarks concerning vector fields on
split supermanifolds. If $\M$ is split, then $\mathcal T$ is a
$\mathbb Z$-graded sheaf of Lie superalgebras, the grading being given by the formula
$$
\mathcal T = \bigoplus_{p\ge -1}\mathcal T_p,
$$
where
\begin{equation}
\mathcal T_p := \mathcal Der_p\mathcal O =
  \{\delta \in\mathcal T\mid \delta (\mathcal O_q)
  \subset\mathcal O_{q+p}\text{ for all }q \in\mathbb Z\}.\label{(2.3)}
\end{equation}
Hence,  $\mathfrak v\M :=\bigoplus_{p\ge-1}\mathfrak v\M_p$ is a~$\mathbb Z$-graded Lie
superalgebra. Moreover, we get a~grading in any cohomology 
$H^p(M,\mathcal T)$, turning $H^{\bcdot}(M,\mathcal T)$ into a~bigraded algebra.
One easily verifies that the filtration \eqref{(2.2)} of $\mathcal T$ coincides with the filtration
associated with the grading \eqref{(2.3)}, so that
$$
\mathcal T_{(p)} = \bigoplus_{r\ge p}\mathcal T_r.
$$

Since $\mathcal O = \bigwedge\mathcal E$, where ${\mathcal E}$ is a~locally free
analytic sheaf on $M = (M,\mathcal F)$, it follows that $\mathcal T$ can be regarded as an analytic
sheaf on the complex manifold $M$. It was useful to interpret
$\mathcal T$ directly in terms of the sheaf ${\mathcal E}$. A partial description of
$\mathcal T_p$ for $p \ge -1$ is given by the following exact
sequence of locally free analytic sheaves on $M$ (see \cite{21}):
\begin{equation}
0 \tto\mathcal E^* \otimes \bigwedge ^{p+1}\mathcal E
\overset i\rightarrow \mathcal T_p \overset\al\rightarrow \Th
\otimes\bigwedge^p\mathcal E\tto 0, \label{(2.4)}
\end{equation}
where $\Th = \mathcal Der\mathcal F$ is the tangent sheaf of the manifold $M$.
The mapping $\al$ is the restriction of the derivation of degree $p$ onto
the subsheaf $\mathcal F$, while $i$ identifies any sheaf
homomorphism $\mathcal E\tto\bigwedge^{p+1}\mathcal E$ with a~derivation of degree
$p$ that vanishes on $\mathcal F$. Clearly, $\Im i$ is the subsheaf of $\mathcal T$
consisting of $\mathcal F$-derivations; they act on the stalks of $\mathcal O$ as
derivations of a~Grassmann algebra (see Example~ \ref{E1.5}).

In particular, in the case $p = -1$, we have an isomorphism
\begin{equation}
\mathcal T_{-1} \simeq \mathcal Hom_{\mathcal F}(\mathcal E,\mathcal F) = \mathcal E^*,\label{(2.5)}
\end{equation}
and in the case $p = 0$, we have the exact sequence
\begin{equation}
0\tto\mathcal E^*\otimes\mathcal E\overset i \rightarrow \mathcal T_0
\overset\al\rightarrow \Th\tto 0.\label{(2.6)}
\end{equation}
Let $\mathbf E$ be the holomorphic vector bundle
over $M$ corresponding to the locally free sheaf $\mathcal E$. Clearly, $\mathcal T_0$ is the sheaf of infinitesimal automorphisms of
$\mathbf E$ and $\mathcal E^*\otimes\mathcal E = \mathcal End\;\mathbf E$ is its subsheaf
consisting of germs of endomorphisms preserving each fiber.

The first terms of the cohomology exact sequence, corresponding to the sequence \eqref{(2.6)},
have the form
\begin{equation}
0\tto\mathfrak{gl}(\mathbf E)\overset i \rightarrow \mathfrak v\M_0
\overset\al\rightarrow \mathfrak v(M),\label{(2.7)}
\end{equation}
where $\mathfrak{gl}(\mathbf E) = \Ga (M,\mathcal E^*\otimes\mathcal E)$ is the Lie algebra
of all endomorphisms of $\mathbf E$ (preserving the fibers) and $\mathfrak v(M) =
\Ga (M,\Th)$ is the Lie algebra of holomorphic vector fields on $M$, whereas $i$
and $\al$ are Lie algebra homomorphisms. If $M$ is compact, then we
have the corresponding exact sequence of complex Lie groups
\begin{equation}
e\tto\GL(\mathbf E)\tto\Aut\mathbf E\tto\Bih M,\label{(2.8)}
\end{equation}
where $\Aut\mathbf E$ is the group of automorphisms of $\mathbf E$ and
$\GL(\mathbf E)$ its normal subgroup that consists of automorphisms
preserving the fibers (the \textit{gauge group} of $\mathbf E$) (see \cite{21}).
Note that $\GL(\mathbf E)$ is the group of invertible elements of
$\Mat(\mathbf E)$ regarded as an associative algebra.

We remark that the Lie algebra $\mathfrak{gl}(\mathbf E)$ is never zero (whenever
$m > 0$). Indeed, it always contains the identity endomorphism $\ep$.
Regarded as a~vector field, $\ep$ coincides with the grading derivation of
the sheaf $\mathcal O$ acting by means of eq.~\eqref{(1.1)}. Clearly, in
any splitting local coordinates $x_i,\xi_j$, it has the form
\begin{equation}
\ep = \sum_{j=1}^m\xi_j\pd{\xi_j}.\label{(2.9)}
\end{equation}

\subsection{The tangent sheaf of $(M,\Om)$. The splitting mapping $l$. The Fr\"olicher and Nijenhuis (FN) bracket}\label{ss2.4}
Clearly,  locally, the exact sequence \eqref{(2.4)} splits. But in the case where 
$\mathbf E = \mathbf T(M)^*$ (see Example \ref{E1.9}), there is a~canonical global
splitting, discovered
by Fr\"olicher and Nijenhuis (see \cite{7} and \cite{16}; see also \cite{Gz}). In this case, we
have $\mathcal E = \Om^1$, and hence the sheaves $\bigwedge\mathcal E\otimes\Th$ and
$\bigwedge\mathcal E\otimes\mathcal E^*$ both coincide
with the sheaf $\Om\otimes\Th$ of holomorphic vector-valued differential
forms. Thus, the exact sequence \eqref{(2.4)} has the form
\begin{equation}
0\tto\Om^{p+1}\otimes\Th\overset i \rightarrow\mathcal T\overset\al
\rightarrow\Om^p\otimes\Th\tto 0.\label{(2.10)}
\end{equation}
The splitting mapping $l:\Om\otimes\Th\tto\mathcal T$ is defined by the formula
\begin{equation}
l(\ph) = [i(\ph),d],\label{(2.11)}
\end{equation}
where $d$ is the exterior derivative which, clearly, is a~section of
$\mathcal T_1$. One verifies that $\al(l(\ph)) = \ph$, so that $l$ is, in fact, a
splitting
of the sequence \eqref{(2.10)}. It follows that there is the following decomposition into the
direct sum of sheaves of vector spaces:
\begin{equation}
\mathcal T = i(\Om\otimes\Th)\oplus l(\Om\otimes\Th).\label{(2.12)}
\end{equation}
More precisely,
\begin{equation}
\mathcal T_p = i(\Om_{p+1}\otimes\Th)\oplus l(\Om_p\otimes\Th)\simeq
(\Om_{p+1}\otimes\Th)\oplus(\Om_p\otimes\Th).\label{(2.13)}
\end{equation}

Note that for $p = 0$ the derivation $l(u),\; u\in\Th$ is the classical Lie
derivative along the vector field $u$.

We recall now the Lie bracket in $\mathcal T$. Clearly,
$i(\Om\otimes\Th)$ is a~sheaf of subalgebras of $\mathcal T$, and hence we get
a bracket $\{ - , -\}$ in the sheaf $\Om\otimes\Th$, which is often
called the \textit{algebraic} bracket and is given by the formula \eqref{(1.4)}.
In \cite{7}, another bracket $[-,-]$ was defined in $\Om\otimes\Th$, namely,
$$
[\ph,\psi] = \al([l(\ph),l(\psi)]).
$$
It is called the \textit{FN-bracket}. Under this bracket and the grading
$$
(\Om\otimes\Th)_p = \Om_p\otimes\Th,
$$
the sheaf $\Om\otimes\Th$ is a~sheaf of graded Lie superalgebras as well. We
also have
\begin{equation}
\begin{aligned}
{}[l(\ph),l(\psi)] &= l([\ph,\psi]),\\
{}[i(\ph),l(\psi)] &= l(\ph\barwedge\psi) + (-1)^qi([\ph,\psi]),\ \
\ph\in\Om\otimes\Th,\; \ps\in\Om^q\otimes\Th.
\end{aligned}\label{(2.14)}
\end{equation}
Thus, $l$ is a~homomorphism
of sheaves of graded Lie superalgebras, and the formula \eqref{(2.12)} describes a~decomposition into
the sum of sheaves of subalgebras (but not ideals!).

In particular, from the isomorphism \eqref{(2.5)} we get the isomorphisms
\begin{equation}
\begin{aligned}
&i: \Th\tto\mathcal T_{-1},\\
&i: \mathfrak v(M)\tto\mathfrak v(M,\Om)_{-1},
\end{aligned}\label{(2.15)}
\end{equation}
and from the sequence \eqref{(2.7)} we get the semi-direct decomposition of Lie algebras
\begin{equation}
\mathfrak v(M,\Om)_0 = i(\mathfrak{gl}(\mathbf T(M)^*)\; +\!\!\!\!\!\!\supset
l(\mathfrak v(M)).\label{(2.16)}
\end{equation}
For a~compact $M$, we have the following global form of this decomposition:
\begin{equation}
\Aut\mathbf T(M)^* = \GL(\mathbf T(M)^*)\rtimes\Bih M,\label{(2.17)}
\end{equation}
where the group $\Bih M$ of biholomorphic
transformations of $M$ acts on differential forms in an obvious way.

Note that the identity
endomorphism $\id\in\mathfrak{gl}(\mathbf T(M)^*)$ gives rise to the vector fields
\[
\text{$d = l(\id)\in\mathfrak v(M,\Om)_1$ and $\ep = i(\id)\in\mathfrak v(M,\Om)_0$,}
\]
 the
first one being the exterior derivative and the second one the grading
derivation.

Let $\ps\in\Ga(M,\Om^p)$ be a~holomorphic $p$-form on $M$. Using formula \eqref{(1.5)} we
can define a~vector-valued form $j(\ps)\in\Ga(M,\Om^{p+1}\otimes\Th)$.

\ssbegin[What is $l\circ j$]{Proposition}\label{P2.2}
We have
$$
l(j(\ps)) = \ps d + (-1)^{p+1}(d\ps)\ep\text{~~ for any~~} \ps\in\Ga(M,\Om^p).
$$
\end{Proposition}

\begin{Proof}
Follows immediately from (1.6) and (2.11).
\end{Proof}

\hspace*{-5.6pt}We are now going to discuss the problem of calculating the cohomology algebra $H^{\bcdot}(M,\mathcal T)$ endowed with the bracket induced by the Lie bracket in $\mathcal T$ (for details, see \cite{25}). First, it follows from the decomposition \eqref{(2.12)} that
$$
H^{\bcdot}(M,\mathcal T) = i^*(H^{\bcdot}(M,\Om\otimes\Th))\oplus l^*(H^{\bcdot}(M,\Om\otimes\Th)).
$$
To calculate $H^{\bcdot}(M,\Om\otimes\Th)$, one can use the standard
Dolbeault--Serre resolution consisting of smooth vector-valued forms on $M$
(actually, it was first considered in \cite{8}). Denote $\Ph := \bigoplus_{p,q\ge 0}
\Ph^{p,q}$, where $\Ph^{p,q}$ is the sheaf of complex-valued smooth
$(p,q)$-forms on $M$. Then,  for any
$p\ge 0$, the differential graded sheaf $(\Ph^{p,*}\otimes\Th,\bpd)$ is a
fine resolution of $\Om^p\otimes\Th$, whence
\begin{equation}
H^q(M,\Om^p\otimes\Th)\simeq H^q(\Ga(M,\Ph^{p,*}\otimes\Th),\bpd).\label{(2.18)}
\end{equation}
The algebraic bracket and the FN-bracket in $\Om\otimes\Th$ induce certain
brackets in the graded vector space $H^{\bcdot}(M,\Om\otimes\Th)$. By the
isomorphism \eqref{(2.18)}, they
correspond, respectively, to the algebraic bracket and the FN-bracket in
$\Ph\otimes\Th$ defined in \cite{7}. As to the cohomology of $\mathcal T$, we obtain
the following result.

\ssbegin[Decomposing $H^{\bcdot}(M,\mathcal T)$ using $i$ and $l$]{Proposition}\label{P2.3}
We have
$$
\begin{aligned}
H^{\bcdot}(M,\mathcal T) &= i^*(H^{\bcdot}(M,\Om\otimes\Th))\oplus l^*(H^{\bcdot}(M,\Om\otimes\Th))\\
&\simeq H(\Ga(M,\Ph\otimes\Th),\dif)\oplus H(\Ga(M,\Ph\otimes\Th),\dif).
\end{aligned}
$$
The bigrading in $H^{\bcdot}(M,\mathcal T)$ is given by the formula
$$
H^q(M,\mathcal T_p)\simeq H^q(\Ga(M,\Ph^{p+1,\bcdot}\otimes\Th),\dif)\oplus
H^q(\Ga(M,\Ph^{p,\bcdot}\otimes\Th),\dif),\ \ p\ge -1,\; q\ge 0,
$$
and the bracket $[\al,\be]$, where $\al,\be\in H^{\bcdot}(M,\mathcal T)$, is determined by the
algebraic bracket of smooth vector-valued forms in the left summand, by the
FN-bracket in the right one and by the formula \eqref{(2.14)} in the case where  $\al,\be$ belong
to different summands.

\end{Proposition}

\subsection{Actions of Lie superalgebras on supermanifolds. Transitive and $\bar 0$-tran\-si\-tive supermanifolds} \label{ss2.7}
Let $\M$ be a~supermanifold and $\mathfrak g$ a~Lie superalgebra. An \textit{action} of $\mathfrak g$ on $\M$ is an arbitrary Lie superalgebra homomorphism $\ph: \mathfrak g \tto\mathfrak v\M$. Actions of Lie superalgebras usually appear as the differentials of actions of Lie supergroups, but we will avoid to consider the general (rather technical) Lie theory for supermanifolds, referring to \cite{18}. Actually, we will deal only with the standard actions of classical Lie supergroups on super-Grassmannians (see Section \ref{S5}).

If an action $\ph: \mathfrak g \tto\mathfrak v\M$ is
given, then with any $x \in M$ the linear mapping 
\[
\ph^{x} = \ev_{x}\ph:
\mathfrak g \tto T_{x}\M
\]
is associated. The set $\mathfrak g_{x} = \operatorname
{Ker}\ph^{x}$ is a~subalgebra of $\mathfrak g$, called the \textit{stabilizer} of
$x$. The action $\ph$ is called \textit{transitive} if $\ph^{x}$
is surjective for any $x \in M$. In this case one also says that $\M$ is a
\textit{homogeneous space} of the Lie superalgebra $\mathfrak g$.

Restricting an action $\ph: \mathfrak g \tto\mathfrak v\M$ to the even component, we
get a~homomorphism $\ph _0: \mathfrak g_{\bar 0}\tto \mathfrak v\M_{\bar 0}$. If $M$ is
compact, then, as in the classical Lie theory,
it is possible to integrate $\ph_0$, getting a~homomorphism $\Ph: G\to
\Aut\M$, where $G$ is the simply connected complex Lie group with tangent
algebra $\mathfrak g_{\bar 0}$. This homomorphism induces an action $\Ph_0: G\to
\Bih\, M$ of $G$ on $M$. We say that the action $\ph$ is $\bar 0$-{\it
transitive} if $\Ph_0$ is a
transitive action in the usual sense. Clearly, this is equivalent to the
following condition: $\ph^x: \mathfrak g_{\bar 0}\tto T_x(M)$ is surjective for
any $x\in M$. Any transitive action is $\bar 0$-transitive.

Let again $\ph$ be an action of $\mathfrak g$ on $\M$, where $M$ is compact. Then, 
$G$ acts on the sheaf $\mathcal T$ by the automorphisms 
\[
\text{$g_*: v\mapsto
(\Ph(g)^{-1})^*v\Ph(g)^*$ for any $g\in G$.}
\]
 One immediately verifies 
that
$$
\ev_{gx}g_* = d_x\Ph_0(g)\ev_x,\ \ g\in G,\; x\in M.
$$
It follows that
$$
\ph^{gx}\operatorname{Ad}_g = d_x\Ph_0(g)\ph^x,\ \ g\in G,\; x\in M.
$$
As a~corollary, we get the following proposition.

\ssbegin[Transitive and $\bar 0$-transitive actions]{Proposition}\label{P2.4}
Let $\ph$ be a~$\bar 0$-transitive action of $\mathfrak g$ on $\M$, where $M$ is
compact. Then, 
\begin{enumerate}
\item[$(1)$]
The stabilizers $\mathfrak g_x$, where $x\in M$, of $\ph$ are conjugate by inner
automorphisms of $\mathfrak g$ \textup{(i.e.,  by the automorphisms of the form
$\operatorname{Ad}_g$ for $g\in G$).}
\item[$(2)$]
The action $\ph$ is transitive if and only if the mapping $\ph^x_0: \mathfrak g
_{\bar 1}\tto T_{x_0}\M_{\bar 1}$ is surjective for a~certain $x_0\in M$.
\end{enumerate}
\end{Proposition}

\subsection{Homogeneous and $\bar 0$-homogeneous supermanifolds}\label{ss2.5}Let now $\M$ be a~supermanifold, where $M$ is compact. Then, there is a~
natural action $\ph = \id$ of the finite-dimensional Lie superalgebra
$\mathfrak v\M$ on $\M$. The supermanifold $\M$ is called \textit{homogeneous}
(respectively, $\bar 0$-\textit{homogeneous}) if this action is transitive (respectively,
$\bar 0$-transitive). This means that the mapping $\ev_x:
\mathfrak v\M \tto T_x\M$ (respectively, the even component of this mapping) is
surjective for any $x \in M$. Proposition \ref{P2.4} implies that a
$\bar 0$-homogeneous supermanifold is homogeneous if and only if the odd
component of the mapping $\ev_{x_0}: \mathfrak v\M \tto T_{x_0}\M$ is surjective
for a~certain point $x_0\in M$.

Suppose that an action $\ph$ of a~Lie
superalgebra $\mathfrak g$ on a~supermanifold $(M,\mathcal O)$ is given. We are
going to define an action on the split supermanifold
$(M,\gr\mathcal O)$. To do this, we note that the filtration \eqref{(2.2)}
gives rise to the filtration
$$
\mathfrak g = \mathfrak g_{(-1)} \supset \mathfrak g_{(0)} \supset \ldots\supset
\mathfrak g_{(m)} \supset\mathfrak g_{(m+1)} = 0,
$$
defined by the formula 
$$
\mathfrak g_{(p)} = \mathfrak g \cap \ph\i(\Ga (M,\mathcal T_{(p)})) =
\{u\in\mathfrak g \mid \ph (u)(\mathcal O) \subset \mathcal J^p,\, \ph(u)(\mathcal J) \subset
\mathcal J^{p+1}\}.
$$
Clearly, $\mathfrak g$ is a~filtered Lie superalgebra, and $\ph$ determines
a homomorphism $\widetilde\ph$ of the correspondent graded Lie superalgebra
$\widetilde{\mathfrak g}$ into the graded Lie superalgebra
$\mathfrak v(M,\operatorname {gr}\mathcal O)$, i.e.,  an action of $\widetilde{\mathfrak g}$ on
$(M,\operatorname {gr} \mathcal O)$.

In particular,  consider the natural action $\ph = \id$
of $\mathfrak g = \mathfrak v\M$ on a~compact supermanifold $\M$. We see that
$\mathfrak g_{(p)} = \Ga(M,\mathcal T_{(p)})$, and $\widetilde\ph$ is an injective
homomorphism $\widetilde{\mathfrak g} \tto \mathfrak v(M,\gr \mathcal O)$  induced by $\si_p$ (see Proposition \ref{P2.1}).

\ssbegin[Super-Grassmannians]{Example}\label{E2.9}
Consider the super-Grassmannian $\Gr^{n|m}_{k|l}$ defined in Example
\ref{E1.9}. The general linear Lie supergroup $\GL_{n|m}(\mathbb C)$ acts on
$\Gr^{n|m}_{k|l}$ by multiplying the coordinate matrix $Z$ (see formula \eqref{(1.14)})
on the left by a~matrix of $\GL_{n|m}(\mathbb C)$. The differential of this
action is an action of the Lie superalgebra $\mathfrak{gl}_{n|m}(\mathbb C)$ on this
supermanifold. One easily checks that it is transitive.

The above action induces transitive actions of the subsuperalgebras
$\mathfrak{osp}_{n|m}(\mathbb C)$, $\mathfrak{pe}_{n|n}(\mathbb C)$, and $\mathfrak{q}_{n}(\mathbb C)$
of $\mathfrak{gl}_{n|m}(\mathbb C)$ on the subsupermanifolds $\I\Gr^{n|m}_{k|l}$, 
$\I_{\odd}\Gr^{n|n}_{k|l}$, and $\Pi\Gr^{n|n}_{s|s}$ of the general
super-Grassmannian, considered in Examples \ref{E1.13}, \ref{E1.14}, and \ref{E1.15}, respectively
(for $\Pi\Gr^{n|n}_{s|s}$, the transitivity will be proved in Proposition \ref{P5.5}).
\end{Example}

\section{Classification of non-split supermanifolds
}\label{S3}

\subsection{Sheaves of automorphisms and the classification theorem
}\label{ss3.1}
Let $\M$ be a~supermanifold. Consider the sheaf $\mathcal Aut\;
\mathcal O$ of automorphisms of the structure sheaf $\mathcal O$ of
$\M$ (as usual, any automorphism is even and maps each stalk $\mathcal O_x$, where
$x\in M$, onto itself). This is a~sheaf of groups. For any $F =
(f,\ph)\in\Aut\M$,
the mapping $\operatorname{Int}F: a\mapsto\ph a\ph\i$ is an
automorphism of $\mathcal Aut\,\mathcal O$ which gives an action
$\operatorname{Int}$ of the group $\Aut\M$ on $\mathcal Aut\;\mathcal O$ by
automorphisms of this sheaf.

Clearly, any $a\in\mathcal Aut\;\mathcal O$ maps $\mathcal J$ onto itself, and hence
preserves the filtration \eqref{(1.12)} and induces a~germ of an automorphism of
$\gr\mathcal O$. By definition, $a$ induces the identity mapping on $\mathcal F =
\mathcal O/\mathcal J$. Consider the filtration
\begin{equation}
\mathcal Aut\;\mathcal O = \mathcal Aut_{(0)}\mathcal O\supset\mathcal Aut_{(2)}\mathcal O\supset\ldots,
\label{(3.1)}
\end{equation}
where
$$
\mathcal Aut_{(2p)}\mathcal O= \{a\in\mathcal Aut\;\mathcal O\mid a(u) - u\in\mathcal J^{2p}
\text{ for all }\; u\in\mathcal O\}.
$$
One easily  sees that the $\mathcal Aut_{(2p)}\mathcal O$ are subsheaves of normal
subgroups of $\mathcal Aut\;\mathcal O$. They also are invariant under the action
$\operatorname{Int}$ of $\Aut\M$ defined above.

Following \cite{9} and \cite{33}, \cite{30}, we will use the sheaves of automorphisms in
order to describe the family of all supermanifolds, having as their retract a
given split supermanifold $(M,\mathcal O_{\gr})$. The  1-cohomology sets $H^1(M,\mathcal Aut_{(2p)}\mathcal O_{\gr})$ for
$p\ge 1$ play the main role in this
description. We recall (see \cite{11}) that for any sheaf of (not necessarily
abelian) groups $\mathcal G$
on $M$ the 1-cohomology set $H^1(M,\mathcal G)$ is defined. It has  no
natural group structure, but has a~distinguished element $e$, also called
the \textit{unit element}. We will express the cohomology class  by its \v Cech
cocycle in a~sufficiently fine open cover of $M$. The unit element is
determined by the unit \v Cech cocycle.

Let $\mathbf E$ be a~holomorphic vector bundle over a~complex manifold $M$
and $\mathcal E$ be the sheaf of holomorphic sections of $\mathbf E$. Then, we can consider the split supermanifold $(M,\mathcal O_{\gr})$, where $\mathcal O_{\gr} = \bigwedge\mathcal E$. Let $\Aut\mathbf E$ be the group of all automorphisms of the vector bundle $\mathbf E$. Clearly, any automorphism of $\mathbf E$ gives rise to an automorphism of $(M,\mathcal O_{\gr})$, and thus we get a~natural inclusion $\Aut\mathbf E\subset\Aut(M,\mathcal O_{\gr})$. It follows that $\Aut\mathbf E$ acts on the sheaves $\mathcal Aut_{(2p)}\mathcal O_{\gr}$ by the action $\operatorname{Int}$. Hence,  this group acts on each 1-cohomology set $H^1(M,\mathcal Aut_{(2p)}\mathcal O_{\gr})$ for $p\ge 0$, leaving the unit element $e$ fixed.

\ssbegin[The role of $H^1(M,\mathcal Aut_{(2)}\mathcal O_{\gr})$]{Theorem}\label{T3.1}
To any supermanifold $\M$ that has $(M,\mathcal O_{\gr})$ as its retract
there corresponds an element of the set $H^1(M,\mathcal Aut_{(2)}\mathcal O_{\gr})$.
This correspondence gives rise to a~bijection between the isomorphism classes
of supermanifolds, satisfying the above condition, and the orbits of the
group $\Aut\mathbf E$ on $H^1(M,\mathcal Aut_{(2)}\mathcal O_{\gr})$ under the action
$\operatorname{Int}$ of the group $\Aut\mathbf E$. The given split
supermanifold $(M,\mathcal O_{\gr})$ corresponds to the unit element
$e\in H^1(M,\mathcal Aut_{(2)}\mathcal O_{\gr})$.
\end{Theorem}

\ssbegin[On splitness]{Corollary}\label{Cor3.3}
The following conditions are equivalent:
\begin{enumerate}
\item[$(1)$] The supermanifold 
is split,
i.e.,  isomorphic to its retract.
\item[$(2)$]  $H^1(M,\mathcal Aut_{(2)}\mathcal O_{\gr}) = \{e\}$.
\end{enumerate}
\end{Corollary}

Let us describe the correspondence mentioned in Theorem \ref{T3.1}. Let $\M$ be a
supermanifold such that $\gr\mathcal O = \mathcal O_{\gr}$. We can choose an
open cover $\mathfrak U = (U_i)_{i\in I}$ of $M$ such that the exact sequences
\eqref{(1.13)} split over each $U_i$. Then, we get isomorphisms
$\si _i: \mathcal O|_{U_i}\tto\mathcal O_{\gr}|_{U_i}$, where $i\in I$, inducing the identity
isomorphisms of the $\mathbb Z$-graded sheaves. Setting $g_{ij} := \si _i
\si_j\i$, we obtain a~1-cocycle $g = (g_{ij})\in Z^1(\mathfrak U,\mathcal Aut_{(2)}
\mathcal O_{\gr})$. Its cohomology class $\ga\in H^1(M,\mathcal Aut_{(2)}
\mathcal O_{\gr})$ does not depend of the choice of $\si_i$; this is the
class  desired.

The above cover $\mathfrak U$ can be chosen in such a~way that $\mathbf E$ is
trivial over $U_i$ for any $i\in I$. Then,  for any $i\in I$, we have an
isomorphism 
\[
\rh_i: \mathcal O_{\gr}|_{U_i}\tto\bigwedge_{\mathcal F_n(x^{(i)})}
(\row{\xi^{(i)}}1m)|_{U_i}, 
\]
providing $U_i$ with the local coordinate system
$\row{x^{(i)}}1n,\row{\xi^{(i)}}1m$. For any pair $i,\;j$ such
that $U_i\cap U_j\ne\emptyset$ we get the isomorphism 
\[
\ph_{ij} = \rh_i
\rh_j\i : \bigwedge_{\mathcal F_n(x^{(j)})}(\row{\xi^{(j)}}1m)|_{U_j}\to
\bigwedge_{\mathcal F_n(x^{(i)})}(\row{\xi^{(i)}}1m)|_{U_i}
\]
 which is expressed by
the transition functions of $(M,\mathcal O_{\gr})$. One can ask: ``\textsl{how to write
the transition functions of $\M$ in terms of the transition functions of the retract and
 the cocycle}~ $g$?" 
 
 To answer this question, we have to consider the
isomorphisms 
\[
\ps_{ij} = \rh_i\si_i\si_j\i\rh_j\i = \rh_i
g_{ij}\rh_j\i. 
\]
Clearly, $\ps_{ij} = (\rh_i g_{ij}\rh_i\i)
\ph_{ij}$. This means that the transition functions of $\M$ can be
obtained from the transition functions of $(M,\mathcal O_{\gr})$ by applying the automorphism
$g_{ij}$ expressed in terms of the coordinates $x^{(i)},\;\xi^{(i)}$.

\subsection{The exponential mapping and its applications
}\label{ss3.4}
In general,  to explicitly describe the set
$H^1(M,\mathcal Aut_{(2)}\mathcal O_{\gr})$ is a~difficult problem. But we will see below that,
under certain strong conditions, this set coincides with the 1-cohomology
of a~locally free analytic sheaf on $M$. This simple case is
sufficient for further applications.

We will use the linearization method proposed in \cite{30}. Let $\M$ be an
arbitrary supermanifold of dimension $n|m$. As in the
classical Lie theory, there exists a~natural relationship between
automorphisms and derivations of the sheaf $\mathcal O$. From formula \eqref{(2.2)} we get the
filtration
\begin{equation}
\mathcal T_{(2)\bar 0}\supset\mathcal T_{(4)\bar 0}\supset\ldots,\label{(3.2)}
\end{equation}
where
$$
\mathcal T_{(2p)\bar 0} = \mathcal T_{(2p)}\cap\mathcal T_{\bar 0} = \mathcal T_{(2p-1)}
\cap\mathcal T_{\bar 0} = \{\de\in\mathcal T_{\bar 0}\mid\de (\mathcal O)\subset
\mathcal J^{2p}\}.
$$
Then, we have the exponential mapping
$$
\exp: \mathcal T_{(2)\bar 0}\tto\mathcal Aut_{(2)}\mathcal O.
$$
It is expressed by the usual exponential series which is actually a
polynomial, since $v^k = 0$ for any $v\in\mathcal T_{(2)\bar 0}$ and  any 
$k > \left [\nfrac m2\right]$. One
proves that $\exp$ is bijective \cite{4} and maps $\mathcal T_{(2p)\bar 0}$ onto
$\mathcal Aut_{(2p)}\mathcal O$ for any $p = 1,2.\ldots$. Thus, it is an isomorphism of
sheaves of sets (but in general not of groups). We denote $\log :=\exp^{-1}$.

\ssbegin[Necessary conditions of splitness]{Proposition}\label{P3.1}
For any $p\ge 1$, there is the following exact sequence of the sheaves of
groups:
\begin{equation}
0\tto\mathcal Aut_{(2p+2)}\mathcal O\tto \mathcal Aut_{(2p)}\mathcal O \overset {\la_{2p}}
\to(\mathcal T_{\gr})_{2p}\tto 0,\label{(3.3)}
\end{equation}
where $\mathcal T_{\gr} = \mathcal Der\,\mathcal O_{\gr}$ is the $\mathbb Z$-graded tangent
sheaf of $(M,\mathcal O_{\gr})$ and $\la_p$ is the composition of the following
mappings:
$$
\la_{2p} : \mathcal Aut_{(2p)}\mathcal O\overset {\log}\tto \mathcal T_{(2p)\bar 0}
\overset {\pi_p}\tto \mathcal T_{(2p)\bar 0}/\mathcal T_{(2p+2)\bar 0}\overset {h_p}
\to(\mathcal T_{\gr})_{2p},
$$
$\pi_p$ being the canonical projection and $h_p$ the natural isomorphism
implied by Proposition $\ref{P2.1}$. If $\M = (M,\mathcal O_{\gr})$ is split, then
$\la_{2p}$ maps any
germ $\exp u\in\mathcal Aut_{(2p)}\mathcal O$ onto the $(2p)$-component of
$u\in\mathcal T_{(2p)}$.
\end{Proposition}

\begin{Proof}
Consider the mapping $\widetilde\la_{2p}=\pi_p\log$. Using the
Campbell--Hausdorff formula, we get
$$
\begin{aligned}
\widetilde\la_{2p}((\exp u)(\exp v))&=\widetilde\la_{2p}(\exp(u+v+\nfrac{1}{2}[u,v]+
\dots))=\pi_p(u)+\pi_p(v)=\\
&=\widetilde\la_{2p}(\exp u)+\widetilde\la_{2p}(\exp v).
\end{aligned}
$$
Hence,  $\widetilde\la_{2p}$ is a~homomorphism of sheaves of groups. Clearly,
$\Ker\widetilde\la_{2p}=\mathcal Aut_{(2p+2)}\mathcal O$, and we get the exact sequence
of sheaves of groups
\begin{equation}
0\tto\mathcal Aut_{(2p+2)}\mathcal O\tto\mathcal Aut_{(2p)}\mathcal O\overset{\widetilde\la_{2p}}
\tto \mathcal T_{(2p)\bar 0}/\mathcal T_{(2p+2)\bar 0}\tto 0.\label{(3.4)}
\end{equation}
Clearly, $\mathcal T_{(2p+2)\bar 0}=\mathcal T_{(2p+1)\bar 0}$. Using Proposition
\ref{P2.1}, we get
$$
\mathcal T_{(2p)\bar 0}/\mathcal T_{(2p+2)\bar 0}\simeq
\mathcal T_{(2p)\bar 0}/\mathcal T_{(2p+1)\bar 0}\simeq (\mathcal T_{\gr})_{2p}.
$$
Now the sequence \eqref{(3.3)} follows from the sequence \eqref{(3.4)}.
\end{Proof}

\sssbegin{Lemma}\label{L3.1}
For any $p\ge 2$, if $H^1(M,(\mathcal T_{\gr})_{2p})=0$,  then
$H^1(M,\mathcal Aut_{(2p)}\mathcal O)=\{e\}$.
\end{Lemma}

\begin{Proof}
We will use the induction on $p$. Clearly, the claim is true for all  $p$
sufficiently big. We have to prove that if it is true for a~certain
$p\ge 3$, then it is true for $p-1$ as well. The exact sequence \eqref{(3.3)} gives
the cohomology exact sequence (see \cite{11})
$$
H^1(M,\mathcal Aut_{(2p)}\mathcal O)\tto H^1(M,\mathcal Aut_{(2p-2)}\mathcal O)\overset
{\la_{2p-2}^*}\tto H^1(M,(\mathcal T_{\gr})_{2p-2}).
$$
Clearly, $H^1(M,\mathcal Aut_{(2p-2)}\mathcal O)=\{ e\}$ follows from
$$
H^1(M,(\mathcal T_{\gr})_{2p-2})=0,\ \ H^1(M,\mathcal Aut_{(2p)}\mathcal O)=\{ e\}.
$$
\end{Proof}

\ssbegin{Proposition}\label{P3.2}
Suppose that $H^1(M,(\mathcal T_{\gr})_{2p}) = 0$ for any $p\ge 2$. Then, the
mapping 
\[
\la_2^*: H^1(M,\mathcal Aut_{(2)}\mathcal O_{\gr})\to
H^1(M,(\mathcal T_{\gr})_{2})
\]
 is injective.
\end{Proposition}

\begin{Proof}
The sequence \eqref{(3.3)} for the sheaf $\mathcal O=\mathcal O_{\gr}$ gives the cohomology
exact sequence
\begin{equation}
H^1(\mathcal Aut_{(4)}\mathcal O_{\gr})\tto H^1(M,\mathcal Aut_{(2)}\mathcal O_{\gr})\overset
{\la_2^*}\tto H^1(M,(\mathcal T_{\gr})_{2})\label{(3.5)}
\end{equation}
Suppose that $\gamma,\eta\in H^1(M,\mathcal Aut_{(2)}\mathcal O_{\gr})$ and that
$\la_2^*(\gamma )=\la_2^*(\eta )$. Let $\gamma$ be determined by the cocycle
$(g_{ij})$ and $\eta$ by the cocycle $(h_{ij})$ in a~cover $\mathfrak U =
(U_i)_{i\in I}$ of $M$. Then, our assumption implies that 
\[
\text{$\la_2(g_{ij}) =
\la_2(h_{ij})+c_j-c_i$, where $c_i\in ((\mathcal T_{\gr})_2)_{U_i}$.}
\]
 We can
assume that $c_i=\la_2(g_i)$, where $g_i\in
(\mathcal Aut_{(2)}\mathcal O_{\gr})_{U_i}$. Then, $\la_2(g_ig_{ij}g_j^{-1}) =
\la_2(h_{ij})$. Thus, we can suppose from the beginning that
$\la_2(g_{ij})=\la_2(h_{ij})$. 

Consider the cochain $f\in
C^1(M,\mathcal Aut_{(2)}\mathcal O_{\gr})$, given by the formula $f_{ij}=h_{ij}g_{ij}^{-1}$. Then, 
${\la_2(f) = 0}$.

Let $\M$ be the supermanifold corresponding to the cohomology class of
$\gamma$ due to Theorem~ \ref{T3.1}. Then, $g_{ij}=h_i h_j^{-1}$, where $h_i:
\mathcal O|_{U_i}\tto \mathcal O_{\gr}|_{U_i}$ for $i\in I$, are certain isomorphisms
of sheaves of superalgebras inducing the identity mappings on
$\mathcal O_{\gr}|_{U_i}$. The equalities $h_{ij}=f_{ij}g_{ij}=
f_{ij}h_ih_j^{-1}$ and $h_{ij}h_{jk}=h_{ik}$ imply that
$$
f_{ij}h_ih_j^{-1}f_{jk}h_jh_k^{-1}=f_{ik}h_ih_k^{-1}
$$
or
$$
(h_i^{-1}f_{ij}h_i)(h_j^{-1}f_{jk}h_j)=h_i^{-1}f_{ik}h_i.
$$
Clearly, $\la_2(h_i^{-1}f_{ij}h_i)=0$, whence $(h_i^{-1}f_{ij}h_i)
\in Z^1(\mathfrak U, \mathcal Aut_{(4)}\mathcal O)$. Therefore, 
\[
(h_i^{-1}f_{ij}h_i)\in Z^1(M,\mathcal Aut_{(4)}\mathcal O). 
\]
By Lemma \ref{L3.1}, this
latter cocycle is cohomologous to $e$, i.e.,  
\[
\text{$h_i^{-1}f_{ij}h_i=
u_iu_j^{-1}$, where $u_i\in (\mathcal Aut_{(4)}\mathcal O)_{U_i}$.}
\]
 Thus,
$f_{ij}=h_iu_iu_j^{-1}h_i^{-1}$. It follows that
$$
\begin{aligned}
h_{ij}&=f_{ij}g_{ij}=h_iu_iu_j^{-1}h_i^{-1}g_{ij}=
h_iu_iu_j^{-1}h_j^{-1}=(h_iu_ih_i^{-1})(h_ih_j^{-1})(h_ju_j^{-1}h_j^{-1})=\\
&=v_ig_{ij}v_j^{-1},
\end{aligned}
$$
where $v_i=h_iu_ih_i^{-1}\in (\mathcal Aut_{(4)}\mathcal O_{\gr})_{U_i}$. This
implies that $\gamma =\eta$.
\end{Proof}

\ssbegin{Theorem}\label{T3.2}
Suppose that $H^1(M,(\mathcal T_{\gr})_{2p}) = H^2(M,(\mathcal T_{\gr})_{2p}) = 0$
for any $p\ge 2$. Then, the mapping $\la_2^*:
H^1(M,\mathcal Aut_{(2)}\mathcal O_{\gr})\tto H^1(M,(\mathcal T_{\gr})_{2})$ is bijective.
\end{Theorem}

\begin{Proof}
The injectivity follows from Proposition \ref{P3.2}, while the surjectivity is
implied by Theorem~ 3 of \cite{30}.
\end{Proof}

To calculate the quotient of $H^1(M,\mathcal Aut_{(2)}\mathcal O_{\gr})$ by
$\Aut\mathbf E$, the following assertion is useful. For any $c\in
\mathbb C^{\times}$, denote by $A_c$ the automorphism of $\mathbf E$ given by
the multiplication by the scalar $c$.

\ssbegin[Technical]{Lemma}\label{L3.2}
We have
$$
\la_2^*\circ\Int\, A_c = c^2\la_2^*,\text{ where }c\in\mathbb C^{\times}.
$$
\end{Lemma}

\begin{Proof}
Consider the grading vector field $\ep\in\Gamma (M,(\mathcal T_{\gr})_0)$ on
$(M,\mathcal O)$. Clearly, $A_c$ gives rise to the automorphism
$\al_c =\exp (b\ep)$, where $c = \exp b$, of the structure sheaf
$\mathcal O_{\gr}$.
Let $g = (g_{ij})\in Z^1(\mathfrak U,\mathcal Aut_{(2)}\mathcal O_{\gr})$ be a~cocycle of
a cover $\mathfrak U = (U_i)_{i\in I}$ of $M$ and $w_{ij} = \log g_{ij}$. Then, 
$$
\begin{aligned}
\al_c g_{ij}\al_c^{-1}&=(\exp (b\ep))(\exp w_{ij})(\exp (b\ep))^{-1} =
\exp((\Ad\exp (b\ep))(w_{ij})) =\\
&=\exp ((\exp (b\,\ad\ep))(w_{ij}))=
\exp(w_{ij}+b[\ep,w_{ij}]+\nfrac1{2!}b^2[\ep,[\ep,w_{ij}]]+\dots ).
\end{aligned}
$$
Applying $\la_2$ and denoting the 2-component of the vector
field $w_{ij}$ by $w_{ij}^{(2)}$, we get
$$
\la_2(\al_c g_{ij}\al_c^{-1})=w_{ij}^{(2)}+2bw_{ij}^{(2)}+
\nfrac{4b^2}{2!}w_{ij}^{(2)} + \ldots = (\exp 2b)w_{ij}^{(2)} =
c^2\la_2(g_{ij}).
$$
Thus, $\la_2^*((\operatorname{Int}\al_c)\gamma ) = c^2\la_2^*(\gamma)$,
where $\ga$ is the cohomology class of $g$.
\end{Proof}

This yields, in particular, the following simple fact.

\ssbegin[Uniqueness of non-split supermanifold with given retract]{Proposition}\label{P3.3}\ 

Suppose that
\begin{gather*}
H^1(M,(\mathcal T_{\gr})_2) \simeq \mathbb C,\\
H^1(M,(\mathcal T_{\gr})_{2p}) = 0 \text{ for any } p\ge 2,\\
H^1(M,\mathcal Aut_{(2)}\mathcal O_{\gr})\neq\{ e\}.
\end{gather*}
Then, $\la_2^*$ is bijective and $H^1(M, \mathcal Aut_{(2)}\mathcal O_{\gr})
/\Aut\mathbf E$ consists of two elements. Thus, there exists precisely one
non-split supermanifold having $(M,\mathcal O_{\gr})$ as its retract.
\end{Proposition}

\begin{Proof}
By Proposition \ref{P3.2}, $\la_2^*: H^1(M, \mathcal Aut_{(2)}\mathcal O_{\gr})\tto\mathbb C$ is
injective. Therefore,  $\Im\la_2^*$ contains a~non-zero element, and Lemma
\ref{L3.2} implies that $\la_2^*$ is surjective. By the same lemma, the group
$\Aut\mathbf E$ has precisely two orbits on $H^1(M,\mathcal Aut_{(2)}\mathcal O_{\gr})$.
\end{Proof}

\ssbegin{Remark}
The conditions of Theorem \ref{T3.2} are fulfilled, in particular, if $m = 2$ or 3.
The corresponding special cases were proved in \cite[Ch. 4]{20}, and \cite{4},
respectively. In the general case, the class $\la_2^*(\ga)$ is closely
related to the first obstruction to splitting the sequences \eqref{(1.13)} considered
in \cite[Ch. 4]{20}. If $\la_2^*(\ga) = 0$, then $\ga\in\Im H^1(M,\mathcal Aut_{(4)})$,
and we can apply $\la_4^*$, and so on. The resulting obstruction theory is
discussed in \cite{2}, \cite{5}, \cite{20}, \cite{28}. On the other hand, any non-split supermanifold
can be
regarded as a~result of deformation of its retract, and $\la_2^*$ can be interpreted as
the corresponding Kodaira--Spencer mapping (for details, see \cite{5}, \cite{6}).
\end{Remark}

\subsection{A family of non-split supermanifolds with retract $(M,\Om)$
}\label{ss3.10}
Here we consider the case where $\mathcal O_{\gr}$ is the sheaf of holomorphic
forms $\Om$ on $M$. Using closed (1,1)-forms on $M$, we will construct an
abelian subsheaf of the sheaf of groups $\mathcal Aut_{(2)}\Om$. The 1-cohomology
of this subsheaf determines a~family of supermanifolds with retract $(M,\Om)$.
This family is non-trivial whenever $M$ is a~compact K\"ahler manifold of
$\dim M > 1$ and $H^{1,1}(M,\mathbb C)\ne 0$.

Let $\mathcal Z\Om^1$ denote the subsheaf of $\Om^1$ consisting of closed forms.
Consider the following sequence of sheaves and their homomorphisms:
$$
\mathcal Z\Om^1\overset\be\tto\Om^1\overset\nu\tto\mathcal T_2\overset{\exp}\to
\mathcal Aut_{(2)}\Om,
$$
where $\be$ is the identical inclusion and $\nu$ is given by the formula
\begin{equation}
\nu(\ps) = \ps d.\label{(3.6)}
\end{equation}
We claim that the
composition mapping $\mu: \mathcal Z\Om^1\tto\mathcal Aut_{(2)}\Om$ is a~homomorphism
of sheaves of groups. By formula \eqref{(3.6)}, we see that $\mu(\ps) = \exp(\ps d)$.
Clearly, for any $\ps_1,\ps_2\in\mathcal Z\Om^1$, we have $(\ps_1d)(\ps_2d) = 0$.
Therefore, 
$$
\mu(\ps) = \id + \ps d,
$$
and
$$
\mu(\ps_1 +\ps_2) = \mu(\ps_1)\mu(\ps_2).
$$

It follows that we have the cohomology mapping taking 0
to the unit element
$$
\mu^*: H^1(M,\mathcal Z\Om^1)\tto H^1(M,\mathcal Aut_{(2)}\Om).
$$
Consider the homomorphism of sheaves of groups $\la_2:
\mathcal Aut_{(2)}\Om\tto\mathcal T_2$ defined in Proposition \ref{P3.1}.

\ssbegin[Technical]{Proposition}\label{P3.4}
The relation $\la_2\mu = \nu\be$ holds.
Suppose that $\dim M > 1$. If $\mu^*(\ze) = \mu^*(\ze')$ for
some $\ze,\ze'\in H^1(M,\mathcal Z\Om^1)$, then $\be^*(\ze) = \be^*(\ze')$.
\end{Proposition}

\begin{Proof}
Since $\la_2\exp = \id$ on $\mathcal T_2$, we see that $\la_2\mu = \nu\be$. Now,
it follows from Proposition \ref{P2.2} that $\nu\be = lj\be$, where $j: \Om^1\to
\Om^2\otimes\Th = \mathcal Hom (\Om^1,\Om^2)$ is given by the formula
\begin{equation}
j(\ps)(\al) = \ps\al.\label{(3.7)}
\end{equation}
Thus, $\la_2\mu = lj\be$, whence $\la_2^*\mu^* = l^*j^*\be^*$.

Suppose that $\dim M > 1$. We claim
that $l^*$ and $j^*$ are injective. Indeed, both $l$ and $j$ are injections
onto a~direct summand. In the first case this follows from formula \eqref{(2.12)}. Further, if
$\dim M > 1$, then
$j(\Om^1)$ admits a~direct complement in $\Om^2\otimes\Th$, namely, the
kernel of the contraction mapping $c: \Om^2\otimes\Th\tto\Om^1$ (see formula \eqref{(1.9)}).
Thus, $l^*j^*$ is injective, which implies our second assertion.
\end{Proof}

Let $\mathfrak U = (U_i)$ be an open cover of $M$ and let $\ps =(\ps_{ij})\in Z^1(\mathfrak U,\mathcal Z\Om^1)$.
Then, the above construction assigns to $\ps$ the supermanifold given by
the cocycle
\begin{equation}\label{(3.8)}
g = (g_{ij})\in Z^1(\mathfrak U,\mathcal Aut_{(2)}\Om),
\quad \textup{where} \quad
g_{ij} = \id + \ps_{ij}d.
\end{equation}
Suppose that $\dim M > 1$. Due to Theorem \ref{T3.1}, we see from Proposition \ref{P3.4}
that this supermanifold is non-split whenever the cohomology class of $\ps$
in $H^1(M,\Om^1)$ is non-zero.

 Now we pass to an important case, where a~``closed cocycle" $\ps$
appears. Let $\om$ be a~(1,1)-form on $M$ satisfying $d\om = 0$. Then, 
clearly, $\bpd\om = 0$, and by the Dolbeault theorem $\om$ determines a
cohomology class in $H^1(M,\Om^1)$. It turns out that it can be given
by a~closed \v Cech cocycle.

\ssbegin[Technical]{Lemma}\label{L3.2bis}
We have the exact sequence of sheaves: 
\begin{equation}
0 \tto \mathcal Z\Om^1 \tto  \Ph^{1,0}_{\hpd} \stackrel{\bpd}{\tto} \mathcal Z\Ph^{1,1} \tto  0,
\label{(3.9)}
\end{equation}
where $\Ph^{1,0}_{\hpd}\subset\Ph^{1,0}$ is the subsheaf of
$\bpd$-closed $(1,0)$-forms and $\mathcal Z\Ph^{1,1}\subset\Ph^{1,1}$ the subsheaf
of $d$-closed $(1,1)$-forms.
\end{Lemma}

\begin{Proof}
Clearly, we have the following exact sequence: 
$$
0 \tto  \Om^1 \tto  \Ph^{1,0}\stackrel{\bpd}{\tto} \Ph^{1,1}_{\bpd} \tto  0,
$$
where $\Ph^{1,1}_{\bpd}\subset\Ph^{1,1}$ is the subsheaf of $\bpd$-closed
$(1,1)$-forms. By definition, $\mathcal Z\Om^1 = \Om^1\cap\Ph^{1,0}_{\hpd}$.
Therefore,  we only have to prove that $\bpd(\Ph^{1,0}_{\hpd}) =
\mathcal Z\Ph^{1,1}$.

If $\ph\in\Ph^{1,0}$ and $\hpd\ph = 0$, then $d\bpd\ph = \pd\bpd\ph =
-\bpd\hpd\ph = 0$. Conversely, suppose that the form $\bpd\ph$, where
$\ph\in\Ph^{1,0}$, satisfies $d\bpd\ph = 0$. Then, $\bpd\hpd\ph = 0$. Since
$\hpd\ph\in\Ph^{2,0}$, this form is holomorphic and closed. Therefore, 
$\hpd\ph = \hpd\ph_1$, where $\ph_1\in\Om^1$. Hence,  $\ph - \ph_1\in
\Ph^{1,0}_{\hpd}$, and $\bpd\ph = \bpd(\ph -\ph_1)$.
\end{Proof}

\subsection{An isomorphism $D$}\label{ss4.12} Now, consider the cohomology exact sequence, corresponding to \eqref{(3.9)}: 
$$
\Ga(M,\Ph^{1,0}_{\hpd}) \stackrel{\bpd}{\tto} \Ga(M,\mathcal Z\Ph^{1,1}) \stackrel{\de^*}{\tto} 
H^1(M,\mathcal Z\Om^1).
$$
Using $\de^*$, we get the mapping
$$
\mu^*\de^*: \Ga(M,\mathcal Z\Ph^{1,1})\tto H^1(M,\mathcal Aut_{(2)}\Om).
$$
Thus, any $(1,1)$-form $\om$ on $M$ such that $d\om = 0$ determines a
supermanifold with retract $(M,\Om)$. To obtain an expression of the
corresponding cocycle $g$, we have to find a~cocycle $\ps$ determining
$\de^*(\om)$. Consider
an open cover $\mathfrak U = (U_i)$ of $M$ such that $\om = \bpd\ps_i$ in any
$U_i$, where $\ps_i\in\Ph^{1,0}_{\hpd}(U_i)$. By definition of the connecting
homomorphism $\de^*$, the desired cocycle is
$\ps = (\ps_{ij})\in Z^1(\mathfrak U,\mathcal Z\Om^1)$, where $\ps_{ij} = \ps_j -
\ps_i$ in $U_i\cap U_j\ne\emptyset$. Finally, the cocycle $g$ is given by the formula
\eqref{(3.8)}.

Note that any $\om\in\Ga(M,\mathcal Z\Ph^{1,1})$ satisfies the condition $\bpd\om = 0$, and
hence determines an element $[\om]$ of the Dolbeault cohomology group
\[
H^{1,1}(M,\mathbb C) = \Ga(M,\mathcal Z\Ph^{1,1})/\bpd\Ga(M,\Ph^{1,0}). 
\]
Further,
by the Dolbeault theorem, we have an isomorphism 
\[
D: H^{1,1}(M,\mathbb C)\to
H^1(M,\Om^1).
\]

\ssbegin[An isomorphism $D$]{Proposition}\label{P3.5}
We have
$$
D([\om]) = \be^*\de^*\om,\text{~~for any~~} \om\in\Ga(M,\mathcal Z\Ph^{1,1}).
$$
If $\dim M > 1$, then, for any two $\om,\om'\in
\Ga(M,\mathcal Z\Ph^{1,1})$, the equation $\mu^*\de^*\om =
\mu^*\de^*\om'$ implies $[\om] = [\om']$. In particular, the
supermanifold corresponding to $[\om]$ is non-split whenever $[\om]\ne 0$.

Any K\"ahler form $\om$ on a~compact manifold $M$ of dimension $> 1$
determines a~non-split supermanifold with retract $(M,\Om)$.
\end{Proposition}

\begin{Proof}
The usual proof of the Dolbeault theorem (see, e.g., \cite{10}) shows that
$D([\om])$ is the cohomology class of the cocycle $(\ps)$ described above.
Thus, $D([\om]) = \be^*\de^*\om$. If $\mu^*\de^*\om = \mu^*\de^*\om'$, then
$\be^*\de^*\om = \be^*\de^*\om'$ by Proposition \ref{P3.4}, and hence $D([\om]) =
D([\om'])$, and $[\om] = [\om']$.

If $M$ is compact and $\om$ is a~K\"ahler form, then the de Rham class of
$\om$ is non-zero. Since $M$ is K\"ahler, this implies that $[\om]\ne 0$.
\end{Proof}

The situation is much more simple in the case where $M$ is a~compact K\"ahler
manifold.

\ssbegin[$M$ is a~compact K\"ahler manifold]{Theorem}\label{T3.3}
If $M$ is a~compact K\"ahler manifold, then we have a~linear mapping
\[
\hat\de: H^{1,1}(M,\mathbb C)\tto H^1(M,\mathcal Z\Om^1)
\]
such that $\be^*\hat\de =
D$. If $\dim M > 1$, then the mapping 
\[
\mu^*\hat\de: H^{1,1}(M,\mathbb C)\to
H^1(M,\mathcal Aut_{(2)}\Om)
\]
is injective.
\end{Theorem}

\begin{Proof}
Since $M$ is compact K\"ahler, any cohomology class in $H^{1,1}(M,\mathbb C)$
contains a~closed $(1,1)$-form $\om$, e.g., a~harmonic one. We set
\[
\hat\de([\om]) := \de^*(\om). 
\]

To check the correctness of this definition,
consider a~closed form 
\[
\text{$\om' = \om + \bpd\al$, where $\al\in\Ga(M,\Ph^{0,1})$.}
\]
Then,  by $\hpd\bpd$-Lemma (see \cite{10}), $\om' -\om = \bpd\hpd\ph$, where
$\ph$ is a~$C^{\infty}$ function on $M$. If $\om = \bpd\ps_i$ in $U_i$, then
$\om' = \bpd(\ps_i + \hpd\ph)$. Now it is clear that $\de^*(\om) =
\de^*(\om')$. By Proposition \ref{P3.5}, $\be^*\hat\de = D$. The
injectivity of $\mu^*\hat\de$ follows from the same proposition.
\end{Proof}

We return to the sheaf homomorphism $j: \Om^1\tto\Om^2\otimes\Th$ defined by the formula 
\eqref{(3.7)}. Due to the formula 
\eqref{(1.8)}, $j$ can be also expressed by
\begin{equation}
j(\ps)(u_1,u_2) = \ps(u_1)u_2 - \ps(u_2)u_1\text{~~for any~~} u_1,u_2\in\Th.\label{(3.10)}
\end{equation}
We want to express the corresponding homomorphism $j^*: H^1(M,\Om^1)\to
H^1(M,\Om^2\otimes\Th)$ in terms of differential forms. Together with the
Dolbeault resolution of $\Om^1$, we use the Dolbeault--Serre resolution
$(\Ph^{2,*}\otimes\Th,\dif)$ of $\Om^2\otimes\Th$ formed by smooth
vector-valued forms. Clearly, $j$ extends to the homomorphism of
resolutions 
\[
\widetilde\jmath = \id\otimes j: \Ph^{0,*}\otimes\Om^1\tto\Ph^{0,*}
\otimes\Om^2\otimes\Th. 
\]
Identifying $\Ph^{0,1}\otimes\Om^1$
 with $\Ph^{1,1}$
and $\Ph^{0,1}\otimes\Om^2\otimes\Th$ with $\Ph^{2,1}\otimes\Th$,
respectively, we see from formula \eqref{(3.8)} that $\widetilde\jmath:\Ph^{1,1}\tto\Ph^{2,1}
\otimes\Th$ is expressed by
\begin{equation}
\widetilde\jmath(\om)(u_1,u_2,v) = \om(u_1,v)u_2 - \om(u_2,v)u_1,\text{~~where~~}u_1,u_2\in
\Th,\;v\in\bar\Th.\label{(3.11)}
\end{equation}
This implies the following Proposition.

\ssbegin[A useful formula]{Proposition}\label{P3.6}
For any $\dif$-closed form $\om\in\Ga(M,\Ph^{1,1})$, the class $j^*(D[\om])
\in  H^1(M,\Om^2\otimes\Th)$ is determined by the $\dif$-closed form
$\widetilde\jmath(\om)$ given by the formula \eqref{(3.11)}.
\end{Proposition}

Now, we apply our construction to the \textit{canonical form} $\om$ defined by Koszul, see \cite{19}.
Let $V$ be a~volume form on a~complex manifold $M$. Then, one associates with
$V$ a~closed (1,1)-form $\om$ in the following way. Let $\mathfrak U = (U_i)$ be
a coordinate cover of $M$ and $x^{(i)}_1,\ldots,x^{(i)}_n$ holomorphic
coordinates in $U_i$. Then, in any $U_i$, we have
$$
V = V_idx^{(i)}_1\ldots dx^{(i)}_nd\overline x^{(i)}_1\ldots
d\overline x^{(i)}_n,
$$
where $V_i$ is a~positive
$C^{\infty}$ function in $U_i$, unique up to a~constant factor and independent of $i$.  Denote
$$
J_{ij} := \nfrac{D(x^{(i)}_1,\ldots,x^{(i)}_n)}{D(x^{(j)}_1,\ldots,x^{(j)}_n)},
$$
then
$$
V_j = |J_{ij}|^2V_i\text{~~in $U_i\cap U_j$. }
$$
The canonical form $\om$ is defined by the formula 
$$
\om = \bpd\hpd\log V_i\;\text{ in }\; U_i
$$
(this definition differs by a~sign from that due to Koszul). Clearly, $d\om = 0$.

\ssbegin[The \, canonical \, supermanifold]{Theorem}\label{T3.4}
The \, supermanifold \, with \, retract $(M,\Om)$, corresponding to the canonical
form $\om$, does not depend of the choice of $V$. It is determined by the
following cocycle $g\in Z^1(\mathfrak U,\mathcal Aut_{(2)}\Om)$:
\begin{equation}
g_{ij} = \id + \nfrac 1{J_{ij}}(dJ_{ij})d.\label{(3.12)}
\end{equation}
\end{Theorem}

\begin{Proof}
The cocycle $\ps$ corresponding to $[\om]$ has the form
$$
\ps_{ij} = \hpd\log V_j - \hpd\log V_i = \hpd\log\nfrac{V_j}{V_i} =
\hpd\log |J_{ij}|^2 = d\log J_{ij} = \nfrac 1{J_{ij}}(dJ_{ij}).
$$
This implies our assertion.
\end{Proof}

The supermanifold, described in Theorem \ref{T3.4}, will be called the \textit{canonical supermanifold}, corresponding to $M$. It is not necessarily
non-split.

\subsection{Lifting of vector fields
}\label{ss3.16}
Let $\M$ be a~supermanifold having $(M,\mathcal O_{\gr})$ as its retract. The
filtration \eqref{(2.2)} of $\mathcal T = \mathcal Der\,\mathcal O$ gives rise to the filtration
$$
\mathfrak v\M = \mathfrak v\M_{(-1)}\supset\mathfrak v\M_{(0)}
\supset\mathfrak v\M_{(1)}\supset\ldots,
$$
where $\mathfrak v\M_{(p)} = \Ga(M,\mathcal T_{(p)})$. By Proposition \ref{P2.1}, we get the
exact sequences
\begin{equation}
0\tto\mathfrak v\M_{(p+1)}\tto\mathfrak v\M_{(p)}\overset{\si_p}\to
\mathfrak v(M,\mathcal O_{\gr})_p\text{ for }p\ge -1.\label{(3.13)}
\end{equation}
We say that a~vector field $u\in\mathfrak v(M,\mathcal O_{\gr})_p$ \textit{lifts to}
$\M$, if $u$
belongs to $\Im\si_p$. In this case, one can suppose that $u = \si_p(\hat u)$,
where $\hat u$ has the same parity as $p$. We are going to express this
property in cohomological terms. Let $\M$ be determined by a~class $\ga\in
H^1(M,\mathcal Aut_{(2)}\mathcal O_{\gr})$.

Suppose that we have an open cover $\mathfrak U = (U_i)_{i\in I}$ of $M$ and a
system of isomorphisms of the sheaves of superalgebras $f_i: \mathcal O|_{U_i}\to
\mathcal O_{\gr}|_{U_i}$ such that $f_i(\ph) = \ph +\mathcal J^{q+1}\in(\mathcal O_{\gr})_q$
for $\ph\in\mathcal J^q$. Then, $g =(g_{ij})$, where $g_{ij} = f_if_j\i$, is a
1-cocycle defining $\ga$.

\ssbegin[Conditions on lifting]{Proposition}\label{P3.7}
\hspace*{-2.7pt}A vector field $v\in\mathfrak v(M,\mathcal O_{\gr})_p$ lifts to $(M, \mathcal O)$ if and only if there exists a~$0$-cochain $(v_i)\in C^0(M,(\mathcal T_{\gr})_{(p)})$ such that
\begin{equation}
v_i\equiv v\bmod(\mathcal T_{\gr})_{(p+1)}(U_i),\label{(3.14)}
\end{equation}
\begin{equation}
g_{ij}v_j = v_ig_{ij}\text{ in } U_i\cap U_j\ne\emptyset.\label{(3.15)}
\end{equation}
In this case, we have
\begin{equation}
[\la_2^*(\ga),v]= 0.\label{(3.16)}
\end{equation}
\end{Proposition}

\begin{Proof}
Suppose that $v$ lifts to $\M$ and  $\hat v\in\mathfrak v(M,\mathcal O)_{(p)}$
satisfies $\si_p(\hat v) = v$. 

Define $v_i\in(\mathcal T_{\gr})_{(p)}(U_i)$ by the formula
$v_i = f_i\hat v f_i\i$. Clearly, $v_i$ satisfies condition \eqref{(3.15)}. Now, for any $\ph
\in\mathcal J^q$, denote $\ph_i := f_i\i(\ph + \mathcal J^{q+1})\in\mathcal J^q$. Then, 
$\hat v(\ph_i) = g_i + h_i$, where $g_i\in f_i\i(\Gr^{p+q}(\mathcal O))$ and $h_i
\in f_i\i(\gr\mathcal O)_{(p+q+1)}) = \mathcal J^{(p+q+1)}$. Hence, 
$$
v_i(\ph + \mathcal J^{q+1}) = (f_i\hat v f_i\i)(\ph + \mathcal J^{q+1}) =
f_i\hat v(\ph_i)\equiv f_i(g_i)\bmod (\gr\mathcal O)_{(p+q+1)}.
$$
On the other hand, we have by definition
$$
v(\ph + \mathcal J^{q+1}) = \hat v(\ph_i) + \mathcal J^{p+q+1} = g_i + \mathcal J^{p+q+1}
= f_i(g_i).
$$
Thus, condition \eqref{(3.14)} is proved.

Conversely, suppose a~cochain $(v_i)\in C^0(M,(\mathcal T_{\gr})_{(p)})$
satisfying conditions \eqref{(3.14)} and \eqref{(3.15)} be given. By condition \eqref{(3.15)}, we have
$$
f_j\i v_jf_j = f_i\i g_{ij}v_jf_j = f_i\i v_ig_{ij}f_j = f_i\i v_if_i.
$$
Then, $\hat v = f_i\i v_if_i$ is a~global section of $\mathcal T_{(p)}$. For any
$\ph\in\mathcal J^q$, we have
$$
\hat v(\ph) = (f_i\i v_if_i)(\ph) = f_i\i v_i(\ph + \mathcal J^{q+1}) =
f_i\i(v(\ph + \mathcal J^{q+1}) + \ps),
$$
where $\ps\in(\gr\mathcal O)_{(p+q+1)}$. It follows that $\hat v(\ph)$ lies in
$v(\ph + \mathcal J^{q+1})\in\mathcal J^{p+q}/\mathcal J^{p+q+1}$. 

Thus,
$\hat v(\ph) + \mathcal J^{p+q+1} = v(\ph + \mathcal J^{q+1})$, and hence $\si_p
(\hat v) = v$.

To prove formula \eqref{(3.16)}, we denote $w_{ij} = \log g_{ij}$ and deduce from condition \eqref{(3.15)} that
\begin{equation}
\begin{aligned}
v_i &= g_{ij}v_jg_{ij}\i = (\exp w_{ij}) v_j (\exp w_{ij})\i =
\Ad_{\exp w_{ij}}v_j \\&= \exp(\ad_{w_{ij}})(v_j)
= v_j + [w_{ij},v_j] + \nfrac1{2!}[w_{ij},[w_{ij},v_j]] +\ldots.
\end{aligned}\label{(3.17)}
\end{equation}
Write $v_i = v_i^{(p)} + v_i^{(p+2)} +\ldots$, where $v_i^{(k)}\in
(\mathcal T_{\gr})_k(U_i)$. By condition \eqref{(3.14)}, $v_i^{(p)} = v$. Then, formula \eqref{(3.17)} implies that
$v_i^{(p+2)} = v_j^{(p+2)} + [\la_2(g)_{ij},v]$. Thus, formula \eqref{(3.16)} is proved.
\end{Proof}

Now, we return to the case where  $\mathcal O_{\gr} = \Om$ and
$\M$ is determined by a~class $\ze\in H^1(M,\mathcal Z\Om^1)$ as in Theorem~\ref{T3.3}. Let a~vector field $u\in\mathfrak v(M,\Om)_p$ be
given. We would like to know, whether $u$ lifts to a~vector field on the
supermanifold $\M$. This problem will be studied in
the following three cases: $u = d$, and $u = l(v)$ as well as $u = i(v)$, where $v\in
\mathfrak v(M)$.

Denote a~cocycle
representing the class $\ze$  by $\ps = (\ps_{ij})\in Z^1(\mathfrak U,\mathcal Z\Om^1)$. Then, $\M$ corresponds to the cohomology class
$\ga = \mu^*(\ze)$ of the cocycle $g = (g_{ij})$ given by the formula \eqref{(3.14)}. We can
suppose that there exist isomorphisms
$f_i: \mathcal O|_{U_i}\tto\Om|_{U_i}$ for $i\in I$, inducing the identity
isomorphisms of the $\mathbb Z$-graded sheaves and such that $g_{ij} =
f_if_j\i$ over $U_i\cap U_j\ne\emptyset$.

\ssbegin[The lift of $d$]{Proposition}\label{P3.8}
The derivation $d\in\mathfrak v(M,\Om)_1$ always lifts to $\M$.
\end{Proposition}

\begin{Proof}
We have
$$
g_{ij}dg_{ij}\i = (\id +\ps_{ij}d)d(\id -\ps_{ij}d) = d.
$$
Applying Proposition \ref{P3.7} (with $v = v_i = d$), we get our assertion.
\end{Proof}

\ssbegin[Technical]{Proposition}\label{P3.9}
If $u\in\mathfrak v(M)$ and $l(u)^*(\ze) = 0$, then $l(u)$ lifts to $\M$.
\end{Proposition}

\begin{Proof}
We can assume that $l(u)(\ps_{ij}) = \al_j - \al_i$ in $U_i\cap U_j$,
where $\al_i\in\mathcal Z\Om^1(U_i)$. Set
$$
v_i = l(u) + \al_i d.
$$
Then, 
$$
\begin{aligned}
g_{ij}v_jg_{ij}\i &= (\id + \ps_{ij}d)v_j(\id - \ps_{ij}d) = v_j +
[\ps_{ij}d,v_j] - \ps_{ij}dv_j(\ps_{ij}d)\\ &= v_j + [\ps_{ij}d,v_j]
= l(u) + \al_j d - l(u)(\ps_{ij})d = l(u) + \al_i d = v_i.
\end{aligned}
$$
Thus, Proposition \ref{P3.7} can be applied.
\end{Proof}

\ssbegin{Corollary}[The lift of $u\in\mathfrak v(M)$ on the canonical $M$]\label{Cor3.20}
If $\M$ is the canonical supermanifold, then $l(u)$ lifts to $\M$ for any
$u\in\mathfrak v(M)$.
\end{Corollary}

\begin{Proof}
We want to prove that, if $\ps_{ij} = d\log J_{ij}$, then $l(u)$ satisfies the condition of Proposition~\ref{P3.9} (see~\eqref{(3.12)}). Denoting $w_i :=
dx^{(i)}_1\ldots dx^{(i)}_n$, we have
$$
w_i = J_{ij}w_j \text{ in } U_i\cap U_j.
$$
Applying $l(u)$, we get
$$
l(u)(w_i) = l(u)(J_{ij})w_j + J_{ij}l(u)(w_j).
$$
Clearly, $l(u)(w_i) = \ph_iw_i$, where $\ph_i\in\mathcal F(U_i)$. It follows
that
$$
\ph_i = \nfrac1{J_{ij}}l(u)(J_{ij}) + \ph_j,
$$
whence
$$
l(u)(\ps_{ij}) = dl(u)(\log J_{ij}) =
d\left(\nfrac1{J_{ij}}l(u)(J_{ij})\right) = d\ph_i - d\ph_j.
$$
This yields our assertion.
\end{Proof}

\subsection{A spectral sequence
}\label{ss3.21}
In this subsection, we consider the following problem, more general than the one
studied in Subsection~\ref{ss3.1}. Let $\M$ be a~supermanifold with retract 
$(M,\mathcal O_{\gr})$. Suppose that $\M$ is determined by a
class $\ga\in H^1(M,\mathcal Aut_{(2)}\mathcal O_{\gr})$. Let us denote $\mathcal T :=
\mathcal Der\,\mathcal O$, $\mathcal T_{\gr} := \mathcal Der\,\mathcal O_{\gr}$. We want to
describe $H^{\bcdot}(M,\mathcal T)$ under the assumption that the bigraded algebra
$H^{\bcdot}(M,\mathcal T_{\gr})$ is known.

We fix an open Stein cover $\mathfrak U = (U_i)_{i\in I}$ of $M$ and
consider the corresponding \v Cech cochain complex $C^*(\mathfrak U,\mathcal T) =
\bigoplus_{p\ge 0} C^p(\mathfrak U,\mathcal T)$. The filtration \eqref{(2.2)} gives
rise to the filtration
\begin{equation}
C^*(\mathfrak U,\mathcal T) = C_{(-1)}\supset C_{(0)}\supset\ldots\supset C_{(p)}
\supset\ldots\supset C_{(m+1)} = 0\label{(3.18)}
\end{equation}
of this complex by the subcomplexes
$$
C_{(p)} := C^*(\mathfrak U,\mathcal T_{(p)}).
$$
Denoting the image of the natural mapping
$H^{\bcdot}(M,\mathcal T_{(p)})\tto H^{\bcdot}(M,\mathcal T)$  by $H(M,\mathcal T)_{(p)}$, we get the filtration
\begin{equation}
H^{\bcdot}(M,\mathcal T) = H(M,\mathcal T)_{(-1)}\supset\ldots
\supset H(M,\mathcal T)_{(p)}\supset
\ldots\supset H(M,\mathcal T)_{(m+1)} = 0.\label{(3.19)}
\end{equation}
Note that the filtration \eqref{(3.18)} is a~filtration of the graded differential algebra
$C^*(\mathfrak U,\mathcal T)$ (under a~bracket determined by the Lie bracket in
$\mathcal T$) by graded differential
subalgebras, and hence the filtration  \eqref{(3.19)} is a~filtration of the graded algebra
$H^{\bcdot}(M,\mathcal T)$ by graded subalgebras. Denote by $\gr H^{\bcdot}(M,\mathcal T)$ the
bigraded algebra associated with the filtration \eqref{(3.19)}; its bigrading is given by the formula
$$
\gr H^{\bcdot}(M,\mathcal T) = \bigoplus _{\substack{p\ge -1\\ q\ge 0}} \Gr^p H^q(M,\mathcal T).
$$
By a~general procedure invented by J. Leray (see \cite[Ch.~III.7]{GM}), the filtration \eqref{(3.18)}
gives rise to a~spectral sequence of bigraded algebras $(E_r,d_r)$
converging to $E_{\infty}\simeq \gr H^{\bcdot}(M,\mathcal T)$. We have
\begin{equation}
d_r(E^{p,q}_r)\subset E^{p+r,q-r+1}_r \text{~~for any $r,\,p,\,q$.}\label{(3.20)}
\end{equation}

The algebra $E_{r+1}$ is identified with the homology algebra $H(E_r,d_r)$.
If we denote $Z_r$ by $\Ker d_r$, then we have the natural homomorphism $\ka^r_{r+1}:
Z_r\to
Z_{r+1}$. For any $s > r$, denote $\ka^r_s := \ka^{s-1}_s\ldots\ka^r_{r+1}$
(this composition is not defined on the entire $Z_r$).

The following theorem is proved in \cite{24}.

\ssbegin[The first three terms of the spectral sequence]{Theorem}\label{T3.5}\ 

\begin{enumerate}
\item[$(1)$]
The first three terms of the spectral sequence $(E_r)$ can be
identified with the following bigraded algebras:
$$
\begin{aligned}
E_0 &= C^*(\mathfrak U,\mathcal T_{\gr}),\\
E_1 &= E_2 = H^{\bcdot}(M,\mathcal T_{\gr}).
\end{aligned}
$$
Here
$$
\begin{aligned}
E_0^{p,q} &= C^{p+q}(\mathfrak U,(\mathcal T_{\gr})_p),\\
E_1^{p,q} &= E_2^{p,q} = H^{p+q}(M,(\mathcal T_{\gr})_p).
\end{aligned}
$$
\item[$(2)$]
$d_{2k+1} = 0$, and hence $E_{2k+1} = E_{2k+2}$ for all $k\ge 0$.
\item[$(3)$]
$d_2 = \ad_{\la_2^*(\ga)}$.
\end{enumerate}
\end{Theorem}

Proposition \ref{P2.1} implies the cohomology exact sequence
$$
H^{p+q}(M,\mathcal T_{(p+1)})\tto H^{p+q}(M,\mathcal T_{(p)})\overset{\si_p^*}
\rightarrow H^{p+q}(M,(\mathcal T_{\gr})_p) = E_2^{p,q}.
$$
We would like to describe the subspace $\Im\si_p^*\subset
H^{p+q}(M,(\mathcal T_{\gr})_p)$ by means of our spectral sequence. An
element
$a\in E_2^{p,q}$ will be called a~\textit{permanent cocycle} if 
\[
\text{$d_2a =
0,\;d_4(\ka^2_4a) = 0,\ \ d_6(\ka^2_6a) = 0$, etc.}
\]
 Let us denote the subspace of permanent cocycles  by
$Z^{p,q}_{\infty}$. 

\ssbegin[Technical, \cite{24}]{Proposition}\label{P3.10}
We have
$$
\begin{aligned}
\si_p^*(H^{p+q}(M,\mathcal T_{(p)}))&\subset Z_{\infty}^{p,q},\\
\si_p(H^0(M,\mathcal T_{(p)})) &= Z_{\infty}^{p,-p}.
\end{aligned}
$$
\end{Proposition}

Thus, a~vector field
$v\in\mathfrak v(M,\mathcal O_{\gr})_p$ lifts to $\M$ if and only if $v$ is a
permanent cocycle of our spectral sequence (and, in particular, satisfies the condition
$d_2v = [\la_2^*(\ga),v] = 0$, cf. Proposition~ \ref{P3.7}).

\section{Applications to flag manifolds
}\label{S4}

\subsection{Flag manifolds and homogeneous vector bundles
}\label{ss4.1}
A flag manifold of a~connected semisimple complex Lie group $G$ is, by
definition, a~complex homogeneous space $M = G/P$, where $P$ is a~parabolic
subgroup of $G$. In this subsection, we fix the notation and summarize some
facts about flag manifolds. For proofs, see \cite{1}, \cite{3}, \cite{23}.

Let $P$ be a~parabolic subgroup of $G$, i.e.,  a~subgroup containing a~Borel
subgroup of $G$. We choose a
maximal algebraic torus $T$ of $G$ lying in $P$ and a~pair of mutually
opposite Borel subgroups $B,\ B_-\supset T$ such that $B_-\subset P$. Let
$\De$ denote the root system of $G$ with respect to $T$, let $\De_+$ be the system
of positive roots corresponding to $B$, and $\Pi\subset\De_+$ the subsystem
of simple roots. Denote
$$
\ga := \nfrac12\sum_{\al\in\De_+}\al.
$$

If $G$ is simple, then we denote by $\de$ the highest (or maximal) root, i.e., 
the highest weight of the adjoint representation $\Ad$ of $G$. This root is
the only maximal element of $\De$ relative the following
partial order in the vector space $\mathfrak t(\mathbb R)^*$:
$$
\la\succeq\mu \text{ if and only if } \la - \mu = \sum_{\al\in\Pi}k_{\al}\al,\;\text
{ where all $k_{\al}$ are non-negative integers}.
$$
In particular, we have the decomposition
\begin{equation}
\de = \sum_{\al\in\Pi}n_{\al}\al,\label{(4.1)}
\end{equation}
where all $n_{\al}$ are positive integers.

Denote by $\De(Q)$ the root system of any Lie subgroup $Q$ of $G$ normalized
by $T$; this is a~subsystem of $\De = \De(G)$. In particular, we have
$$
\De(B_{\pm}) = \De_{\pm}
$$
and
$$
\De(P) = \De_-\cup [S],
$$
where $[S]$ is the
set of all roots that can be expressed as linear combination of elements in the subset $S\subset\Pi$.
Here, $S\ne\Pi$ if $\dim M > 0$.

We have the semidirect decomposition
$$
P = R\rtimes N_-,
$$
where $R$ is the maximal reductive subgroup and $N_-$ is the nilradical of $P$.
Here,
$$
\begin{aligned}
\De(R) &= [S],\\
\De(N_-) &= \De_-\setminus [S],
\end{aligned}
$$
and $S$ is the  system of simple roots of $R$, corresponding to $B\cap R$.
Denote by $N_+$ the unipotent subgroup generated by the root vectors,
the roots of which belong to the set $\De(N_+) = -\De(N_-)$. Then,  for the corresponding
Lie algebras, we have the following decompositions:
\begin{equation}
\begin{aligned}
\mathfrak g &= \mathfrak n_-\oplus\mathfrak r\oplus\mathfrak n_+ = \mathfrak p\oplus\mathfrak n_+,\\
\mathfrak p &= \mathfrak n_-\oplus\mathfrak r.
\end{aligned}\label{(4.2)}
\end{equation}

Denote by $o$ the point $P\in M = G/P$. Due to formulas \eqref{(4.2)}, the holomorphic tangent
space $T_o(M) = T^{1,0}_o(M) = \mathfrak g/\mathfrak p$
can be identified with $\mathfrak n_+$. The isotropy representation $\ta$ of $P$
in $T_o(M)$ is induced by the adjoint representation $\Ad_P$ of $P$ in
$\mathfrak g$. Since $\mathfrak n_+$ is invariant under $\Ad_R$, then $\ta|_R$ is
identified with the representation $\Ad_R$ in $\mathfrak n_+$. It follows that
$\De(N_+)$ is the system of weights of $\ta$ relative to $\mathfrak t$.

On the other hand, it is usual to identify $T^{0,1}_o(M)$ with $\mathfrak n_-$
(see \cite{3}).

Denote by $(-,-)$ the Killing form on $\mathfrak g$. We suppose that in any
root subspace
$\mathfrak g_{\al}$ of $\mathfrak g$, a~basis vector $e_{\al}$ is chosen so that
$(e_{\al},e_{-\al}) = 1$ for $\al\in\De$. Then, the $h_{\al} = [e_{\al},e_{-\al}]$ for 
$\al\in\Pi$, form a~basis of $\mathfrak t$. We also will use the notation
$$
\langle\al,\be\rangle = \nfrac{2(\al,\be)}{(\be,\be)} \text{~~for any~~} \al,\be\in\De.
$$
The Killing form determines an $R$-invariant duality between $\mathfrak n_+$ and
$\mathfrak n_-$. Identifying $\mathfrak n_-$ with $T^{0,1}_o(M)$, we see that the
isotropy representation of $R$ in $T^{0,1}_o(M)$ is induced by $\Ad_R$ and
coincides with $\ta^*|_R$. Its system of weights is $\De(N_-)$.

\subsection{Montgomery's theorem}\label{ssMont} Since $M$ is compact and simply connected, any maximal compact subgroup $K$
of $G$ acts on $M$ transitively, due to Montgomery's theorem, see \cite{Mont}: ``\textit{If $G$ is a connected Lie group which acts transitively on a compact manifold $M$, and if the stabilizer $G_x$ of the point $x\in M$ is connected, then $G$ contains a compact subgroup which acts transitively on $M$.}" 

Montgomery's theorem implies the following

\textbf{Corollary}. \textit{If $G$ is a connected Lie group which acts transitively on
a compact simply connected manifold $M$, then $G$ contains a compact
subgroup which also acts transitively on $M$.}
See also 
\url{https://en.wikipedia.org/wiki/Maximal_compact_subgroup#Existence_and_uniqueness}.

The Cartan-Iwasawa-Malcev theorem asserts that every connected Lie
group (and indeed every connected locally compact group) admits
maximal compact subgroups and that they are all conjugate to one
another. For any semisimple Lie group, uniqueness is a consequence of the
Cartan fixed point theorem, which asserts that if a compact group acts
by isometries on a complete simply connected negatively curved
Riemannian manifold, then it has a fixed point.

Maximal compact subgroups of connected Lie groups are usually not
unique, but they are unique up to conjugation.

Therefore, if $G$ contains a compact subgroup that acts transitively, it
also contains a maximal (under inclusion) compact subgroup which acts
transitively. Now, we have one maximal compact subgroup $K$ which acts
transitively. Any other maximal compact subgroup has the form
$K'=gKg^{-1}$, where $g \in G$. The groups $K'$ also acts transitively. Indeed, for any $x\in M$ we
have $Kx=M$. Therefore,
\[
K'(x) = gK(g^{-1}x) = gK(y) =M. 
\]

Then, $M =
G/P = K/L$, where $L = P\cap K$. We can suppose
that the corresponding real Lie subalgebra $\mathfrak k\subset\mathfrak g$ is spanned
by $ih_{\al}$ for $\al\in\Pi$, and $e_{\al} - e_{-\al},\; i(e_{\al} + e_{-\al})$ for $\al\in\De_+$. Then, we have
$$
\mathfrak g = \mathfrak k(\mathbb C),\ \ \mathfrak r = \mathfrak l(\mathbb C).
$$

The subgroup $L$ is the centralizer of its center in $K$, and hence is a
subgroup of maximal rank. Hence,  the Poincar\'e polynomial of $M$ is expressed
by the \textit{Hirsch formula}, see \cite{Meng}. On the other hand, the Dolbeault cohomology
groups of $M$ satisfy
\begin{equation}
H^{p,q}(M,\mathbb C) = 0 \ \ \text{for } p\ne q\label{(4.3)}
\end{equation}
(see, e.g., \cite{3}). Since $M$ is a~K\"ahler manifold, the Hodge decomposition
yields
$$
H^s(M,\mathbb C) \simeq\begin{cases} H^{p,p}(M,\mathbb C)&\text{for } s = 2p\\ 
0
&\text{for } s = 2p+1.\end{cases}
$$
It follows, in particular, that
\begin{equation}\label{(4.4)}
H^{1,1}(M,\mathbb C)\simeq H^2(M,\mathbb C)\simeq\mathfrak z(\mathfrak r)\simeq\mathbb C^r,
\text{~~where $r = |\Pi\setminus S|$.}
\end{equation}

We can suppose that $G = (\Bih M)^{\circ}$. Then, (see \cite{23})
\begin{equation}
\Bih M = G\rtimes\Si,\label{(4.5)}
\end{equation}
where $\Si$ is a~finite group, naturally isomorphic to the subgroup of
$\Aut\Pi$ leaving $S$ invariant.

It is well known that with any holomorphic linear representation $\ph: P\to
\GL(E)$ one can associate a~holomorphic vector bundle $\mathbf E_{\ph}$ over
$M$. The total space of this bundle is the quotient 
\[
G\times_{\ph}E := (G\times E)/P
\]
of
$G\times E$ by the diagonal action of $P$. The group $G$ acts on
$\mathbf E_{\ph}$ by automorphisms, covering the given action on $M$. This
bundle is called the \textit{homogeneous vector bundle} determined by $\ph$.
For example, take the bundle $\mathbf E_{\ta}$  isomorphic to $\mathbf T(M)$.

The cohomology  $H^{\bcdot}(M,\mathcal E_{\ph}) = \bigoplus_{p\ge 0}
H^{\bcdot}(M,\mathcal E_{\ph})$ admits a~natural $G$-module structure. The corresponding
representation $\Ph$ of $G$ is called \textit{induced}. If $\ph$ is
irreducible (or completely reducible), then the induced representation can
be calculated with the help of an algorithm found by Bott.

Denote by $W$ the
Weyl group of $G$. This group is generated by reflections $\si_{\al}$,
corresponding to the roots $\al\in\De$, but as a~system of generators one
can choose $\{\si_{\al} \mid \al\in\Pi\}$.
As usual, we call a~weight $\la$ of $G$ \textit{dominant} (resp. \textit{strictly
dominant}) if $(\la,\al)\ge 0$ (resp. $> 0$) for all $\al\in\Pi$. The
Bott algorithm is the following operation $\xi\mapsto\xi^*$:
\begin{equation}
\xi^* = \si(\xi + \ga) - \ga,\label{(4.6)}
\end{equation}
where $\xi + \ga$ is regular and $\si\in W$ is chosen in such a~way that
$\xi^*$ is strictly dominant (or $\si(\xi + \ga)$ is dominant, which is the
same). The \textit{index} of $\xi + \ga$ is the number of roots in $\Ph_{\si} =
\si\De_-\cap\De_+$ or, which is the same, the minimal number of factors in a
decomposition of $\si$ into the product of $\si_{\al}$ for $\al\in\Pi$. The index is also
equal to the number of positive roots
$\al$ such that $(\xi + \ga,\al) < 0$.

The result of Bott is as follows (see \cite{1}, \cite{3}, \cite{17}):

\ssbegin[Bott's theorem]{Theorem}\label{T4.1}
Let $\ph: P\tto\GL(E)$ be an irreducible holomorphic representation with
highest weight $\La$. Then, the induced representation can be determined as
follows:
\begin{enumerate}
\item[$(1)$]
If $\La + \ga$ is singular, then $H^{\bcdot}(M,\mathcal E_{\ph}) = 0$.
\item[$(2)$]
If $\La + \ga$ is regular of index $p$, then $H^q(M,\mathcal E_{\ph}) = 0$ for
$q\ne p$ and $H^p(M,\mathcal E_{\ph})$ is an irreducible $G$-module with highest
weight $\La^*$.
\end{enumerate}
\end{Theorem}

This theorem gives, in particular, a~description of the vector space
\[
\Ga(\mathbf E_{\ph}) = \Ga(M,\mathcal E_{\ph}) = H^0(M,\mathcal E_{\ph})
\] of
holomorphic sections of $\mathbf E_{\ph}$. Note that the induced representation
$\Ph: G\tto\GL(\Ga(\mathbf E_{\ph}))$ acts as follows:
\begin{equation}
(\Ph(g)s)(x) = gs(g\i x)\text{ for any }g\in G, \ \ s\in\Ga(\mathbf E_{\ph})\text{ and }x\in M.
\label{(4.7)}
\end{equation}

\ssbegin{Corollary}\label{Cor4.3}
Under the assumptions of Theorem $\ref{T4.1}$, $\Ga(\mathbf E_{\ph})\ne 0$ if and only
if $\La$ is dominant, and in this case $\Ga(\mathbf E_{\ph})$ is an irreducible
$G$-module with highest weight $\La$.
\end{Corollary}

If $\ph$ is completely reducible, then the induced representation can be
calculated as well, by decomposing $\ph$ into irreducible components and
applying Theorem  \ref{T4.1} to the corresponding homogeneous vector bundles. As to
the general case, we only make the following useful remark.

\ssbegin{Corollary}\label{Cor4.4}
Let $\ph: P\tto\GL(E)$ be an arbitrary holomorphic representation and let
$\La$ be a~highest weight of the induced representation $\Ph$ of $G$ in
$\Ga(\mathbf E_{\ph})$. Then, $\La$ is a~highest weight of $\ph$.
\end{Corollary}

\begin{Proof}
Note that a~highest weight of $\ph$ is the same as a~highest weight of the
completely reducible representation $\ph|_R$. By Corollary \ref{Cor4.3}, our assertion is
true whenever $\ph$ is irreducible. Suppose that it is true for $\dim E < m$
and let us prove it for $\dim E = m$. Let $E'$ be a~maximal $P$-submodule
of~ $E$ and denote $E'' := E/E'$. Then, we have the exact sequence of
$G$-sheaves
$$
0\tto\mathcal E'\tto\mathcal E\tto\mathcal E''\tto 0
$$
and the corresponding exact sequence of cohomology with $G$-equivariant
mappings
$$
0\tto\Ga(\mathcal E')\tto\Ga(\mathcal E)\tto\Ga(\mathcal E'').
$$
Let $\La$ be a~highest weight of the $G$-module $\Ga(\mathcal E)$. Since $\Ph$
is completely reducible, $\La$ is a~highest weight of $\Ga(\mathcal E')$ or
$\Ga(\mathcal E'')$. Using the inductive hypothesis and the complete reducibility
of $\ph|_R$, we see that $\La$ is a~highest weight of $\ph$.
\end{Proof}

\subsection{Vector fields on $(M,\Om)$
}\label{ss4.5}
Here we will study the split supermanifold $(M,\Om)$, assuming that $G$ is
simple. Our goal is to calculate the graded Lie superalgebra of vector fields
$\mathfrak v(M,\Om)$ (see \cite{21}).

It is known (see \cite{1}, \cite{23}) that the Lie group $(\Bih\,M)^{\circ}$ is simple
and its isotropy subgroup is parabolic. Thus, we can assume that $G =
(\Bih\,M)^{\circ}$ and $\mathfrak g = \mathfrak v(M)$. Thanks to  \eqref{(2.15)} and \eqref{(2.16)}, we have
\begin{equation}
\begin{aligned}
\mathfrak v(M,\Om)_{-1} &= i(\mathfrak g),\\
\mathfrak v(M,\Om)_0 &= i(\mathfrak{gl}(\mathbf T(M)^*))\; +\!\!\!\!\!\!\supset
l(\mathfrak g).
\end{aligned}\label{(4.8)}
\end{equation}
As in Subsection~\ref{ss1.1}, denote by $\ad_p$ the adjoint representation of
$\mathfrak v(M,\Om)_0$ in
$\mathfrak v(M,\Om)_p$. The following lemma was first proved in \cite{14}.

\sssbegin{Lemma}\label{L4.1}
If $\mathfrak g$ is simple, then $\mathfrak{gl}(\mathbf T(M)^*) = \mathfrak{gl}(\mathbf T(M))
= \langle\id\rangle$. The representation $\ad_{-1}$ is irreducible and
faithful.
\end{Lemma}

\begin{Proof}
From the classical relation (see also \eqref{(2.14)})
$$
[l(u),i(v)] = i([u,v]),\ \ u,v\in\mathfrak g,
$$
we see that $\ad_{-1}$ is irreducible and faithful on $l(\mathfrak g)$. 

Further,
let us regard $\mathfrak{gl}(\mathbf T(M)^*)$ as $H^0(M,\Om^1\otimes\Th) =
\mathfrak{gl}(\mathbf T(M))$. Then, \eqref{(1.4)} implies that
$$
[i(\et),i(v)] = i(\{\et,v\}) = -\et\barwedge v,\ \text{~~ for any $\et\in
\mathfrak{gl}(\mathbf T(M)^*)$ and $v\in\mathfrak g$}.
$$
If $\ad_{-1}i(\et) = 0$, then $\et(v) = 0$ for any $v\in\mathfrak g$. Since
$G$ acts on $M$ transitively, we have 
\[
\text{$\ev_x(\mathfrak g) = T_x(M)$ for all $x
\in M$,}
\]
 and hence $\et = 0$. Thus, $\ad_{-1}$ is faithful on
$i(\mathfrak{gl}(\mathbf T(M)^*))$.

Let $\mathfrak a$ denote the radical of $i(\mathfrak{gl}(\mathbf T(M)^*))$. It is
non-trivial: it contains $\langle\ep\rangle = \langle i(\id)\rangle$. Since
$\ad_{-1}$ is irreducible, its image is a~reductive Lie algebra with
radical $\langle\id\rangle = \langle\ad_{-1}\ep\rangle$. By \eqref{(4.8)},
$\mathfrak a$ coincides with the radical of $\mathfrak v(M,\Om)_0$, and hence
$\mathfrak a~= \langle\ep\rangle$. It follows that $i(\mathfrak{gl}(\mathbf T(M)^*))$
is reductive, and
$$
i(\mathfrak{gl}(\mathbf T(M)^*)) = \langle\ep\rangle\oplus\mathfrak s,
$$
where $\mathfrak s$ is a~semisimple Lie algebra. We have to prove that $\mathfrak s
= 0$.

Clearly, $\mathfrak s$ is invariant under $\ad_0(l(\mathfrak g))$, and hence we get
the homomorphism 
\[
\ad_0(l): \mathfrak g\tto\mathfrak{der}\,\mathfrak s = \ad\mathfrak s
\]
which is injective if $[l(\mathfrak g),\mathfrak s]\ne 0$. In this latter case, we
obtain therefore an injective homomorphism $h: \mathfrak g\tto\mathfrak s$ satisfying
$$
[h(u),z] = [l(u),z]\ \text{~for any $u\in\mathfrak g,\;z\in\mathfrak s$}.
$$
In particular,
\begin{equation}
h([u,v]) = [h(u),h(v)] = [l(u),h(v)]\ \text{~for any $ u,v\in\mathfrak g$}.\label{(4.9)}
\end{equation}

Now, we note that $\mathfrak{gl}(\mathbf T(M)^*) = \Ga(\mathbf T(M)^*\otimes
\mathbf T(M))$ is the vector space of holomorphic sections of the homogeneous
vector bundle $\mathbf T(M)^*\otimes\mathbf T(M) = \mathbf E_{\ph}$, where $\ph =
\ta^*\ta$.   From \eqref{(4.7)} we deduce that the induced representation $\Ph$ of
$G$ in $\Ga(\mathbf T(M)^*\otimes\mathbf T(M))$ satisfies the following condition
\begin{equation}
i(d\Ph(u)\et) = [l(u),i(\et)]\ \ \text{~for any $ u\in\mathfrak g,\;\et\in\mathfrak{gl}(\mathbf T(M)^*)$}.
\label{(4.10)}
\end{equation}
Suppose that $[l(\mathfrak g),\mathfrak s]\ne 0$. Then, eqs.~ \eqref{(4.9)} and
\eqref{(4.10)} imply that $\Im h$ determines a~$G$-sub\-module of $\mathfrak{gl}(\mathbf T(M)^*)$,
where the adjoint representation of $G$ is realized. Thus, the highest root
$\de$ is a~highest weight of $\Ph$. By Corollary \ref{Cor4.4}, $\de$ is
a highest weight of $\ph$. But the weights of $\ph$ have the form $\al -
\be$, where $\al,\be\in\De_+$. This yields a~contradiction.

Thus, we have proved that $[l(\mathfrak g),\mathfrak s] = 0$. It follows that
$\ad_{-1}(\mathfrak s)$ commutes with the irreducible linear Lie algebra
$\ad_{-1}l(\mathfrak g)$. By
the Schur lemma, $\ad_{-1}(\mathfrak s) = 0$, and hence $\mathfrak s = 0$.
\end{Proof}

A graded Lie superalgebra of the form
\begin{equation}
\mathfrak v = \bigoplus_{p\ge 1}\mathfrak v_p\label{(4.11)}
\end{equation}
is called \textit{transitive} if for any $p\ge 1$ it satisfies
\begin{equation}
\{x\in\mathfrak v_p \mid [x,\mathfrak v_{-1}] = 0\} = 0.\label{(4.12)}
\end{equation}
A graded Lie superalgebra of the form \eqref{(4.11)} is called \textit{irreducible} if
the representation $\ad_{-1}$ of $\mathfrak v_0$ is irreducible. All irreducible
transitive complex graded Lie superalgebras of finite dimension were
classified in \cite{15} (see also \cite{31}).

\ssbegin[$\mathfrak v(G/P,\Om)$ for $G$ simple]{Theorem}\label{T4.2}\ 

\begin{enumerate}
\item[$1)$]
For any flag manifold $M = G/P$ of a~simple complex Lie group $G$, the
graded Lie superalgebra $\mathfrak v(M,\Om)$ is transitive and irreducible.
\item[$2)$]
Under the above assumptions, suppose that $\mathfrak v\M = \mathfrak g$. Then, 
$$
\begin{aligned}
\mathfrak v(M,\Om)_{-1} &= i(\mathfrak g),\\
\mathfrak v(M,\Om)_0 &= \langle\ep\rangle\oplus l(\mathfrak g),\\
\mathfrak v(M,\Om)_1 &= \langle d\rangle,\\
\mathfrak v(M,\Om)_p &= 0\text{ for any } p\ge 2.
\end{aligned}
$$
\end{enumerate}
\end{Theorem}

\begin{Proof}
By Lemma \ref{L4.1}, $\mathfrak v(M,\Om)$ is irreducible and satisfies condition \eqref{(4.12)} for
$p = 0$. Thus, we need to prove \eqref{(4.12)} for any $p > 0$.

We will use the following fact: if $\ph\in\Om^p$, where $p > 0$, and if
$i(v)\ph = 0$ for all $v\in\mathfrak g$, then $\ph = 0$. To prove this, we note
that
$$
(i(v)\ph)_x(v_1,\ldots,v_p) = \ph_x(\ev_x(v),v_1,\ldots,v_p),\ \ v_i\in
T_x(M)\text{ for }x\in M.
$$
Since $\mathfrak g$ acts transitively, the condition $i(v)\ph = 0$ yields
$\ph_x = 0$ for any $x\in M$.

Now, suppose that a~vector field $u\in\mathfrak v(M,\Om)_p$ for $p >0$  satisfies
$[u,i(v)] = 0$ for all $v\in\mathfrak g$. Then, for any $f\in\mathcal F$, we have
$$
[i(v),u](f) = i(v)u(f) = 0.
$$
By the above, we have $u(f) = 0$. Therefore,  for any $\ph\in\Om^1$, we get
$$
[i(v),u](\ph) = i(v)u(\ph) + (-1)^{p+1}u(i(v)(\ph)) = i(v)u(\ph) = 0.
$$
Since $u(\ph)\in\Om^{p+1}$, this implies $u(\ph) = 0$. Thus, $u = 0$, and item
(1) is proved.

The item (2) for $p = -1,0$ follows from \eqref{(4.8)} and Lemma \ref{L4.1}. It is
also clear that $\langle d\rangle\subset\mathfrak v(M,\Om)_1$. In particular, we
see that the representation $\ad_{-1}$ of $\mathfrak v(M,\Om)_0\simeq
\mathfrak g$ is the adjoint one, while the representation $\ad_1$ of this
Lie algebra contains a~trivial component of dimension 1. The classification
of transitive irreducible graded Lie superalgebras $\mathfrak v$ given in \cite[Theorem 4]{15},
 shows that if $\mathfrak v$ satisfies the above conditions, then
$\dim\mathfrak v_1 = 1$ and $\mathfrak v_p = 0$ for $p\ge 2$. Thus, item (2) follows from
item (1).
\end{Proof}

\subsection{A family of non-split supermanifolds
}\label{ss4.7}
Here, we apply the construction of Subsection~ \ref{Cor4.3} to the case where $M = G/P$ is
a flag manifold of a~simple complex Lie group $G$. We show that in this
situation one always obtains a~non-empty family of non-split supermanifolds
having $(M,\Om)$ as their retract. We also study holomorphic vector fields on
these supermanifolds.

\sssbegin{Theorem}\label{T4.3}
Let $M = G/P$ be a~flag manifold, where $G$ is simple and $\dim M > 1$, and
denote $r := |\Pi\setminus S|$. Then, there exists a~family of distinct
non-split supermanifolds that have $(M,\Om)$ as their retract, parametrized
by $\mathbb{CP}^{r-1}/\Si$. Here $\Si$ is the finite group from \eqref{(4.5)}.

If $P$ is maximal, then this family consists of a~unique supermanifold,
which is isomorphic to the canonical one.
\end{Theorem}

\begin{Proof}
The group $\Aut\mathbf T(M)^*$ naturally acts on $H^{1,1}(M,\mathbb C)$,
and the mapping 
\[
\mu^*\hat\de: H^{1,1}(M,\mathbb C)\tto H^1(M,\mathcal Aut_{(2)}\Om)
\]
is equivariant. By Theorems \ref{T3.1} and  \ref{T3.3}, this mapping
determines a~family of distinct non-split supermanifolds having $(M,\Om)$ as
their retract which is parametrized by the set
$\left( H^{1,1}(M,\mathbb C) \setminus \{0\}\right) / \Aut\mathbf T(M)^*$. On the other hand,
$\GL(\mathbf T(M)^*) = \mathbb C^{\times}$ (see Lemma  \ref{L4.1}), and \eqref{(2.17)} yields
\[
\Aut\mathbf T(M)^* = \mathbb C^{\times}\times\Bih M. 
\]
Thanks to \eqref{(4.5)}, we see that
$$
\Aut\mathbf T(M)^* = \mathbb C^{\times}\times (G\rtimes\Si).
$$
Clearly, the action of $G$ on $H^{1,1}(M,\mathbb C)$ is
trivial. Using Lemma \ref{L3.2bis}, we deduce that
$$
(H^{1,1}(M,\mathbb C)\setminus\{0\})/\Aut\mathbf T(M)^* =
(H^{1,1}(M,\mathbb C)\setminus\{0\})/(\mathbb C^{\times}\times\Si) =
\mathbb P(H^{1,1}(M,\mathbb C))/\Si.
$$
Due to \eqref{(4.4)}, this implies our first assertion.

To prove the second claim, we note that the canonical supermanifold
corresponding to $M$
is non-split, since in our case the canonical form $\om$ is
positive-definite \cite{19} and hence $[\om]\ne 0$. Thus, it enters
the family just constructed. But if $P$ is maximal, then $r = 1$, and
$\mathbb{CP}^{r-1}$ contains only one point.
\end{Proof}

Now we will study holomorphic vector fields on supermanifolds $\M$ of the
family constructed above, applying Proposition \ref{P2.1} to $\mathcal O_{\gr} = \Om$.
We have to settle what derivations of $\Om$ described in Theorem \ref{T4.2} can be lifted to
$\M$. We need the following lemma.

\sssbegin{Lemma}\label{L4.2}
Let $M = G/P$ be a~flag manifold. Then, both homomorphisms of the sequence
$$
H^{1,1}(M,\mathbb C)\overset{\hat\de}\tto H^1(M,\mathcal Z\Om^1)\overset{\be^*}
\tto H^1(M,\Om^1)
$$
\textup{(see Theorem \ref{T3.3})} are isomorphisms, and $H^1(M,\mathcal Z\Om^1)$ is a~trivial
$G$-module.
\end{Lemma}

\begin{Proof}
Clearly, we have the exact sequence of sheaves
$$
0\tto\mathcal Z\Om^1\overset\be\tto\Om^1\overset d\tto\mathcal Z\Om^2\tto 0,
$$
where $\mathcal Z\Om^2$ is the sheaf of closed forms from $\Om^2$. Consider the
corresponding cohomology exact sequence:
$$
H^0(M,\mathcal Z\Om^2)\tto H^1(M,\mathcal Z\Om^1)\overset{\be^*}\tto H^1(M,\Om^1).
$$
Since $M$ is K\"ahler, all holomorphic forms on it are closed, and hence
\[
H^0(M,\mathcal Z\Om^2) = H^0(M,\Om^2) \simeq H^{2,0}(M,\mathbb C). 
\]
By \eqref{(4.3)}, this
group is trivial, and therefore $\be^*$ is injective. It is also surjective,
since $\be^*\hat\de$ is the Dolbeault isomorphism. It follows that $\be^*$
and $\hat\de$ are isomorphisms. The natural $G$-action  on
$H^1(M,\mathcal Z\Om^1)$ is trivial, since this is true for $H^{1,1}(M,\mathbb C)$.
\end{Proof}

\ssbegin[Technical]{Proposition}\label{P4.1}
Let $M = G/P$ be as in Theorem $\ref{T4.3}$, $\mathfrak g = \mathfrak v(M)$, and let $\M$ be
any non-split supermanifold of the family described in Theorem $\ref{T4.3}$. Then, 
$l(v)$ for $v\in\mathfrak g$, and $d$ can be lifted to $\M$, and we have
$$
\begin{aligned}
\mathfrak v\M_{(0)} &= \mathfrak v\M_{\bar 0}\oplus\langle\hat d\rangle,\\
\mathfrak v\M_{(1)} &= \langle\hat d\rangle,\\
\mathfrak v\M_{(p)} &= 0\text{ for } p\ge 2.
\end{aligned}
$$
Here $\si_0: \mathfrak v\M_{\bar 0}\tto\mathfrak g$ is an isomorphism, $\hat d\ne 0$,
$\si_1(\hat d) = d$ and $[\hat d,\hat d] = [\hat d,v] = 0$ for all
$v\in\mathfrak v\M_{\bar 0}$.
\end{Proposition}

\begin{Proof}
Consider the exact sequence \eqref{(3.13)} for $\mathcal O_{\gr} = \Om$. From Theorem \ref{T4.2}
we deduce that $\mathfrak v\M_{(p)} = 0$ for $p\ge 2$ and that $\si_1:
\mathfrak v\M_{(1)}\tto\mathfrak v(M,\Om)_1 = \langle d\rangle$ is injective. By
Proposition \ref{P3.8}, we see that $\mathfrak v\M_{(1)} = \langle\hat d\rangle$,
where $\hat d$ is odd and $\si_1(\hat d) = d$. For $p = 0$, the exact
sequence has the form
$$
0\tto\mathfrak v\M_{(1)}\tto\mathfrak v\M_{(0)}\overset{\si_0}\tto\mathfrak v(M,\Om)_0 =
\langle\ep\rangle\oplus l(\mathfrak g).
$$
By Lemma \ref{L4.2} and Corollary  \ref{Cor3.20}, 
any $v\in l(\mathfrak g)$ lifts to
$\M$. On the other hand, $\ep$ does not lift by Proposition \ref{P3.7}, since $[\ep,\la_2(\ga)]
= 2\la_2(\ga)\ne 0$  by eq.~
\eqref{(1.2)}. Hence,  $\Im\si_0 = l(\mathfrak g)$. This implies
our assertion concerning $\mathfrak v\M_{(0)}$.

Since $\mathfrak v\M_{(0)}$ is a~subalgebra of $\mathfrak v\M$, it follows that $\hat d$ is a~weight
vector of the representation $\ad_{\bar 1}$ of $\mathfrak v\M_{\bar 0}$, but the
corresponding weight is 0, since $\mathfrak g$ is simple. Thus, $[\hat d,v] = 0$
for all $v\in\mathfrak v\M_{\bar 0}$. It follows that $[\hat d,\hat d]$ lies in
the center of $\mathfrak v\M_{\bar 0}$, whence $[\hat d,\hat d] = 0$.
\end{Proof}

One can ask whether $\mathfrak v\M$ coincides with it subalgebra
$\mathfrak v\M_{(0)}$ calculated in Proposition \ref{P4.1}. This is not true in general,
and  in Theorem \ref{T4.7}  I give the complete answer for the case
where $M$ is an irreducible Hermitian symmetric space.

\subsection{Irreducible Hermitian symmetric spaces
}\label{ss4.10}
A \textit{Hermitian symmetric space} is, by definition, a~connected complex
manifold $M$, endowed with a~Hermitian structure and satisfying the following
condition: for any $x\in M$, there exists a~holomorphic isometry $s_x$ of $M$
such that $d_xs_x = -\id$.

Let $M$ be a~compact Hermitian symmetric space. Let $K$ be the
identity component of the group of all holomorphic isometries of $M$; this is
a compact Lie group. It is known (see \cite{12}) that $M$ is a~homogeneous space of $K$,
and hence can be regarded as the coset space $K/L$, where $L$ is the
stabilizer $K_o$ of a~point $o\in M$.

In what follows, we suppose that $M$ is simply connected and irreducible (as
a Hermitian space). It is known (see \cite{12}) that if $M$ is simply connected, then
$L$ is the centralizer of a~torus in $K$, containing the symmetry $s_o$. Now,
$G = (\Bih\, M)^{\circ}$ is the complexification $G = K(\mathbb C)$, and $M =
G/P$, where $P = G_o$ is a~parabolic subgroup of $G$. Thus, $M$ is a~flag
manifold of a~special type. Now, a~simply connected compact Hermitian
symmetric space $M$ is irreducible if and only if $K$ and $G$ are simple.
In this case, $P$ is maximal.

Let $G$ be a~connected simple complex Lie group. We retain the notation of
Subsection~\ref{ss4.1} and suppose that a~maximal torus $T$ and a~Borel subgroup $B
\supset T$ of $G$ are chosen. Consider the decomposition \eqref{(4.1)} of the
highest root $\de$. A simple root $\al\in\Pi$ will be called \textit{special} if
$n_{\al} = 1$.

Let $P$ be a~parabolic subgroup of $G$ containing the Borel subgroup $B_-$.
In the above notation, the flag manifold $M = G/P$ is Hermitian symmetric
if and only if the subset $S\subset\Pi$ defining $P$ has the form $S =
\Pi\setminus\{\al_0\}$, where $\al_0$ is a~special simple root. Thus, in this
case,
$$
\De(P) = \De_-\cup [\Pi\setminus\{\al_0\}].
$$
It follows that
$$
\De(N_-) = \De_-\setminus [\Pi\setminus\{\al_0\}],
$$
i.e.,  this is the set of those negative roots $-\be$ of $G$, whose expression
through simple roots contains $\al_0$ (necessarily with coefficient 1).
The subgroups $N_+$ and $N_-$ are commutative. The isotropy representation
$\ta: P\tto\GL(\mathfrak n_+)$ is irreducible; in particular, $\ta|N_-$ is
trivial, and $\ta$ is completely determined by its restriction onto $R$.

For all simple Lie groups $G$ of rank $> 1$ that have special simple roots, \cite[Table~6]{OV} shows the Dynkin diagrams of extended systems of simple
roots $\widetilde\Pi = \Pi\cup\{-\de\}$, where for any $\al\in\Pi$ the coefficient $n_{\al}$ is indicated. 

Note that a~simple root is special if and only if it lies in the
same orbit as $-\de$ under the symmetry group of $\widetilde\Pi$. 

On the
other hand, the nontrivial symmetries of $\Pi$, existing for the types
$A_l,\;D_l,\;E_6$, transform special roots into special roots, and we have to
consider special roots up to these symmetries. In all cases, we have
$(\al_0,\al_0) = 2$ for a~special root $\al_0$.

\subsection{The three cases}\label{ss3cases} We also see that any irreducible symmetric Hermitian space of dimension
$\ge 2$ satisfies to one
of the following three conditions, depending on the choice of $G$ and of a
special simple root $\al_0$:

I. $(\de,\al_0) = 0$, and there exists a~unique $\al_1\in\Pi$ such that
$(\al_0,\al_1)\ne 0$.

This case occurs for the groups $G$ of types $B_l,\;C_l,\;D_l,\;E_6,\;E_7$
and for any special root$\al_0$; we have $n_{\al_1} = 2$.

II. $(\de,\al_0) = 0$, and there exist two different simple roots $\al_1,\;
\al_2$, such that $(\al_0,\al_i)\ne 0$ for $i = 1,2$.

This case occurs for the groups $G$ of type $A_l,\;l\ge 3$, for any $\al_0$
corresponding to the interior vertices of the (non-extended) Dynkin diagram.
Here, we have
$n_{\al_1} = n_{\al_2} = 1$. The manifolds $M$ are the Grassmannians
$\Gr^{l+1}_{s}$ for $1 < s < l$.

III. $(\de,\al_0)\ne 0$.

This case occurs for the groups $G$ of type $A_l$ for $l\ge 2$, for any of the two
roots $\al_0$ corresponding to the end vertices of the (non-extended) Dynkin
diagram. There exists a~unique
$\al_1\in\Pi$ such that $(\al_0,\al_1)\ne 0$, and we have $n_{\al_1} = 1$.
The manifolds $M$ are the projective spaces $\mathbb{CP}^l$, where $l\ge 2$.

The roots $\al_1,\;\al_2\in\Pi$, described above, will be called the {\it
neighbors} of $\al_0$. We admit a~numbering of simple roots
$\Pi = \{\al_0,\al_1,\ldots,\al_{l-1}\}$ using this notation. For a~weight
$\la$ of $G$, we will denote by $m(\la)$ the coefficient at $\al_1$ (or the
sum of the coefficients of $\al_1,\;\al_2$) in the expression of $\la$
in terms of $\Pi$. In particular, we see that
$$
m(\de) =\begin{cases} 2 \text{ in the cases I, II}\\ 1 \text{ in the case III}.
\end{cases}
$$

Clearly, the weight system of the irreducible representation $\ta$
coincides with $\De(N_+)$, the highest weight being $\de$ and the lowest
one $\al_0$. Similarly, the weight system of $\ta^*$ is $\De(N_-)$, the
highest weight being $-\al_0$ and the lowest one $-\de$.

\subsection{Invariant vector-valued forms
}\label{ss4.11}
In this Subsection, we discuss invariant vector-valued forms on flag manifolds
and the invariant cohomology $H^{\bcdot}(M,\Om\otimes\Th)^G$ of irreducible
Hermitian symmetric spaces.

Retaining the notation of Subsection~\ref{ss4.1}, consider a~flag manifold $M = G/P =
K/L$. Clearly, $K$ naturally acts on the vector space
$\Ga(M,\Ph\otimes\Th)$ of all smooth vector-valued forms on $M$. The
well-known \'E.~Cartan principle of reducing invariants of a~transitive
action to invariants of the isotropy group (see, e.g., \cite[Theorem ~4.2]{23})
gives

\sssbegin{Proposition}\label{P4.2}
The evaluation mapping of $\Ga(M,\Ph\otimes\Th)$ onto
\[
\bigwedge (T^{1,0}_o(M)\oplus T^{0,1}_o(M))^*\otimes T^{1,0}_o(M)
\]
 given by $\ph\mapsto\ph_o$
determines an isomorphism of the bigraded vector spaces
$$
\begin{aligned}
\Ga(M,\Ph\otimes\Th)^K&\to(\bigwedge (T^{1,0}_o(M)\oplus
T^{0,1}_o(M))^*\otimes T^{1,0}_o(M))^L\\ &= (\bigwedge(\mathfrak n_+\oplus
\mathfrak n_-))^*\otimes\mathfrak n_+)^L\\
&= (\bigwedge(\mathfrak n_+\oplus\mathfrak n_-))^*\otimes\mathfrak n_+)^R
\end{aligned}
$$
preserving the operations $\barwedge$ and $\{-,-\}$. This isomorphism maps
$\Ga(M,\Ph^{p,q}\otimes\Th)^K$ onto 
\[
((\bigwedge^p\mathfrak n_+\otimes
\bigwedge^q\mathfrak n_-)^*\otimes\mathfrak n_+)^R =(\bigwedge^p\mathfrak n_+^*\otimes
\bigwedge^q\mathfrak n_-^*\otimes\mathfrak n_+)^R = (\bigwedge^p\mathfrak n_-\otimes
\bigwedge^q\mathfrak n_+\otimes\mathfrak n_+)^R.
\]
\end{Proposition}

Now we give examples of invariant vector-valued forms.

\ssbegin[$\om\in\Ga(M,\Ph^{1,1})$ on a~complex manifold $M$]{Example}\label{E4.13}
Let \, $M$ \, be \, a \, complex manifold and $\om\in\Ga(M,\Ph^{1,1})$. Consider the
form 
\[
\th_2 = \widetilde\jmath(\om)\in\Ga(M,\Ph^{2,1}\otimes\Th)
\]
given by the formula
\eqref{(3.11)}. Thus,
\begin{equation}
\th_2(u_1,u_2,v) = \om(u_1,v)u_2 - \om(u_2,v)u_1, \ \ u_1,u_2\in\Th,\;v\in
\bar\Th.\label{(4.13)}
\end{equation}
More generally, we can construct the
following vector-valued $(p,p-1)$-form $\th_p$ for $p\ge 1$:
\begin{equation}
\th_p(u_1,\ldots,u_p,v_1,\ldots,v_{p-1}) =
(p-1)!\vmatrix \om(u_1,v_1) & \hdots & \om(u_1,v_{p-1}) & u_1 \\
\om(u_2,v_1) & \hdots & \om(u_2,v_{p-1}) & u_2 \\
\vdots & & \vdots & \vdots \\
\om(u_p,v_1) & \hdots & \om(u_p,v_{p-1}) & u_p\endvmatrix,\label{(4.14)}
\end{equation}
where $u_i\in\Th$ and $v_j\in\bar\Th$. In particular, $\th_1 = \id$ and $\th_2$
is as in \eqref{(4.13)}. Clearly, $\th_p\ne 0$ for $p\le n = \dim M$.
We note that (see \cite{25})
\begin{equation}
\th_p\barwedge\th_q = p\th_{p+q-1}.\label{(4.15)}
\end{equation}

By Proposition \ref{P4.2}, the form $\th_p$ is completely determined by its value
at $o\in M$ which is expressed by the same formula \eqref{(4.14)} through the value
$\om_o$ at $o$.
As $\om_o$, we can choose any $L$-invariant $(1,1)$-form at $o$. For example,
the Killing form on $\mathfrak g$ determines an invariant form $\om$ satisfying
\begin{equation}
\om_o(u,v) = (u,v),\ \ u\in\mathfrak n_+,\;v\in\mathfrak n_-.\label{(4.16)}
\end{equation}
In what follows, we consider the case where $M$ is an irreducible compact
Hermitian symmetric space. Then, the isotropy representation is irreducible,
and hence the form \eqref{(4.16)} is the only (up to a~constant factor) $L$-invariant
$(1,1)$-form on $T_o(M)$.
\end{Example}

\ssbegin[$M = \Gr^{n}_{s}$]{Example}\label{E4.14}
Consider the complex Grassmannian $M = \Gr^{n}_{s}$, where $s$ is an integer in $\{1, \dotsc, n-1\}$. This is
an irreducible compact Hermitian symmetric space with $G = \SL_n(\mathbb C)$ and $
K = \SU_n$. It
is convenient to regard $M$ as a~homogeneous space of the group $G_0 =
\GL_n(\mathbb C)$ with a natural action on $M$ (this action is not effective). As usual, we choose in $G_0$ the maximal torus $T$ consisting of all
diagonal matrices and the Borel subgroup $B$ consisting of all upper triangular
matrices. Then, $B_-$ is the subgroup of all lower triangular matrices. Denote
$r := n - s$. For $o$ we take the point $\langle e_{r+1},\ldots,e_n\rangle\in
\Gr^{n}_{s}$. Then, the isotropy subgroup $P$ of $G$ at $o$ is parabolic and
contains $B_-$; it consists of all matrices of the form
\begin{equation}
\begin{pmatrix} A_1 & 0 \\
V & A_2\end{pmatrix},\label{(4.17)}
\end{equation}
where $A_1\in\GL_r(\mathbb C)$ and $A_2\in\GL_s(\mathbb C)$. Its maximal
reductive subgroup $R$ consists of matrices of the form
\eqref{(4.17)} with $V = 0$ and can be identified  with $\GL_r(\mathbb C)\times\GL_s(\mathbb C)$,
while the unipotent radical $N_-$ is abelian and consists of matrices of the
form \eqref{(4.17)} with $A_1 = I_r$ and $ A_2 = I_s$. The subalgebras $\mathfrak n_-$ and
$\mathfrak n_+$ consist of matrices of the form
$$
\begin{pmatrix} 0 & 0 \\v & 0 \end{pmatrix},\ \ \begin{pmatrix} 0 & u \\0 & 0 \end{pmatrix},
$$
respectively, $v$ being an $(s\x r)$-matrix and $u$ an $(r\x s)$-matrix. We
will  identify $\mathfrak n_-$ and $\mathfrak n_+$ with the vector spaces of matrices
$\operatorname M_{s,r}(\mathbb C)$ and $\operatorname M_{r,s}(\mathbb C)$,
respectively. The isotropy representation $\ta$ of $P$ on $\mathfrak n_+ =
T_o(M)$ is as follows:
\begin{equation}
\ta\left(\begin{pmatrix} A_1 & 0 \\V & A_2\end{pmatrix}\right)(u) = A_1uA_2\i.\label{(4.18)}
\end{equation}

Let us replace the Killing form by the following invariant bilinear form
on $\mathfrak{gl}_n(\mathbb C)$:
$$
(X,Y) = \tr\,XY.
$$
Then,  using \eqref{(4.16)}, we can define the $K$-invariant vector-valued forms
$\th_p$ on $M$ by \eqref{(4.14)}. Now, we construct new examples of $K$-invariant
vector-valued (2,1)- and (3,2)-forms. Note that the same
method permits to construct certain invariant vector-valued $(p,p-1)$-forms
for any $p\ge 1$.

Define the $K$-invariant vector-valued (2,1)-form $\et$ by its $L$-invariant
value
\begin{equation}
\et_o(u_1,u_2,v) = u_1vu_2 - u_2vu_1,\ \ u_1,u_2\in\mathfrak n_+,v\in\mathfrak n_-.
\label{(4.19)}
\end{equation}
The forms $\th_2$ and $\et$ are linearly independent whenever $1 < s < n-1$,
and they coincide for $s = 1$ or $s = n-1$.

Similarly, we define the $K$-in\-var\-i\-ant vector-valued (3,2)-forms $\et_1,
\et_2,\et_3$, whose $L$-in\-var\-i\-ant values at $o$ are as follows:
\begin{equation}
\begin{aligned}
(\et_1)_o &= \Alt\,\tr(u_1v_1u_2v_2)u_3 = \Alt\,(u_1v_1,u_2v_2)u_3\\
&= 2((u_1v_1,u_2v_2)u_3 + (u_2v_1,u_3v_2)u_1 + (u_3v_1,u_1v_2)u_2\\
&- (u_2v_1,u_1v_2)u_3 - (u_3v_1,u_2v_2)u_1 - (u_1v_1,u_3v_2)u_2),
\end{aligned}\label{(4.20)}
\end{equation}
\begin{equation}
\begin{aligned}
(\et_2)_o &= \Alt\,\tr(u_1v_1)u_2v_2u_3 = \Alt\,(u_1,v_1)u_2v_2u_3 \\
&= (u_1,v_1)u_2v_2u_3 + (u_2,v_1)u_3v_2u_1 + (u_3,v_1)u_1v_2u_2\\
&- (u_2,v_1)u_1v_2u_3 - (u_3,v_1)u_2v_2u_1 - (u_1,v_1)u_3v_2u_2\\
&- (u_1,v_2)u_2v_1u_3 - (u_2,v_2)u_3v_1u_1 - (u_3,v_2)u_1v_1u_2 \\
&+ (u_2,v_2)u_1v_1u_3 + (u_3,v_2)u_2v_1u_1 + (u_1,v_2)u_3v_1u_2,
\end{aligned}\label{(4.21)}
\end{equation}
\begin{equation}
\begin{aligned}
(\et_3)_o = \Alt\,u_1v_1u_2v_2u_3 &= u_1v_1u_2v_2u_3 + u_2v_1u_3v_2u_1 +
u_3v_1u_1v_2u_2\\ &- u_2v_1u_1v_2u_3 - u_3v_1u_2v_2u_1 - u_1v_1u_3v_2u_2 \\
&- u_1v_2u_2v_1u_3 - u_2v_2u_3v_1u_1 - u_3v_2u_1v_1u_2 \\&+ u_2v_2u_1v_1u_3 +
u_3v_2u_2v_1u_1 + u_1v_2u_3v_1u_2.
\end{aligned}\label{(4.22)}
\end{equation}
Here, $u_1,u_2,u_3\in\mathfrak n_+$ for $v_1,v_2\in\mathfrak n_-$.
\end{Example}

We will  need the following properties of the forms introduced in Examples
\ref{E4.13} and \ref{E4.14}.

\sssbegin{Lemma}\label{L4.2-}
Suppose that $M = \Gr^{n}_{s}$ and \text{(recall, $r:=n-s$)}
\begin{enumerate}
\item[$(1)$]
The forms $\th_3,\,\et_1,\,\et_2,\,\et_3$ are linearly independent whenever
$s,r\ge 3$.
\item[$(2)$]
If $r = 2$ for $s\ge 3$, then $\th_3,\,\et_1,\,\et_2$ are linearly independent,
while
\begin{equation}
\et_3 =\et_2 + \nfrac12\et_1 - \nfrac12\th_3.\label{(4.23)}
\end{equation}
\item[$(3)$]
If $s = 2$ for $r\ge 3$, then $\th_3,\,\et_1,\,\et_2$ are linearly independent,
while
\begin{equation}
\et_3 = -\et_2 - \nfrac12\et_1 - \nfrac12\th_3.\label{(4.24)}
\end{equation}
\item[$(4)$]
If $s = r = 2$, then $\th_3,\,\et_1$ are linearly independent, while
\begin{equation}
\et_2 = -\nfrac12\et_1,\ \ \et_3 = -\nfrac12\th_3.\label{(4.25)}
\end{equation}
\end{enumerate}
\end{Lemma}

\begin{Proof}
To check the relations \eqref{(4.23)}, \eqref{(4.24)}, \eqref{(4.25)}, we use the following simple
fact: for any two $2\x 2$-matrices $A,B$ we have
\begin{equation}
AB + BA = (\tr A)B + (\tr B)A + (\tr AB - (\tr A)(\tr B))I.\label{(4.26)}
\end{equation}
In the case (2), we obtain \eqref{(4.23)} by applying \eqref{(4.26)} to $A = u_1v_1$ and $B =
u_2v_2$, and alternating the resulting expression of 
\[
u_1v_1u_2v_2u_3 +
u_2v_2u_1v_1u_3 = (u_1v_1u_2v_2 + u_2v_2u_1v_1)u_3. 
\]
Similarly, in
the case (3) we apply \eqref{(4.26)} to $A = v_1u_2$ and $B = v_2u_3$ and alternate
the resulting expression of 
\[
u_1v_1u_2v_2u_3 + u_1v_2u_3v_1u_2 =
u_1(v_1u_2v_2u_3 + v_2u_3v_1u_2). 
\]
Now, \eqref{(4.23)} and \eqref{(4.24)} imply \eqref{(4.25)} in
the case (4).

I skip the proof of linear independence.
\end{Proof}

It is well known (thanks to \'E.~ Cartan) that the real cohomology algebra of a~compact Riemannian
symmetric space $M = K/L$ is naturally isomorphic to the algebra of
$K$-invariant differential forms on $M$ (see, e.g., \cite[Corollary of Theorem 9.7]{23}). We want to prove a~similar assertion concerning
the invariant cohomology $H^{\bcdot}(M,\Om\otimes\Th)^K$ of a~simply connected
compact Hermitian symmetric space $M$.

We use the fine resolution $(\Ph\otimes\Th,\dif)$ of the sheaf
$\Om\otimes\Th$. By the Dolbeault-Serre theorem, the sheaf
cohomology $H^{\bcdot}(M,\Om\otimes\Th)$ and the cohomology of the complex
$(\Ga(M,\Ph\otimes\Th),\dif)$ are isomorphic. Actually, we have 
\[
H^q(M,\Om^p\otimes\Th)
\simeq
H^{p,q}(\Ga(M,\Ph\otimes\Th),\dif). 
\]
Under this isomorphism, the algebraic
and the FN-bracket in $H^{\bcdot}(M,\Om\otimes\Th)$ are induced by the same
operations in $\Ga(M,\Ph\otimes\Th)$. Denote  the operator in
$\Ga(M,\Ph\otimes\Th)$ conjugate to $\dif$  with respect to the given
$K$-invariant Hermitian metric on $M$ by $\dif^*$ and  the Beltrami--Laplace operator by $\square
= \dif\dif^* + \dif^*\dif$. As usual, a~form
$\ph\in\Ga(M,\Ph\otimes\Th)$ is called \textit{harmonic} if $\square\ph = 0$.
For a~harmonic $\ph$, we have $\dif\ph = 0$, and any cohomology class
contains precisely one harmonic form.

\ssbegin{Proposition}\label{P4.3}
Let $M$ be a~simply connected compact Hermitian symmetric space, $K$ the
identity component of the group of all holomorphic isometries of $M$. Then, 
$$
\Ga(M,\Ph^r\otimes\Th)^K = 0\text{~~whenever $r$ is even.}
$$
Moreover, any $\ph\in\Ga(M,\Ph\otimes\Th)^K$ is harmonic, and hence $\dif$-closed.
Assigning to a~form $\ph\in
\Ga(M,\Ph\otimes\Th)^K$ its cohomology class, we get an isomorphism of
bigraded algebras 
\[
\la: \Ga(M,\Ph\otimes\Th)^K\tto H^{\bcdot}(M,\Om\otimes\Th)^G
\]
both for the algebraic and the FN-brackets.

The FN-bracket in $H^{\bcdot}(M,\Om\otimes\Th)^G$ is identically 0.
\end{Proposition}

\begin{Proof}
For any form $\ph\in\Ga(M,\Ph^r\otimes\Th)^K$ we have $s^*\ph = \ph$. Since
$ds_o = -\id$, we see that ${(s^*\ph)_o = (-1)^{r+1}\ph_o}$. If $r$ is even, then
$\ph_o = 0$, and hence $\ph = 0$. This proves the first assertion.

Moreover, in the same situation we have $\dif\ph\in
\Ga(M,\Ph^{r+1}\otimes\Th)^K$. If $r$ is odd, then $\dif\ph = 0$. Similarly,
$\dif^*\ph = 0$, and hence $\ph$ is harmonic. It follows that 
\[
\la:
\Ga(M,\Ph\otimes\Th)^K\tto H^{\bcdot}(M,\Om\otimes\Th)^K =
H^{\bcdot}(M,\Om\otimes\Th)^G
\]
 is defined and injective.
To prove that $\la$ is surjective, suppose that $\ph\in\Ga(M,\Ph\otimes\Th)$
is a~harmonic form representing a~$G$-invariant cohomology class. Then,  for
any $k\in K$, the form $k^*\ph$ is harmonic and lies in the same cohomology
class as $\ph$. Therefore,  $k^*\ph = \ph$ for $k\in K$, so ${\ph\in
\Ga(M,\Ph\otimes\Th)^K}$.

Clearly, $\Ga(M,\Ph\otimes\Th)^K$ is a~subalgebra under both brackets and
$\la$ is an isomorphism of algebras. The FN-bracket is 0, since
$H^q(M,\Om^p\otimes\Th)^G = 0$ whenever $p + q$ is even.
\end{Proof}

\ssbegin{Corollary}\label{Cor4.16}
Under assumptions of Proposition $\ref{P4.3}$, we have \textup{(recall definition of $R$ under \eqref{(4.17)})}
$$
H^q(M,\Om^p\otimes\Th)^G\simeq (\bigwedge^p\mathfrak n_-\otimes\bigwedge^q
\mathfrak n_+\otimes\mathfrak n_+)^R.
$$
\end{Corollary}

Now we are going to calculate certain invariant cohomology groups
assuming that $M$ is irreducible. First of
all, we will find the degrees for which they are non-zero.

\ssbegin{Proposition}\label{P4.4}
For any simply connected irreducible compact Hermitian symmetric space
$M$ of dimension $n\ge 2$ we have $H^q(M,\Om^p\otimes\Th)^G\ne 0$ if and only
if $q = p-1$ for any $p = 1,\ldots,n$.
\end{Proposition}

\begin{Proof}
Let $\om$ be the K\"ahler form on $M$ corresponding to the given K\"ahler
metric. Consider the invariant forms $\th_p$ for $p = 1,\ldots,n$, given by the formula
{(4.14)}. By Proposition \ref{P4.3}, they determine non-zero cohomology classes in
$H^{p-1}(M,\Om^p\otimes\Th)^G$.

By Corollary \ref{Cor4.16}, 
it is sufficient to show that the representation
of $R$ induced in
$\bigwedge^p\mathfrak n_-\otimes\bigwedge^q\mathfrak n_+\otimes\mathfrak n_+$ has no zero
weights whenever $q\ne p-1$. But each weight of this representation has the
form 
\[
\la = (-p+q+1)\al_0 + \sum_{j=1}^{l-1}k_j\al_j. 
\]
If $\la = 0$, then
$q = p-1$.
\end{Proof}

It was proved in \cite{25} that
$$
H^{p-1}(M,\Om^p\otimes\Th)^G \simeq\mathbb C,\ \text{where~~} p = 1,\ldots,n,
$$
for $M = \mathbb{CP}^n$ (case III). We investigate now the degrees $p = 1,2,3$
in the general case.

By Lemma \ref{L4.1}, we have
$$
H^0(M,\Om^1\otimes\Th)^G \simeq\mathbb C.
$$
For the case $p=2$ we need the following fact, implied by a~result of Kostant
(see \cite{17})]. We will denote by $\si_j$ the reflection $\si_{\al_j}\in W$
corresponding to the simple root $\al_j$.

\sssbegin{Lemma}\label{L4.3}
The irreducible components of the $R$-module $\bigwedge^2\mathfrak n_+$
correspond one-to-one to those simple roots $\al_k$ of $G$ that are
neighbors of $\al_0$. The component that corresponds to $\al_k$ has the
lowest weight $2\al_0 + \al_k$ and the lowest weight vector $e_{\al_0}\wedge
e_{\al_0+\al_k}$.

Thus, $\bigwedge^2\mathfrak n_+$ is irreducible in the cases I, III and has
two irreducible components in the case~ II.
\end{Lemma}

\begin{Proof}
By \cite[Section 8]{17}, the irreducible components of $\bigwedge^2\mathfrak n_+$
correspond to those elements $\si\in W$ satisfying
\begin{enumerate}
\item $\si = \si_j\si_k$ for $j\ne k$,
\item $\Ph_{\si} = \si\De_-\cap\De_+\subset\De(N_+)$.
\end{enumerate}
The set $\Ph_{\si} = \{\al,\be\}$ can be determined from the relation
$$
\si\ga = \ga - \al - \be.
$$
The lowest weight of the component corresponding to $\si$
is $\al +\be$, and the lowest weight vector is $e_{\al}\wedge e_{\be}$.

Clearly,
$$
\si\ga = \si_j\si_k\ga = \si_j(\ga - \al_k) = \ga - \al_j - \si_j\al_k.
$$
Hence,  $\al +\be = \al_j + \si_j\al_k$, where $\si_j\al_k = \al_k
-\langle\al_k,\al_j\rangle\al_j$. Since $\al$ and $\be$ contain $\al_0$ with
coefficient 1, the same property must have the roots $\al_j$ and $\si_j\al_k$.
It follows that $j = 0$ and $\langle\al_k,\al_0\rangle\ne 0$, i.e.,  $\al_k$ is
a neighbor of $\al_0$. Since $\langle\al_k,\al_0\rangle = -1$, we have
$$
\al +\be = \al_0 + (\al_0 + \al_k).
$$
This easily implies that $\al = \al_0$ and $\be = \al_0 + \al_k$ (or vice
versa).
\end{Proof}

\ssbegin[$H^1(M,\Om^2\otimes\Th)^G$]{Proposition}\label{P4.5} 
We have
$$
H^1(M,\Om^2\otimes\Th)^G \simeq\begin{cases}\mathbb C&\text{in the cases I, III}\\
\mathbb C^2&\text{in the case II. }\end{cases}
$$
\end{Proposition}

\begin{Proof}
By Corollary \ref{Cor4.3},
$$
H^1(M,\Om^2\otimes\Th)^G\simeq (\bigwedge^2\mathfrak n_-\otimes\mathfrak n_+\otimes
\mathfrak n_+)^R.
$$
Now,
$$
\mathfrak n_+\otimes\mathfrak n_+ = \bigwedge^2\mathfrak n_+\oplus\S^2\mathfrak n_+.
$$
By Lemma \ref{L4.3}, $\bigwedge^2\mathfrak n_+$ is the irreducible
$R$-module with lowest weight $2\al_0 + \al_1$ in the cases I, III. It is
easy to prove that
it is not isomorphic to any submodule of $\S^2\mathfrak n_+$. Indeed, the
lowest weight vector of such a~submodule must be $e_{\al_0}e_{\al_0+\al_1}$, which
is impossible. Thus, $\mathfrak n_+\otimes\mathfrak n_+$ contains precisely one
component dual to $\bigwedge^2\mathfrak n_-$, implying the result. The case II
is considered similarly.
\end{Proof}

\ssbegin{Remark}
Clearly, in the cases I, III, a~basic element of $H^1(M,\Om^2\otimes\Th)^G$
is determined by the invariant form $\th_2$ given by the formula \eqref{(4.13)}, $\om$ being
determined by the formula \eqref{(4.16)}. In the case II,
a basis of $H^1(M,\Om^2\otimes\Th)^G$ is formed by the cohomology classes of
$\th_2$ and $\et$, where $\et$ is given by the formula \eqref{(4.19)}, see Example \ref{E4.14}.
\end{Remark}

The following proposition can be proved by case-by-case verification using
the decompositions into irreducible components. We omit the proof, since
we will not use the result.

\ssbegin[$\dim H^2(M,\Om^3\otimes\Th)^G$]{Proposition}\label{P4.6}
The dimension $k = \dim H^2(M,\Om^3\otimes\Th)^G$ is as follows:
\begin{enumerate}
\item[$(1)$]
$k = 0$ in the case III for $l = 2$ ($M = \mathbb{CP}^2$);
\item[$(2)$]
$k = 1$ in the case I for types $B_l, E_l$ and $D_l$ for $l > 4$, if $M$ is a
quadric, and in the case III for $l > 2$;
\item[$(3)$]
$k = 2$ in the case I for types $C_l$ and $D_l$, if $M$ is the isotropic
Grassmannian of maximal type, and in the case II, if $l = s +1 = 3$ ($M =
\Gr^{4}_{2}$);
\item[$(4)$]
$k = 3$ in the case II whenever $2 < s = l-1$ or $2 = s < l-1$;
\item[$(5)$]
$k = 4$ in the case II whenever $2 < s < l-1$.
\end{enumerate}
\end{Proposition}

\ssbegin{Remark}
In the case II, a~basis of $H^2(M,\Om^3\otimes\Th)^G$ is given by the
cohomology classes of the following forms:
$$
\begin{aligned}
&\th_2,\et_1\text{ for } s = t = 2,\\
&\th_2,\et_1,\et_2\text{ for } s = 2\text{ for }t\ge 3\text{ or } s\ge 3,\; t = 2,\\
&\th_2,\et_1,\et_2,\et_3\text{ for } s,t\ge 3
\end{aligned}
$$
(see Example \ref{E4.14} and Lemma  \ref{L4.2-}).
\end{Remark}

\subsection{An application of a theorem of Bott}\label{ss4.22}
Let again $M$ be an irreducible simply connected compact Hermitian symmetric
space. In this subsection, we apply Theorem \ref{T4.1} to calculation of the
cohomology  $H^q(M,\Om^p\otimes\Th)$ for $q =
1,2$. We regard $\Om^p\otimes\Th$ as the sheaf of holomorphic
sections of the homogeneous vector bundle $\bigwedge^p\mathbf T(M)^*\otimes
\mathbf T(M)$
corresponding to the completely reducible representation $\ta\bigwedge^p
\ta^*$ of $P$.

The following well-known property of dominant weights will be used (see \cite[$\S$~13, Exercise~ 8]{13}).

\sssbegin{Lemma}\label{L4.4}
If $\la$ is a~non-zero dominant weight of a~simple group $G$, then in the expression
$$
\la = \sum_{i=0}^{l-1}k_i\al_i
$$
we have $k_i > 0$ for all $i = 0,\ldots,l-1$.
\end{Lemma}

A weight $\la$ of $G$ will be called $R$-\textit{dominant} if $(\la,\al_i)\ge 0$ for
all $i = 1,\ldots,l-1$. Any highest weight of a~representation of $R$ is,
evidently, $R$-dominant.

Recall that in the theory of Bott the operation $\xi\mapsto\xi^*$ given by the formula
\eqref{(4.6)} is essential. Note that if $\si = \si_i$ is the reflection
corresponding to the simple root $\al_i$, then
\begin{equation}
\xi^* = \si_i\xi - \al_i = \xi - (1 + \langle\xi,\al_i\rangle)\al_i.\label{(4.27)}
\end{equation}

We also need the following lemmas.

\sssbegin{Lemma}\label{L4.5}
Let $\la$ be an $R$-dominant weight of $G$. The weight $\la + \ga$ has
index 1 if and only if $\la^* = \si_0(\la + \ga) - \ga$ is dominant.
\end{Lemma}

\begin{Proof}
Clearly, the condition is sufficient. Now suppose that $\la + \ga$ has
index 1. Then, $\al_0$ is the only positive root of $G$ such that $(\la +
\ga,\al_0) < 0$. For any $i > 0$, we have
$$
(\la^*,\al_i) = (\si_0(\la + \ga),\al_i) - 1 = (\la + \ga,\si_0\al_i) - 1.
$$
Since $\si_0\al_i = \al_i - \langle\al_i,\al_0\rangle\al_0$ is a~positive
root not equal to $\al_0$, this number is non-negative. Also
$$
(\la^*,\al_0) = (\si_0(\la + \ga),\al_0) - 1 = -(\la + \ga,\al_0) - 1\ge 0.
$$
Thus, $\la^*$ is dominant.
\end{Proof}

\sssbegin{Lemma}\label{L4.6}
\begin{enumerate}
\item[$(1)$]
A root $\al\in\De(N_+)$ satisfies $m(\al) = 0$ if and only if
$\al = \al_0$.
\item[$(2)$]
Let $\la$ be a~weight of the representation $\ta\bigwedge^p\ta^*$ for $
p\ge 1$, i.e., 
\begin{equation}
\la = \al - \be_1 -\ldots -\be_p,\label{(4.28)}
\end{equation}
where $\al, \be_i\in\De(N_+),\;\be_i$ are all distinct. Then, 
$$
m(\la)\le\begin{cases} 3-p \text{ in the cases }I, II\\ 2-p \text{ in the case }III.
\end{cases}
$$
If the equality takes place here, then $m(\al) = m(\de)$, one of $\be_i$
coincides with $\al_0$, and we have $m(\be_i) = 1$ for all $\be_i\ne\al_0$.
\end{enumerate}
\end{Lemma}

\begin{Proof}
(1) If $\al$ does not coincide with the lowest root $\al_0$ of $\ta$, then
there exists a~sequence of simple roots $\al_{j_1},\ldots,\al_{j_k}$ such 
that
$$
\al = (\ldots((\al_0+\al_{j_1})+\al_{j_2})+\ldots)+\al_{j_k},
$$
where any sum in parentheses is a~root. In particular, we have that $\al_0+\al_{j_1}\in
\De(N_+)$; whence $(\al_0,\al_{j_1}) < 0$, and $\al_{j_1}$ is a~neighbor of
$\al_0$.

(2) The number $m(\la)$ attains its maximum whenever $m(\al)$ is maximal
(that is, whenever $m(\al) = m(\de)$) and $m(\be_i)$ are minimal (that is, $= 0,1$). Due to
item (1), $m(\be_i) = 0$ for only one root $\be_i = \al_0$.
Therefore, 
$$
m(\la)\le\begin{cases} 2-(p-1) = 3-p \text{ in the cases I, II}\\ 1 - (p-1) = 2-p
\text{ in the case III},
\end{cases}
$$
and the equality takes place in the situation described above.
\end{Proof}

\ssbegin[$H^p(M,\Om^1\otimes\Th)$]{Proposition}\label{P4.7}
We have
$$
\begin{aligned}
H^1(M,\Om^1\otimes\Th)&\simeq\begin{cases}\mathfrak g \text{ in the cases I, II,}\\
0\text{ in the case III,}\end{cases}\\
H^p(M,\Om^1\otimes\Th)&= 0 \text{ for } p\ge 2.
\end{aligned}
$$
\end{Proposition}

\begin{Proof}
The representation $\ta^*\ta$ contains a~unique irreducible component with
highest weight $\la_0 = \de - \al_0$. By \eqref{(4.27)},
$$
\la_0^* = \si_0(\la_0 + \ga) - \ga = \si_0\de = \de - \langle\de,\al_0\rangle
\al_0.
$$
In the cases I and II, $\la_0^* = \de$ is dominant. By Lemma \ref{L4.5}, $\la_0 +\ga$
has index 1, and, by Bott's theorem, we get a~unique $G$-submodule of
$H^1(M,\Om^1\otimes\Th)$ isomorphic to $\mathfrak g$. In the case III, we have
$\langle\la_0,\al_0\rangle = -1$. Therefore,  $\la_0 +\ga$ is singular, and,
by Bott's theorem, our component gives nothing to the cohomology.

Now, it suffices to prove that any non-dominant highest weight $\la$ of
$\ta^*\ta$, such that $\la +\ga$ is regular, coincides with $\la_0$.
Clearly, $\la = \al - \be$, where $\al,\be\in\De(N_+)$ for $\al\ne\be$. Since
$\la$ does not contain $\al_0$, we have
$$
\la = \sum_{j=1}^{l-1}k_j\al_j,
$$
where $k_j\in\mathbb Z$. Since $\la\prec\al\preceq\de$, we have $k_j\le
n_{\al_j}$ for $j = 1,\ldots,l-1$. In particular, $m(\la)\le m(\de)$. Since
$\la$ is $R$-dominant, but not dominant, and $\la +\ga$ is regular, it follows that 
$\langle\la,\al_0\rangle\le -2$. On the other hand,
$\langle\la,\al_0\rangle = -m(\la)$, whence $m(\la)\ge 2$. We see that the
case III is impossible and that in the cases I and II we have $m(\la) =
m(\de) = 2$. Then, Lemma \ref{L4.6}(2) implies that $\be = \al_0$.

Thus, $\la = \al - \al_0$ is the only expression of the weight $\la$ as a
difference of two roots from $\De(N_+)$. It follows that the corresponding
highest vector $v\in\mathfrak n_+\otimes\mathfrak n_-$ has the form
$$
v = e_{\al}\otimes e_{-\al_0}.
$$
But this vector cannot be a~highest one if $\al\ne\de$. Thus, $\la = \la_0$.
\end{Proof}

The next proposition reduces calculation of $H^1(M,\Om^p\otimes\Th)$ for $p\ge 2$ to the results of Subsection~\ref{ss4.5}, where its invariant part has
been calculated.

\ssbegin[$H^1(M,\Om^p\otimes\Th) = H^1(M,\Om^p\otimes\Th)^G$]{Proposition}\label{P4.8}
For $p\ge 2$, we have
$$
H^1(M,\Om^p\otimes\Th) = H^1(M,\Om^p\otimes\Th)^G.
$$
\end{Proposition}

\begin{Proof}
Let $\la$ be a~highest weight of $\ta\bigwedge^p\ta^*$ for $p\ge 2$. Then, $\la$
has the form \eqref{(4.28)}. Hence, 
\begin{equation}
\la = (1-p)\al_0 + \mu,\label{(4.29)}
\end{equation}
where
\begin{equation}
\mu = \sum_{j=1}^{l-1}k_j\al_j\text{ for $k_j\in\mathbb Z,\; k_j\le n_{\al_j},\;
j = 1,\ldots,l-1$}.\label{(4.30)}
\end{equation}
Since $1-p < 0$, it follows that $\la$ is not dominant due to Lemma \ref{L4.5}. Hence
$(\la,\al_0) < 0$. But it is $R$-dominant, and hence $\la + \ga$ has index 1 if and only if $\la^* = \si_0(\la + \ga) -
\ga$ is dominant (see Lemma~\ref{L4.4}). Clearly,
$$
\la^* = \si_0\la - \al_0 = (p-2)\al_0 + \si_0\mu =
(p-2 - \langle\mu,\al_0\rangle)\al_0 + \mu =
(p-2 + m(\la))\al_0 + \mu.
$$
By Bott's theorem we have to show that $\la^*$ cannot be dominant and
non-zero.

Suppose that the weight $\la^*$ is dominant and non-zero. Then,  by Lemma \ref{L4.4},
$$
k_j > 0\text{~~for $j = 1,\ldots, l-1;\ \ m(\la) > 2 - p$}.
$$
Applying Lemma \ref{L4.6}(2), we see that the case III is impossible and that in the
cases I and II we have $m(\la) = 3 - p$. Since $m(\la) > 0$, it follows
that $p = 2$ and $m(\la) = 1$. Thus,
$$
\la^* = \al_0 + \mu,\ \ \la = -\al_0 + \mu.
$$
In the case II, we have $m(\la) = k_1 + k_2 = 2$ which gives a~contradiction.
Now me must consider the case I only.

Clearly, $\al_0 + \al_1 = \si_0\al_1$ is a~positive root of $G$. Since
$\la + \ga$ has index 1 and $(\la,\al_0) < 0$, we get $(\la,\al_0 +
\al_1)\ge 0$. On the other hand,
$$
(\la,\al_0 + \al_1) = (-\al_0 + \al_1 + \sum_{j=2}^{l-1}k_j\al_j,
\al_0 + \al_1) = -2 + (\al_1,\al_1) + \sum_{j=2}^{l-1}k_j(\al_j,\al_1).
$$
If $l\ge 3$, we get $(\la,\al_0 + \al_1) < 0$ which is a~contradiction.
If $l = 2$, then $G$ is of type $B_2$, and $(\la,\al_0 + \al_1) = -1 < 0$,
too.
\end{Proof}

\ssbegin[$H^2(M,\Om^p\otimes\Th)$]{Proposition}\label{P4.9}
For $p = 2$ or $p\ge 4$, we have
$$
H^2(M,\Om^p\otimes\Th) = 0.
$$
Also,
$$
H^2(M,\Om^3\otimes\Th) = H^2(M,\Om^3\otimes\Th)^G.
$$
\end{Proposition}

\begin{Proof}
Let $\la$ be a~highest weight of $\ta\bigwedge^p\ta^*$ for $p\ge 2$. Then,  as in
Proposition \ref{P4.9}, statements \eqref{(4.28)}, \eqref{(4.29)} and \eqref{(4.30)} hold. Similarly, $\la$ is $R$-dominant,
but not dominant, and hence ${(\la,\al_0) < 0}$.
Suppose that the index of $\la + \ga$ is 2. As in the proof of Proposition
\ref{P4.8},
$$
\si_0(\la + \ga) = (p-2 + m(\la))\al_0 + \mu + \ga.
$$
We have
$$
\begin{aligned}
(\si_0(\la + \ga),\al_0) &= - (\la + \ga,\al_0) > 0,\\
(\si_0(\la + \ga),\al_j) &= (\la + \ga,\si_0\al_j) = (\la + \ga,\al_j)> 0,
\end{aligned}
$$
if $\al_j$ is not a~neighbor of $\al_0$. Since the index is equal to 2,
$\si_0(\la + \ga)$ is regular and non-dominant, and hence
$$
(\si_0(\la + \ga),\al_1) < 0
$$
for a~neighbor $\al_1$ of $\al_0$. Then, the weight
$$
\la^* = \si_1\si_0(\la + \ga) - \ga
$$
must be dominant. Using \eqref{(4.27)}, we get
\begin{equation}
\la^* = (p-2+ m(\la))\al_0 + ((-p+2- m(\la))\langle\al_0,\al_1\rangle -
\langle\mu,\al_1\rangle + k_1 - 1)\al_1 + \mu',\label{(4.31)}
\end{equation}
where
$$
\mu' = \sum_{j=2}^{l-1}k_j\al_j.
$$
By Proposition \ref{P4.4}, $\la^*\ne 0$, if $p\ne 3$. Suppose that $\la^*\ne 0$ for
$p = 3$, too. Then,  by Lemma~ \ref{L4.4}, all the coefficients in \eqref{(4.31)} are
positive. In particular,
$$
k_j > 0\text{ for $j = 2,\ldots, l-1$ and $m(\la) > 2 - p$}.
$$
Applying Lemma \ref{L4.6}(2), we see that the case III is impossible and that in the
cases I and II we have $m(\la) = 3 - p$.

Now consider the weight
$$
\widetilde\la = \si_0(\la + \ga) - \ga.
$$
Clearly, $\widetilde\la + \ga$ is of index 1. As we saw, $(\widetilde\la + \ga,\al_1)
< 0$, and hence $\al_1$ is the only positive root with this property. It
follows from formula \eqref{(4.27)} that
$$
\widetilde\la = \al_0 + \mu.
$$
Therefore, 
$$
(\widetilde\la,\al_0) = 2 - m(\la) = p - 1.
$$

To get a~contradiction, we consider separately three cases.

1) Case II. We have, evidently, $k_2 =\ldots = k_{l-1} = 1$, and hence
$$
\widetilde\la = \al_0 + (2-p)\al_1 + \al_2 +\ldots + \al_{l-1}.
$$
Therefore, 
$$
(\widetilde\la,\al_1) = \begin{cases} 3-2p \text{ if $\al_1$ corresponds to an end vertex
of the Dynkin diagram,}\\ 2-2p \text{ otherwise}.\end{cases}
$$
Hence, 
$$
(\widetilde\la,\al_1) < 0
$$
for all $p\ge 2$. If $p = 2$, then the first case is impossible, because
$\widetilde\la + \ga$ is singular.
Now,
$$
(\widetilde\la,\al_1 +\al_0) = (p-1) + (2-2p) = 1 - p < 0
$$
for $p\ge 2$. This gives a~contradiction.

2) Case I, the type of $G$ is not $C_l$. We have
\begin{equation}
\widetilde\la = \al_0 + (3-p)\al_1 + \mu',\label{(4.32)}
\end{equation}
where $\mu'\ne 0$ (since $l\ge 3$) and $\al_1$ is long. Hence, 
$$
(\widetilde\la,\al_1) = 5 - 2p + \sum_{j=2}^{l-1}k_j(\al_j,\al_1) < 5-2p,
$$
and
$$
(\widetilde\la,\al_0 + \al_1) < 4 - p.
$$
This gives a~contradiction whenever $p\ge 4$.

If $p = 2$, then
$$
(\widetilde\la,\al_0 + \al_1) = 1 + (\widetilde\la,\al_1)\le -1,
$$
since $(\widetilde\la,\al_1)\le -2$, and we get a~contradiction as well. 

For
$p = 3$, the same argument shows that $(\widetilde\la,\al_1) = -2$. Then, we see
from \eqref{(4.32)} that there exists precisely one root $\al_j$ (say, for $j = 2$)
such
that $(\al_j,\al_1)\ne 0$, and we have $k_2 = 1$ for $(\al_2,\al_1) = -1$. Then, 
$\al_1 + \al_2\in\De_+$ and hence 
\[
0\le (\widetilde\la,\al_1 + \al_2) = -2 +
(\widetilde\la,\al_2). 
\]
Thus, $(\widetilde\la,\al_2)\ge 2$. But
$$
\widetilde\la = \al_0 + \al_1 +\al_2 +\sum_{j=3}^{l-1}k_j\al_j,
$$
whence 
\[
(\widetilde\la,\al_2) = 1 +\sum_{j=3}^{l-1}k_j(\al_j,\al_2)\le 1.
\]

3) Case I, the type of $G$ is $C_l$. Here, equality \eqref{(4.32)} holds as well, but $\al_1$
is short. Hence, 
$$
(\widetilde\la,\al_1) = 2- p + \sum_{j=2}^{l-1}k_j(\al_j,\al_1)\le 2-p.
$$
But in this case, $\al_0 + 2\al_1\in\De(N_+)$, and
$$
(\widetilde\la,\al_0 + 2\al_1)\le 3-p.
$$
This gives a~contradiction whenever $p\ge 4$.

On the other hand, $\langle\widetilde\la,\al_1\rangle\le -2$, whence
$(\widetilde\la,\al_1)\le -1$, and $(\widetilde\la,\al_0 + 2\al_1)\le p-3$, which
gives a~contradiction for $p = 2$. 

For $p = 3$, we see that the equality
$(\widetilde\la,\al_1)\le -1$ is compatible with \eqref{(4.32)} only for $\widetilde\la =
\al_0$ for $l = 2$. But then $\la = -2\al_0$, and it is easy to see that in this case
\[
\text{$-2\al_0\ne\al -\al_0 -\be_1 -\be_2$ for any $\al,\be_1,\be_2\in\De(N_+)$.}
\]
\end{Proof}

\subsection{Cohomology of $\mathcal T$}\label{ss4.26}
Summarizing the results of Subsections~\ref{ss4.5} and ~\ref{T4.6}, we now
describe the structure of the cohomology  $H^q(M,\Om^p\otimes\Th)$ for
$q = 0,1,2$ under our assumptions about $M$.

\sssbegin{Proposition}\label{P4.10}
Suppose that $M$ is a~simply connected irreducible compact Hermitian
symmetric space of dimension $\ge 2$. The $G$-modules
$H^q(M,\Om^p\otimes\Th)$ for $q = 0,1,2$, are listed in the following tables:

Case I:
$$
\begin{matrix}
& p &\vrule & 0 & 1 & 2 & 3 & 4\ldots \\
q & &\vrule\\
 - & - & - & -- & -- & -- & --- & ---\\
0 & &\vrule & \mathfrak g & \mathbb C & 0 & 0 & 0 \\
1 & &\vrule & 0 & \mathfrak g & \mathbb C & 0 & 0 \\
2 & &\vrule & 0 & 0 & 0 & \mathbb C^k & 0 \\
\end{matrix}
$$

Case II:
$$
\begin{matrix}
& p &\vrule & 0 & 1 & 2 & 3 & 4\ldots \\
q &&\vrule \\
 - & - & - & -- & -- & -- & --- & ---\\
0 & & \vrule& \mathfrak g & \mathbb C & 0 & 0 & 0 \\
1 & & \vrule& 0 & \mathfrak g & \mathbb C^2 & 0 & 0 \\
2 & & \vrule& 0 & 0 & 0 & \mathbb C^k & 0 \\
\end{matrix}
$$

Case III:
$$
\begin{matrix}
& p &\vrule & 0 & 1 & 2 & 3 & 4\ldots \\
q &&\vrule\\
- & - & - & -- & -- & -- & --- & ---\\
0 & &\vrule& \mathfrak g & \mathbb C & 0 & 0 & 0 \\
1 & &\vrule& 0 & 0 & \mathbb C & 0 & 0 \\
2 & &\vrule& 0 & 0 & 0 & \mathbb C^k & 0 \\
\end{matrix}
$$
where we denote by $\mathbb C$ the trivial $G$-module and by $\mathfrak g$ the
adjoint one, and the number $k$ is to be found in Proposition $\ref{P4.6}$.
\end{Proposition}

Due to Proposition \ref{P2.3}, this result permits us to describe $H^q(M,\mathcal T_p)$ for $q = 0,1,2$.

\ssbegin[The $G$-modules
$H^q(M,\mathcal T_p)$]{Theorem}\label{T4.4}
Suppose that $M$ is a~simply connected irreducible compact Hermitian
symmetric space of dimension $\ge 2$. The $G$-modules
$H^q(M,\mathcal T_p)$ for $q = 0,1,2$, where $\mathcal T = \mathcal Der\,\Om$, are listed in
the following tables:

Case I:
$$
\begin{matrix}
& p &\vrule &-1 & 0 & 1 & 2 & 3 & 4\ldots \\
q &&\vrule \\
- & - & - & ---- & ---- & ---- & --- & --- & ---\\
0 & &\vrule & i^*(\mathfrak g) & l^*(\mathfrak g)\oplus i^*(\mathbb C) & l^*(\mathbb C) & 0 & 0
& 0 \\
1 & &\vrule & 0 & i^*(\mathfrak g) & l^*(\mathfrak g)\oplus i^*(\mathbb C) & l^*(\mathbb C) & 0 &
0 \\
2 & &\vrule & 0 & 0 & 0 & i^*(\mathbb C^k) & l^*(\mathbb C^k) & 0 \\
\end{matrix}
$$

Case II:
$$
\begin{matrix}
& p &\vrule &-1 & 0 & 1 & 2 & 3 & 4\ldots \\
q &&\vrule \\
- & - & - & ---- & ---- & ---- & --- & --- & ---\\
0 & &\vrule & i^*(\mathfrak g) & l^*(\mathfrak g)\oplus i^*(\mathbb C) & l^*(\mathbb C) & 0 & 0
& 0 \\
1 & &\vrule & 0 & i^*(\mathfrak g) & l^*(\mathfrak g)\oplus i^*(\mathbb C^2) & l^*(\mathbb C^2) &
0 & 0 \\
2 & &\vrule & 0 & 0 & 0 & i^*(\mathbb C^k) & l^*(\mathbb C^k) & 0 \\
\end{matrix}
$$

Case III:
$$
\begin{matrix}
& p &\vrule &-1 & 0 & 1 & 2 & 3 & 4\ldots \\
q &&\vrule \\
- & - & - & ---- & ---- & ---- & --- & --- & ---\\
0 & &\vrule & i^*(\mathfrak g) & l^*(\mathfrak g)\oplus i^*(\mathbb C) & l^*(\mathbb C) & 0 & 0
& 0 \\
1 & &\vrule & 0 & 0 & i^*(\mathbb C) & l^*(\mathbb C) & 0 & 0 \\
2 & &\vrule & 0 & 0 & 0 & i^*(\mathbb C) & l^*(\mathbb C^k) & 0 \\
\end{matrix}
$$
where we denote by $\mathbb C$ the trivial $G$-module and by $\mathfrak g$ the
adjoint one, and the number $k$ is to be found in Proposition $\ref{P4.6}$.
\end{Theorem}

Using Proposition \ref{P2.3}, it is also possible to calculate
the Lie bracket $[-,-]$ for the part of the algebra $H^{\bcdot}(M,\mathcal T)$ that
is described in Theorem \ref{T4.4}. Here we calculate only the adjoint operator
$\ad\ze$, where $\ze\in H^1(M,\mathcal T_2)$.

Recall a~result of Bott (see \cite[Theorem I and Corollary 2 of Theorem W]{3}, and
also \cite{17}) that describes
the cohomology of a~flag manifold $M = G/P$ with values in the sheaf of
holomorphic sections of a~homogeneous vector bundle $\mathbf E\tto M$ in terms
of the cohomology of the Lie algebra $\mathfrak n_-$. Suppose that $\mathbf E =
\mathbf E_{\ph}$, where $\ph$ is a~holomorphic representation of $P$. In
contrast to Theorem~ \ref{T4.1}, this description is valid for arbitrary $\ph$.

\ssbegin{Theorem}\label{T4.5}
Let a~holomorphic representation
of $G$ in a~finite-dimensional vector space $V$ be given. Then, 
\begin{equation}
\Hom_G(V,H^q(M,\mathcal E))\simeq H^q(\mathfrak n_-,\Hom(V,E_o))^R, \label{(4.33)}
\end{equation}
where the representation of $\mathfrak n_-$ in $V$ is the restriction of the
differential of the given representation of $G$, and that in $E_o$ is the
restriction of $\ph$.
\end{Theorem}

\begin{proof}{\nopoint Proof (a sketch of)} 
By the Dolbeault--Serre theorem,
$H^q(M,\mathcal E)$ can be identified  with the $q$-th cohomology of the complex
$(\Ga(M,\Ph^{0,*}\otimes\mathcal E),\bpd)$ of $\mathbf E$-valued forms of type
$(0,*)$. The vector
space $\Ga(M,\Ph^{0,q}\otimes\mathcal E)$ is the space of smooth sections of the
homogeneous vector bundle $\bigwedge^q\mathbf T^{0,1}(M)^*\otimes\mathbf E$,
whose fiber at $o$ can be identified  with $\bigwedge^q\mathfrak n_-^*\otimes E_o$. By
the Frobenius reciprocity law,
$$
\begin{aligned}
\Hom_K(V,\Ga(M,\Ph^{0,q}\otimes\mathcal E)&\simeq\Hom_L(V,\bigwedge^q\mathfrak n_-^*
\otimes E_o)\\
&=\Hom_R(V,\bigwedge^q\mathfrak n_-^*\otimes E_o) = C^q(\mathfrak n_-,\Hom (V,E_o))^R.
\end{aligned}
$$
The isomorphism here is defined by the formula 
\begin{equation}
h\mapsto \widetilde h, \text{ where $\widetilde h(v)= h(v)(o)$ for any $v\in V$},\label{(4.34)}
\end{equation}
and we denote by $C^q(\mathfrak n_-,\Hom (V,E_o))$ the vector space of
$q$-cochains of the Lie algebra $\mathfrak n_-$ with values in $\Hom (V,E_o)$.
Passing to the cohomology, we get the isomorphism \eqref{(4.33)}.\end{proof}

\ssbegin{Proposition}\label{P4.11}
Let $\ze = l^*([\th])\in H^1(M,\mathcal T_2)$, where $[\th]\in
H^1(M,\Om^2\otimes\Th)$ is the cohomology class of the form $\th\in
\Ga(M,\Ph^{2,1}\otimes\Th)^K$. The map $\ad_\ze: H^0(M,\mathcal T_{-1})\to
H^1(M,\mathcal T_1)$ is as follows:
\begin{enumerate}
\item[$(1)$]
an isomorphism of the $G$-modules
$$
H^0(M,\mathcal T_{-1}) = i^*(H^0(M,\Th))\tto l^*(H^1(M,\Om^1\otimes\Th))
$$
for any $\th\ne 0$ in the case I
and for any $\th = a\th_2 + b\et$, where $a\ne 0$, in the case II;
\item [$(2)$]
0 for $\th = b\et$ in the cases II and III.
\end{enumerate}
\end{Proposition}

\begin{Proof}
For any $w\in\mathfrak g$ we have, by \eqref{(2.14)},
$$
[l(\th),i(w)] = [i(w),l(\th)] = l(\th\barwedge w) - i([w,\th]).
$$
Since $\th$ is $K$-invariant, we see that $[w,\th] = 0$. By Proposition \ref{P2.3},
$[l^*([\th]),i^*(w)]$ is determined by the cocycle $l(\th\barwedge w)$.
Thus, our problem is reduced to the study of the mapping 
\[
H^0(M,\Th)\to
H^1(M,\Om^1\otimes\Th)
\]
defined on the cochain level by $w\mapsto\th\barwedge
w$. Recall that the form $\th\barwedge w\in\Ga(M,\Ph^{1,1}\otimes\Th)$ is
given by the formula
$$
(\th\barwedge w)(u,v) = \th (w,u,v)\text{ for }u\in\Th, v\in\bar\Th.
$$

We will use the isomorphism
\begin{equation}
\Hom_G(\mathfrak g,H^1(M,\Om^1\otimes\Th))\simeq H^1(\mathfrak n_-,
\Hom(\mathfrak g,\mathfrak n_+^*\otimes\mathfrak n_+))^R \label{(4.35)}
\end{equation}
that follows from \eqref{(4.33)} if we identify the fiber $E_o$ of the bundle
$\mathbf E =
\mathbf T(M)^*\otimes\mathbf T(M)$ with $\mathfrak n_+^*\otimes\mathfrak n_+$.
As it was noticed above, this isomorphism on the cochain level is determined
by \eqref{(4.34)}. Let 
\[
h: \mathfrak g\tto\Ga(M,\Ph^{1,1}\otimes\Th)
\;(= \Ga(M,\Ph^{0,1}\otimes\Om^1\otimes\Th))
\]
be given by the formula $h(w) = \th
\barwedge w$. Then, $h$ determines the mapping 
\[
\widetilde h: \mathfrak g\to
(\mathfrak n_+^*\otimes\mathfrak n_-^*)\otimes\mathfrak n_+ = \Hom (\mathfrak n_+\otimes
\mathfrak n_-,\mathfrak n_+)
\]
 given by the formula
$$
\widetilde h(w)(u,v) = \th_o(\pi(w),u,v),\ \ u\in\mathfrak n_+\text{ for }v\in\mathfrak n_-,
$$
where we identify the value $w(o)$ of the vector field $w$ at $o$ with
$\pi(w)$, where $\pi: \mathfrak g\tto\mathfrak n_+$ is the projection along $\mathfrak p$
in the decomposition \eqref{(4.2)}. In order to interprete $\widetilde h(w)$ as an
element of 
\[
\mathfrak n_-^*\otimes(\mathfrak n_+^*\otimes\mathfrak n_+) =
\Hom (\mathfrak n_-,\mathfrak n_+^*\otimes\mathfrak n_+), 
\]
we choose a~basis
$e_1,\ldots,e_n$ of $\mathfrak n_+$ and denote by $e_1^*,\ldots,e_n^*$ the dual
basis of $\mathfrak n_+^*$. Then, 
$$
\widetilde h(w)(v) = \sum_{i=1}^n e_i^*\otimes\widetilde h(w)(u,v) =
\sum_{i=1}^n e_i^*\otimes\th_o(\pi(w),e_i,v),\ \ v\in\mathfrak n_-.
$$
Now, this form is viewed as the following cochain $c_{\ga}\in C^1(\mathfrak n_+^*,
\Hom(\mathfrak g,\mathfrak n_-\otimes\mathfrak n_+))$:
\begin{equation}
c_{\th}(v)(w) = \sum_{i=1}^n e_i^*\otimes\th_o(\pi(w),e_i,v),\ \ v\in
\mathfrak n_-,\; w\in\mathfrak g.\label{(4.36)}
\end{equation}
This cochain is an $R$-invariant cocycle of $\mathfrak n_-$, and we have to
understand what is its cohomology class. By Proposition \ref{P4.10}, we have
$$
H^1(\mathfrak n_-,\Hom(\mathfrak g,\mathfrak n_+^*\otimes\mathfrak n_+))^R\simeq\begin{cases}\mathbb C
&\text{in the cases I, II, }\\ 0&\text{in the case III. }\end{cases}
$$

It is convenient to identify $\mathfrak n_+^*$ with $\mathfrak n_-$ using the
Killing form. Then, we have to consider the cochain complex $C^*(\mathfrak n_-,
\Hom(\mathfrak g,\mathfrak n_-\otimes\mathfrak n_+))^R,\de)$. Let us describe the space
of 1-coboundaries
$\de C^0(\mathfrak n_-,\Hom(\mathfrak g,\mathfrak n_-\otimes\mathfrak n_+))^R$. Clearly,
$$
C^0(\mathfrak n_-,\Hom(\mathfrak g,\mathfrak n_-\otimes\mathfrak n_+))^R = \Hom_R
(\mathfrak g,\mathfrak n_-\otimes\mathfrak n_+).
$$
For any $c\in\Hom_R(\mathfrak g,\mathfrak n_-\otimes\mathfrak n_+)$, we have $\de c(y) =
yc$ for any $y\in\mathfrak n_-$, i.e., 
$$
\de c(y)(z) = c([y,z]),\text{~~for any $y\in\mathfrak n_-,\; z\in\mathfrak g$},
$$
since $d\ta (\mathfrak n_-) = 0$. Clearly, $[\mathfrak n_-,\mathfrak g] = \mathfrak n_-
\oplus\mathfrak r$. Since $c$ is a~homomorphism of $R$-modules, it follows that
$(\de c)(\mathfrak n_-)(\mathfrak g)$ is contained in the vector subspace of
$\mathfrak n_-\oplus\mathfrak n_+$ spanned by all $e_{-\al}\otimes e_{\be}$, where
$\al,\be
\in\De(N_+)$ and $\be - \al\in\De(R)$ or $\al =\be$. 

In the cases I and II
this subspace does not coincide with $\mathfrak n_-\oplus\mathfrak n_+$, i.e.,  there
exist $\al,\be\in\De(N_+)$ such  that $\be - \al\notin\De(R)$ and $\al\ne\be$.
Indeed, we can take $\be = \de$ and $\al = \al_0$.

Suppose that $\th = \th_2$. By \eqref{(4.36)},
$$
\begin{aligned}
c_{\th_2}(v)(w) &= \sum_{i=1}^n e_i^*\otimes((e_i,v)\pi(w) - (\pi(w),v)e_i) \\
&= \sum_{i=1}^n (e_i,v)e_i^*\otimes\pi(w) - (\pi(w),v)\sum_{i=1}^n e_i^*
\otimes e_i\\
&= v\otimes\pi(w) - (\pi(w),v)\sum_{i=1}^n e_i^*\otimes e_i.
\end{aligned}
$$
In particular,
$$
c_{\th_2}(e_{-\al_0})(e_{\de}) = e_{-\al_0}\otimes e_{\de}.
$$
It follows that $c_{\th_2}\notin\de C^0(\mathfrak n_-,\Hom(\mathfrak g,\mathfrak n_-
\otimes\mathfrak n_+))^R$. Thus, $\th = \th_2$ defines a~non-zero homomorphism in
the cases I and II.

Now consider the case II, i.e.,  suppose that $M = \Gr^{n}_{s}(\mathbb C)$,
where $1 < s < n-1$. We will use the notation of Example \ref{E4.14}. Then, 
$\mathfrak n_+$ and $\mathfrak n_-$ are the following subspaces of
$\mathfrak{gl}_n(\mathbb C)$:
\begin{equation}\label{k}
\begin{aligned}
\mathfrak n_+ &= \langle E_{i\be} \mid 1\le i\le s,\; s+1\le\be\le m\rangle,\\
\mathfrak n_- &= \langle E_{\al j} \mid  s+1\le\al\le m,\;1\le j\le s\rangle.
\end{aligned}
\end{equation}
Here, $E_{\al j} = E_{j\al}^*$ form the basis dual to $E_{j\al}$. If $\th =
\et$, then the cochain \eqref{(4.36)} has the form
$$
c_{\et}(v)(w) = \sum_{i,\al} E_{\al i}\otimes (E_{i\al}v\pi(w) - \pi(w)
vE_{i\al}).
$$
We write $v = \sum_{\be j}v_{\be j}E_{\be j},\; \pi(w) = \sum_{j\be}w_{j\be}
E_{j\be}$. Then, 
$$
\begin{aligned}
E_{i\al}v\pi(w) &= E_{i\al}(\sum_{\be j\rh}v_{\be j}w_{j\rh})E_{\be\rh}
= \sum_{j\rh}v_{\al j}w_{j\rh}E_{i\rh},\\
\pi(w)vE_{i\al}) &= (\sum_{jk\rh}w_{j\rh}v_{\rh k})E_{jk} E_{i\al}
= (\sum_{j\rh}w_{j\rh}v_{\rh i})E_{j\al}.
\end{aligned}
$$
Hence, 
$$
c_{\et}(v)(w) = \sum_{ij\al\rh} E_{\al i}\otimes
\sum_{j\rh} (v_{\al j}w_{j\rh})E_{i\rh} - \sum_{ij\al\rh} E_{\al i}\otimes
(\sum_{j\rh}w_{j\rh}v_{\rh i})E_{j\al}.
$$
Consider the 0-cochain $c\in\Hom_R(\mathfrak{gl}_n(\mathbb C),
\mathfrak n_-\otimes\mathfrak n_+)$ given by the formula
$$
\begin{aligned}
c(\mathfrak n_+) &= c(\mathfrak n_-) = 0,\\
c(E_{ij}) &= \sum_{\rh}E_{\rh j}\otimes E_{i\rh},\\
c(E_{\al\be}) &= \sum_kE_{\al k}\otimes E_{k\rh}
\end{aligned}
$$
and restrict it to $\mathfrak g = \mathfrak{sl}_n(\mathbb C)$. Then,  for any
$v\in\mathfrak n_-$ and $w\in\mathfrak g$, we have
$$
\begin{aligned}
\de c(v)(w) &= c([v,w]) = c([v,\pi(w)])\\
&= c(\sum_{j\al\rh}v_{\al j}w_{j\rh}
E_{\al\rh} - \sum_{ij\rh}v_{\rh i}w_{j\rh}E_{ji})\\
&= \sum_{jk\al\rh}v_{\al j}w_{j\rh}E_{\al k}\otimes E_{k\rh} - \sum_{ij\al\rh}
v_{\rh i}w_{j\rh}E_{\al i}\otimes E_{j\al}\\
&= c_{\et}(v)(w).
\end{aligned}
$$
Thus, $c_{\et} = \de c$, and $\et$ defines the zero homomorphism. Evidently,
this implies the statement.
\end{Proof}

\ssbegin{Proposition}\label{P4.12}
Let $M = \Gr^{n}_{s}$ for $2\le s\le n-2$, and let $\th,\ph\in\Ga(M,\Ph^{2,1}
\otimes\Th)^K$. If $n\ge 5$, then $\th\barwedge\ph = 0$ implies $\th = 0$ or
$\ph = 0$. For $M = \Gr^{4}_{2}$, the only solutions of $\th\barwedge\ph = 0$,
up to a~constant factor, are $\th = \sqrt 2\th_2\pm\et,\; \ph = \th_2\pm
\sqrt 2\et$.
\end{Proposition}

\begin{Proof}
By \eqref{(4.15)},
$$
\th_2\barwedge\th_2 = 2\th_3.
$$
From \eqref{(4.20)}, \eqref{(4.21)}, and \eqref{(4.22)} we easily  deduce  the following relations:
$$
\begin{aligned}
\th_2\barwedge\et &= 2(\et_1 + \et_2),\\
\et\barwedge\th_2 &= 4\et_2,\\
\et\barwedge\et &= 4\et_3.
\end{aligned}
$$
Write $\th = a\th_2 + b\et,\; \ph = c\th_2 + d\et$ with $a,b,c,d\in\mathbb C$.
It follows that
$$
\th\barwedge\ph = 2ac\th_3 + 2ad\et_1 + 2(ad + 2bc)\et_2 + 4bd\et_2.
$$

Suppose that $3\le s\le n-3$. By Lemma \ref{L4.2-}, $\th\barwedge\ph = 0$ yields
\[
ac = ad = ad + bc = bd = 0.
\]
Clearly, this implies $(a,b) = 0$ or $(c,d) =
0$.

If $n-s = 2$ and $s\ge 3$, then, by \eqref{(4.23)},
$$
\th\barwedge\ph = 2(ac - bd)\th_3 + 2(ad - bd)\et_1 + 2(ad + 2bc + 2bd)\et_2.
$$
By Lemma \ref{L4.2-}, $\th\barwedge\ph = 0$ yields $ac - bd = ad + bd = ad + 2bc +
2bd = 0$. Clearly, this implies $(a,b) = 0$ or $(c,d) = 0$. The case $s = 2,\;
n - s\ge 3$ is considered similarly.

Suppose now that $ n = 4$ and $k = 2$. It follows from \eqref{(4.25)} that
$$
\th\barwedge\ph = 2(ac - bd)\th_3 + (ad - 2bc)\et_1.
$$
If $\th\barwedge\ph = 0$, then $ac - bd = ad- 2bc = 0$. If $(a,b)\ne 0$, then
this implies 
\[
\vmatrix c & -d\\ d & -2c\endvmatrix = -2c^2 + d^2 = 0, 
\]
whence
$d = \pm\sqrt 2c$. If $(c,d)\ne 0$, then $a = \pm\sqrt 2b$.
\end{Proof}

\subsection{Non-split supermanifolds
}\label{ss4.32}
In this subsection, we apply our results to the problem of classification of
non-split supermanifolds.
Theorem \ref{T4.4} implies that the split supermanifold $(M,\Om)$ satisfies the
conditions of Theorem \ref{T3.2}. Thus, in this case the mapping 
\[
\la_2^*:
H^1(M,\mathcal Aut_{(2)}\Om)\tto H^1(M,\mathcal T_2)
\]
 is bijective. By Theorem \ref{T3.1},
we can parametrize non-split supermanifolds with retract $(M,\Om)$ (up to isomorphism) by
orbits of the group $\Aut\,\mathbf T(M)^*$ in $H^1(M,\mathcal T_2)\setminus\{0\}$.
By Propositions \ref{P4.8} and \ref{P4.2}, one can identify $H^1(M,\mathcal T_2) =
l^*(H^1(M,\Om^2\otimes\Th))$ with the
vector space of $K$-invariant vector-valued (2,1)-forms
$\Ga(M,\Ph^{2,1}\otimes\Th)^K$ using the Dolbeault--Serre isomorphism. We
use this parametrization in the statement of the following classification
theorem.

\sssbegin{Theorem}\label{T4.6}
Suppose that $M$ is a~simply connected irreducible compact Hermitian
symmetric space of dimension $\ge 2$.
\begin{enumerate}
\item[$1)$]
If $M$ is of type I or III, then there exists (up to an isomorphism) precisely
one non-split supermanifold with retract $(M,\Om)$, namely, the canonical
one. The corresponding invariant vector-valued $(2,1)$-form is the form
$\th_2$ given by the formula $\eqref{(4.13)}$, $\om$ being determined by $\eqref{(4.16)}$.
\item[$2)$]
If $M = \Gr^{n}_{s},\; 1 < s < n-1$ is of type II, then non-split
supermanifolds with retract $(M,\Om)$ are parametrized by $\mathbb{CP}^1/\Si$,
where
$$
\Si = \begin{cases}\mathbb Z_2&\text{if } n = 2s\\ \{e\}&\text{otherwise}.\end{cases}
$$
The corresponding invariant vector-valued $(2,1)$-forms are $a\th_2 + b\et$,
where $\et_o$ is given by the formula $\eqref{(4.19)}$ and $a,b\in\mathbb C$ serve as homogeneous
coordinates in $\mathbb{CP}^1$. For $n = 2s$, the action of the generator $\si$
of $\Si$ is expressed in these coordinates as follows: $\si(a\,:\,b) =
(a\,:\,-b)$.
\end{enumerate}
\end{Theorem}

\begin{Proof}
Similarly to the proof of Theorem \ref{T4.3}, we have
$$
(H^1(M,\mathcal T_2)\setminus\{0\})/\Aut\,\mathbf T(M)^* = \mathbb P(H^1(M,\mathcal T_2))/
\Si = \mathbb P(\Ga(M,\Ph^{2,1}\otimes\Th)^K)/\Si.
$$
Then, one applies Proposition \ref{P4.5}.

Suppose that $M = \Gr^{2s}_{s}$ for $s\ge 2$. It is known (one deduces this from
\cite[$\S$~15, Theorem 3]{23}) that the generator
$\si$ of $\Si$, being regarded as a~biholomorphic transformation of $M$, acts
as follows:
$$
\si(gP) = A(g)P,\ \ g\in G = \SL_{2s}(\mathbb C),
$$
where $A$ is the automorphism of $G$ given by the formula
$$
A(g) = \begin{pmatrix} 0 & I_s \\ I_s & 0\end{pmatrix} (g^{\top})\i
\begin{pmatrix} 0 & I_s \\ I_s & 0\end{pmatrix}.
$$
We easily check that the automorphism $d_eA$ acts on $\mathfrak n_{\pm}$ by
$$
d_eA(u) = - u^{\top},\ \ u\in\mathfrak n_{\pm}.
$$
By formula \eqref{(4.19)},
$$
\et_o(-u_1^{\top},-u_2^{\top},-v^{\top}) = \et(u_1,u_2,v)^{\top}.
$$
Therefore,  $\si^*\et = -\et$. Clearly, $\si^*\th_2 = \th_2$. Thus,
$\si^*(a\th_2 + b\et) = a\th_2 + (-b)\et$.
\end{Proof}

Comparing Theorem \ref{T4.6} with Theorem \ref{T4.3}, we see that the construction of
Subsection~\ref{T4.3} gives all non-split supermanifolds with retract $\Om$ in the
cases I and III, while this is not true in the case II (one uses Proposition
{\ref{P3.6}, see Example \ref{E4.13}).

Let us now fix a~non-split supermanifold $\M$ with retract $(M,\Om)$, where
$M$ is a~compact irreducible Hermitian symmetric space. Changing the
notation, we will denote by $\mathcal T$ the tangent sheaf $\mathcal Der\,\mathcal O$ of
$\M$, setting $\mathcal T_{\gr} = \mathcal Der\,\Om$. Our goal is to calculate the
cohomology groups $H^q(M,\mathcal T)$ for $q = 0,1$. These groups depend on the
non-zero form $\th\in\Ga(M,\Ph^{2,1}\otimes\Th)^K$ which parametrizes the
supermanifolds $\M$, as it has been described above.

\ssbegin{Theorem}\label{T4.7}
Let $M$ be a~simply connected irreducible compact Hermitian
symmetric space of dimension $\ge 2$. Let $\M$ be a~non-split supermanifold
with retract $(M,\Om)$, the  tangent sheaf $\mathcal T$, and the corresponding vector-valued form  $\th\in
\Ga(M,\Ph^{2,1}\otimes\Th)^K$.
\begin{enumerate}
\item[$(1)$]
Let $M$ be of type I and $\th = \th_2$, where $\om$ is determined by
$\eqref{(4.16)}$. Then, 
$$
H^0(M,\mathcal T_{\bar 0}) = \mathfrak v\M_{\bar 0}\simeq\mathfrak g\textup{~~(as Lie algebras),}
$$
 while \textup{(as $\mathfrak g$-modules)}
$$
\begin{aligned}
H^0(M,\mathcal T_{\bar 1}) &= \mathfrak v\M_{\bar 1}\simeq\mathbb C\\
H^1(M,\mathcal T) &= H^1(M,\mathcal T_{\bar 0})\simeq\mathfrak g. 
\end{aligned}
$$
 The basic element $\hat d\in\mathfrak v\M_{\bar 1}$
satisfies $[\hat d,\hat d] = 0$.
\item[$(2)$]
Let $M$ be of type II, i.e., $M = \Gr^{n}_{s}$ for $2\le s\le n-2$, and $\th = a\th_2
+ b\et$, where $a\ne 0$. If $n\ge 5$, or $n = 4$ and $(a,b)$ is not
proportional to $(\sqrt 2,\pm1)$, then $H^0(M,\mathcal T)$ is as in $(1)$, while
$$
H^1(M,\mathcal T) = H^1(M,\mathcal T_{\bar 0})\simeq\mathfrak g\oplus\mathbb C \textup{~~(as $\mathfrak g$-modules).}
$$
\item[$(3)$]
Let $M = \Gr^{4}_{2},\; \th =\sqrt 2\th_2 +\et$. Then, $H^0(M,\mathcal T)$ is as in
$(1)$, while \textup{(as $\mathfrak g$-modules)}
$$
\begin{aligned}
H^1(M,\mathcal T_{\bar 0})&\simeq\mathfrak g\oplus\mathbb C,\\
H^1(M,\mathcal T_{\bar 1})&\simeq\mathbb C
\end{aligned}
$$
\item[$(4)$]
Let $M$ be of type II and $\th = \et$. Then, 
$$
H^0(M,\mathcal T_{\bar 0}) = \mathfrak v\M_{\bar 0}\simeq\mathfrak g\textup{~~(as Lie algebras),}
$$
while \textup{(as $\mathfrak g$-modules)}
$$
\begin{aligned}
H^0(M,\mathcal T_{\bar 1}) &= \mathfrak v\M_{\bar 1}\simeq\mathfrak g\oplus\mathbb C,\\
H^1(M,\mathcal T_{\bar 0})&\simeq\mathfrak g\oplus\mathbb C,\\
H^1(M,\mathcal T_{\bar 1})&\simeq\mathfrak g.
\end{aligned}
$$
\item[$(5)$]
Let $M$ be of type III, i.e., $M = \mathbb{CP}^{n-1}$ for $n\ge 3$, and $\th = \th_2 =
\et$. Then, $H^0(M,\mathcal T)$ is as in $(4)$, while \textup{(as $\mathfrak g$-modules)}
$$
H^1(M,\mathcal T) = H^1(M,\mathcal T_{\bar 1})\simeq\begin{cases} 0&\text{for } n\ge 4,\\
\mathbb C&\text{for } n = 3.\end{cases}
$$
\end{enumerate}
\end{Theorem}

\sssec{Comment and Open problem} 
\hspace*{-2.6pt}Due to the isomorphism between $gr_pH^q(M,\mathcal T)$ and $E_{\infty} ^{p,q-p}$, in the proof we need $
H^q(M,\mathcal T)$ only for $q=0$ and 1, so for our purposes
it is not necessary to compute $E^{p,q-p}$ for $q=2$ for any $p$.
Therefore, some terms (denoted by "?" in the tables below) remain unknown and should be calculated for completeness.

\begin{Proof}
Consider the spectral sequence $(E_r)$ associated with $\M$ due to Theorem
\ref{T3.5}. By this theorem, $E_2^{p,q-p} = H^q(M,(\mathcal T_{\gr})_p)$ and $d_2 =
\ad_l^*([\th])$, where $[\th]\in H^1(M,\Om^2\otimes\Th)^G$ is the cohomology
class of $\th$. Clearly, $d_2$ is $G$-equivariant.

We are going to calculate $d_2$ on $E_2^{p,q-p}$ for $q = 0,1$. The case
$q = 0,\; p = -1$ is settled by Proposition \ref{P4.11}. In the case where $q = p = 0$,
we see that
$$
\begin{aligned}
{}[l^*([\th]),l(v)] &= l^*([\th,v]) = l^*([[\th,v]]) = 0,\\
{}[l^*([\th]),\ep] &= - 2l^*([\th]).
\end{aligned}
$$

Clearly, $d_2(E_2^{0,1}) = d_2(E_2^{2,-1}) = 0$. The mapping $d_2: E_2^{0,1}
\tto E_2^{2,0}$ is 0, too, since $E_2^{2,0}$ is a~trivial $G$-module.
Similarly, $d_2 = 0$ on $l^*(H^1(M,\Om^1\otimes\Th)\subset E_2^{1,0}$. 

Now,
for any $\ph\in\Ga(M,\Ph^{2,1}\otimes\Th)^K$ we have
$$
[l(\th),i(\ph)] = [i(\ph),l(\th)] = l(\th\barwedge\ph),
$$
due to \eqref{(2.14)}, since $[\ph,\th] = 0$ by Proposition \ref{P4.2}. By Theorem \ref{T3.3},
$d_2i^*(\ph) = l^*([\th\barwedge\ph])$. This class can be calculated with
the help of Proposition \ref{P4.12} (note that, by Theorem \ref{T4.5}, the forms
$\sqrt 2\th_2\pm\et$ determine isomorphic non-split supermanifolds). This
settles the case $p = q = 1$.

Summarizing, we see that the terms $E_3^{p,q-p} = E_4^{p,q-p}$ for $q = 0,1,2$,
are as follows  (for the definition of $s$, see \eqref{k}):

Case I,\; $\th = \th_2$:
$$
\begin{matrix}
& p &\vrule &-1 & 0 & 1 & 2 & 3 & 4\ldots \\
q &&\vrule \\
- & - & - & - & -- & -- & -- & --- & ---\\
0 & &\vrule & 0 & \mathfrak g & \mathbb C & 0 & 0 & 0 \\
1 & &\vrule & 0 & \mathfrak g & 0 & 0 & 0 & 0 \\
2 & &\vrule & 0 & 0 & 0 & ? & \mathbb C^{s-1} & 0 \\
\end{matrix}
$$

Case II,\; $\th = a\th_2 + b\et,\; a\ne 0$ :
$$
\begin{matrix}
& p &\vrule &-1 & 0 & 1 & 2 & 3 & 4\ldots \\
q &&\vrule \\
- & - & - & -- & -- & -- & -- & --- & ---\\
0 & &\vrule & 0 & \mathfrak g & \mathbb C & 0 & 0 & 0 \\
1 & &\vrule & 0 & \mathfrak g & 0 & \mathbb C & 0 & 0 \\
2 & &\vrule & 0 & 0 & 0 & ? & \mathbb C^{s-2} & 0 \\
\end{matrix}
$$

Case II,\; $\th = \et$ :
$$
\begin{matrix}
& p &\vrule &-1 & 0 & 1 & 2 & 3 & 4\ldots \\
q &&\vrule \\
- & - & - & -- & -- & -- & -- & --- & ---\\
0 & &\vrule & \mathfrak g & \mathfrak g & \mathbb C & 0 & 0 & 0 \\
1 & &\vrule & 0 & \mathfrak g & \mathfrak g & \mathbb C & 0 & 0 \\
2 & &\vrule & 0 & 0 & 0 & ? & \mathbb C^{s-2} & 0 \\
\end{matrix}
$$

Case III,\; $n\ge 4,\;\th = \th_2 =\et$ :
$$
\begin{matrix}
& p &\vrule &-1 & 0 & 1 & 2 & 3 & 4\ldots \\
q &&\vrule \\
- & - & - & -- & -- & -- & -- & -- & ---\\
0 & &\vrule & \mathfrak g & \mathfrak g & \mathbb C & 0 & 0 & 0 \\
1 & &\vrule & 0 & 0 & 0 & 0 & 0 & 0 \\
2 & &\vrule & 0 & 0 & 0 & ? & 0 & 0 \\
\end{matrix}
$$

Case III,\; $n = 3,\;\th = \th_2 =\et$ :
$$
\begin{matrix}
& p &\vrule &-1 & 0 & 1 & 2\ldots \\
q &&\vrule \\
- & - & - & -- & -- & -- & --- \\
0 & &\vrule & \mathfrak g & \mathfrak g & \mathbb C & 0 \\
1 & &\vrule & 0 & 0 & \mathbb C & 0 \\
2 & &\vrule & 0 & 0 & 0 & 0 \\
\end{matrix}
$$

Clearly, for $q = 0,1$ we have $d_4 = d_6 =\ldots = 0$, and hence
$E_3^{p,q-p} = E_{\infty}^{p,q-p}$ for all $p\ge 0$. This implies
our theorem.
\end{Proof}

\ssbegin{Corollary}\label{Cor4.35}
Under assumptions of Theorem $\ref{T4.7}$, we have
$$
\begin{aligned}
\mathfrak v\M_{\bar 0} &\simeq\mathfrak g,\\
\mathfrak v\M_{(0)} &\simeq\mathfrak g\oplus\mathbb C,\\
\mathfrak v\M_{(1)} &\simeq\mathbb C,\\
\mathfrak v\M_{(p)} &= 0\text{ for } p\ge 2.
\end{aligned}
$$

In the cases $(1),\,(2),\,(3)$, $\mathfrak v\M = \mathfrak v\M_{(0)}$, and the
supermanifold $\M$ is not homogeneous. In the remaining cases, $\mathfrak v\M\ne
\mathfrak v\M_{(0)}$.
\end{Corollary}

\begin{Proof}
The claims about $\mathfrak v\M_{(p)}$ are implied by the calculation of the
spectral sequence $(E_r)$. It follows that $\mathfrak v\M_{\bar 0}\simeq\mathfrak g$.

In the cases $(1),\,(2),\,(3)$, we see that $\mathfrak v\M = \mathfrak v\M_{(0)}$.
Therefore,  $\ev_x(v) = 0$ for all $v\in\mathfrak v\M_{\bar 1},\; x\in M$, and
hence $\M$ is not homogeneous.
\end{Proof}

\section{The $\Pi$-symmetric super-Grassmannian
}\label{S5}

Consider the supermanifold $\Pi\Gr^{n|n}_{s|s}$ defined in Example \ref{E1.12}.
Its reduction is the submanifold $M$ of $\Gr^{n}_{s}\x\Gr^{n}_{s}$ consisting of
the vector subsuperspaces $L\subset\mathbb C^{n|n}$ of dimension $s|s$
satisfying $L_{\bar 1} = \Pi(L_{\bar 0})$.
Projecting $M$ onto the first factor, we identify this manifold with
$\Gr^{n}_{s}$. Denoting $r = n - s$, we suppose that $r,s\ge 1$.
Assume that $\Pi$ is given in the standard basis by the matrix
$$
\Pi = \begin{pmatrix} 0 & I_n\\I_n & 0\end{pmatrix},
$$
and define local coordinates in a~neighborhood of the point $o =
\langle e_{r+1},\ldots,e_n, f_1,\ldots,f_s\rangle$ in $M$ identified
with $\langle e_{r+1},\ldots,e_n\rangle\in\Gr^{n}_{s}$. Clearly,
the subsupermanifold $\Pi\Gr^{n|n}_{s|s}$ of $\Gr^{n|n}_{s|s}$
is defined in terms of the coordinate matrix \eqref{(1.14)} by the equations
$$
Y = X,\ \ \H = \Xi.
$$
Thus, the coordinate matrix has the form
\begin{equation}
Z = \begin{pmatrix} X & \Xi\\ I_{s} & 0 \\ \Xi & X \\0 & I_s \end{pmatrix},\label{(5.1)}
\end{equation}
where $X$ and $\Xi$ are $(r\x s)$-matrices. Denoting $X := (x_{i\al})$ and $
\Xi := (\xi_{i\al})$, we get the even local coordinates $x_{i\al}$ and
the odd ones $\xi_{i\al},\; 1\le i\le r,\; 1\le\al\le s$, in a~neighborhood
of the point $o$.

Denote by $Q_{n}(\mathbb C)$ the subsupergroup of $\GL_{n|n}(\mathbb C)$
that preserves $\Pi$. Its coordinate matrix has
the form
\begin{equation}
\begin{pmatrix} A & B\\ B & A\end{pmatrix},\label{(5.2)}
\end{equation}
where $A$ and $B$ are $(n\x n)$-matrices of even and odd coordinates,
respectively, $\det A\ne 0$. The reduction $G_0$ of $Q_{n}(\mathbb C)$ can be
identified, in an obvious way, with $\GL_n(\mathbb C)$. The Lie superalgebra
$\mathfrak{q}_{n}(\mathbb C)$ of $Q_{n}(\mathbb C)$ consists of all complex matrices
of the form \eqref{(5.2)} with arbitrary $(n\x n)$-matrices $A$ and $B$, and its
even part $\mathfrak g_0$ can be identified  with $\mathfrak{gl}_n(\mathbb C)$.

The supermanifold $\Pi\Gr^{n|n}_{s|s}$ admits the standard action of
$Q_{n}(\mathbb C)$, which is expressed in coordinates as the
multiplication of $Z$ from the left by the coordinate matrix \eqref{(5.2)}. This
action induces, clearly, the standard
transitive action of the Lie group $G_0 = \GL_n(\mathbb C)$ on $M =
\Gr^{n}_{s}$. Let $P$ denote the isotropy subgroup of $G_0$ at the
point $o\in M$; it consists of all matrices of the form (4.7). We will use
the notation introduced in Example \ref{E4.14}.

Let us denote by $a\mapsto a^*$ the differential of the standard action of
$Q_{n}(\mathbb C)$ on $\Pi\Gr^{n|n}_{s|s}$. This is a~homomorphism of the Lie
superalgebra $\mathfrak{q}_{n}(\mathbb C)$ into the Lie superalgebra
$\mathfrak v(\Pi\Gr^{n|n}_{s|s})$ of holomorphic vector fields on
$\Pi\Gr^{n|n}_{s|s}$.
In what follows, we need the expression of this homomorphism restricted
to $\mathfrak p$. The holomorphic vector fields on
$\Pi\Gr^{n|n}_{s|s}$ will be written in terms of the local coordinates in a
neighborhood of $o$ given by the matrix \eqref{(5.1)}.
Denote the elements of $\mathfrak p$ by
$$
a_1 = \begin{pmatrix} (a_{ij}) & 0\\ 0 & 0\end{pmatrix},\;
a_2 = \begin{pmatrix} 0 & 0\\ 0 & (b_{\al\be})\end{pmatrix},\;
v = \begin{pmatrix} 0 & 0\\ (v_{\al j}) & 0\end{pmatrix},
$$
where $(a_{ij})\in\mathfrak{gl}_{r}(\mathbb C), \; (b_{\al\be})\in
\mathfrak{gl}_s(\mathbb C)$, and $(v_{\al j})$ is an $(s\x r)$-matrix.
We want to calculate the corresponding fundamental vector fields.

Clearly,
$$
\begin{pmatrix} A_1 & 0 & 0 & 0\\ 0 & I_s & 0 & 0\\
0 & 0 & A_1 & 0\\ 0 & 0 & 0 & I_s\end{pmatrix}
\begin{pmatrix} X & \Xi\\ I_{s} & 0 \\ \Xi & X \\0 & I_s \end{pmatrix} =
\begin{pmatrix} A_1X & A_1\Xi\\ I_s & 0\\A_1\Xi & A_1X\\ 0 & I_s
\end{pmatrix}.
$$
By substituting $A_1 = \exp ta_1$ with $t\in\mathbb C$, by differentiating
at $t = 0$ and changing the signs, we get
\begin{equation}
a_1^*(x_{i\al}) = -(a_1X)_{i\al},\ \ a_1^*(\xi_{i\al}) = -(a_1\Xi)_{i\al},
\label{(5.3)}
\end{equation}
where we identify $a_1$ with $(a_{ij})$. Similarly, we find that
$$
\begin{pmatrix} I_r & 0 & 0 & 0\\ 0 & A_2 & 0 & 0\\
0 & 0 & I_r & 0\\ 0 & 0 & 0 & A_2\end{pmatrix}
\begin{pmatrix} X & \Xi\\ I_{s} & 0 \\ \Xi & X \\0 & I_s \end{pmatrix}\sim
\begin{pmatrix} XA_2\i & \Xi A_2\i\\ I_s & 0\\\Xi A_2\i & XA_2\i\\ 0 & I_r
\end{pmatrix},
$$
whence
\begin{equation}
a_2^*(x_{i\al}) = (Xa_2)_{i\al},\ \ a_2^*(\xi_{i\al}) = (\Xi a_2)_{i\al},
\label{(5.4)}
\end{equation}
where we identify $a_2$ with $(b_{\al\be})$.

Further, for any $t\in\mathbb C$, we get
$$
\begin{pmatrix} I_r & 0 & 0 & 0\\ tv & I_s & 0 & 0\\
0 & 0 & I_r & 0\\ 0 & 0 & tv & I_s\end{pmatrix}
\begin{pmatrix} X & \Xi\\ I_{s} & 0 \\ \Xi & X \\0 & I_s \end{pmatrix} =
\begin{pmatrix} X & \Xi\\ I_s + tvX & tv\Xi\\ \Xi & X\\ tv\Xi & I_r +
tvX\end{pmatrix}.
$$
Multiplying the result from the right by
$$
\begin{pmatrix} I_s + tvX & tv\Xi\\ tv\Xi & I_r + tvX\end{pmatrix}\i =
\begin{pmatrix} I_s - tvX +\ldots & -tv\Xi +\ldots\\ -tv\Xi +\ldots
& I_r - tvX +\ldots\end{pmatrix},
$$
where the omitted terms are of order $> 1$ in $t$, we get the matrix
$$
\begin{pmatrix} X - t(XvX + \Xi v\Xi) +\ldots & \Xi - t(\Xi vX + Xv\Xi)
+\ldots\\ I_s & 0\\ \Xi - t(\Xi vX + Xv\Xi) +\ldots &
X - t(XvX + \Xi v\Xi) +\ldots\\ 0 & I_r\end{pmatrix}.
$$
Therefore, 
\begin{equation}
\begin{aligned}
v^*(x_{i\al}) &= (XvX + \Xi v\Xi)_{i\al},\\
v^*(\xi_{i\al}) &= (\Xi vX + Xv\Xi)_{i\al},
\end{aligned}\label{(5.5)}
\end{equation}
where we identify $v$ with $(v_{\al j})$.

From \eqref{(5.3)}, \eqref{(5.4)}, and\eqref{(5.5)} we get

\ssbegin[Explicit formulas of vector fields]{Proposition}\label{P5.1}
We have
$$
\begin{aligned}
a_1^* &= -\sum_{i,k=1}^r\sum_{\al = 1}^sa_{ik}x_{k\al}\pd{x_{i\al}} -
\sum_{i,k=1}^r\sum_{\al = 1}^sa_{ik}\xi_{k\al}\pd{\xi_{i\al}},\\
a_2^* &= \sum_{\al,\be=1}^s\sum_{i=1}^rb_{\be\al}x_{i\be}\pd{x_{i\al}} +
\sum_{\al\be=1}^s\sum_{i=1}^rb_{\be\al}\xi_{i\be}\pd{\xi_{i\al}},\\
v^* &= \sum_{i,j=1}^r\sum_{\al,\be = 1}^s v_{\be j}(x_{i\be}x_{j\al} +
\xi_{i\be}\xi_{j\al})\pd{x_{i\al}} \\
&+ \sum_{i,j=1}^r\sum_{\al,\be = 1}^s v_{\be j}(\xi_{i\be}x_{j\al} +
x_{i\be}\xi_{j\al})\pd{\xi_{i\al}}.
\end{aligned}
$$
\end{Proposition}

Let $\mathcal O$ denote the structure sheaf of the supermanifold
$\Pi\Gr^{n|n}_{s|s}$. Clearly, the action of $G$
on $\M$ determines a~linear representation of the group $P$ by automorphisms
of the superalgebra $\mathcal O_o$, which gives a~linear representation
$\ch = \ch_{\bar 0} + \ch_{\bar 1}$ of this group in $T_o(M,\mathcal O)$, called
the \textit{isotropy representation}. Proposition \ref{P5.1} easily implies its explicit expression.

Indeed, denote  the tautological representations of
$\GL_r(\mathbb C)$ and $\GL_s(\mathbb C)$ by $\rho_1$ and $\rho_2$, respectively. Let $\widetilde m_o$ be the
linear span of germs at $o$ of all coordinate functions $x_{i\al},\;
\xi_{i\al}$ in $m_o$. Then, $m_o = \widetilde m_o\oplus m_o^2$. As Proposition \ref{P5.1}
shows, $v^*(\widetilde m_o)\subset m_o^2$ for all $v\in\mathfrak n_-$, and hence
$\mathfrak n_-$ trivially acts on $m_o/m_o^2$. The same proposition implies that
$\widetilde m_o$ is invariant under $\mathfrak r$ (or $R$), inducing in both
components $(\widetilde m_o)_{\bar 0}$ and $(\widetilde m_o)_{\bar 1}$ the
representation $\rh_1^*\otimes\rh_2$ of $R$.

As in Example \ref{E4.14}, we consider the maximal algebraic torus $T$ of $R$ and
$G_0$ consisting of all diagonal matrices. We will write the matrices of
the corresponding Cartan subalgebra $\mathfrak t$ in the form
$$
H = \diag(\la_1,\ldots,\la_t,\la_{t+1},\ldots,\la_n),\; \la_i\in\mathbb C.
$$
Proposition \ref{P5.1} also implies that the germs of $x_{i\al},\;\xi_{\al i},\;
\et_{i\al}$ form a~weight basis for the representation $\ch^*$ in
$\widetilde m_o\simeq m_o/m_o^2$ with respect to $T$, the corresponding weights
being $-\la_i +\la_{t+\al}$, where $1\le i\le r,\;1\le\al\le s$ (with multiplicity
2). Thus, we got

\ssbegin{Proposition}\label{P5.2}
\begin{enumerate}
\item[$(1)$] The isotropy representation $\ch$ is completely reducible, and the
restrictions of its even and odd components onto $R$ are as follows:
$$
\ch_{\bar 0 }|_R\simeq\ch_{\bar 1}|_R\simeq\rh_1\otimes\rh_2^*.
$$
\item[$(2)$] The germs of $x_{i\al},\;\xi_{i\al}$ form a~weight basis
with respect to $T$ in their linear span $\widetilde m_o$, the corresponding
weights being in both cases $-\la_i +\la_{t+\al}$, where $1\le i\le r$ and $
1\le\al\le s$.
\end{enumerate}
\end{Proposition}

Note that $\ch_{\bar 0}$ coincides with the isotropy representation $\ta$ of
the homogeneous space $\Gr^{n}_{s}(\mathbb C)$ (see \eqref{(4.18)}).

Clearly, the action of $G_0$ on the sheaf $\mathcal O$ leaves invariant the
filtration \eqref{(1.12)} and induces an action of this group on the locally free
sheaf $\mathcal E = \mathcal J/\mathcal J^2$, and hence on the corresponding vector
bundle $\mathbf E$, covering its standard action on $M$. Thus, $\mathbf E$ is a
homogeneous vector bundle over $M$.

\ssbegin{Proposition}\label{P5.3}
The vector bundle $\mathbf E$ is isomorphic to the cotangent bundle
$\mathbf T(M)^*$. The retract of the super-Grassmannian $\M$ is isomorphic
to the supermanifold $(M,\Om)$ from Example~$\ref{E1.9}$.
\end{Proposition}

\begin{Proof}
By Proposition \ref{P5.2}, the representation of $P$ in $E_o = T_o\M_{\bar 1}^*$
is isomorphic to~ $\ta^*$. Hence,  $\mathbf E \simeq\mathbf E_{\ta^*} =
(\mathbf E_{\ta})^*\simeq \mathbf T(M)^*$.
\end{Proof}

Next, I want to prove  that our super-Grassmannian is, as a~rule, non-split.
Note that the canonical action of $G_0$ on $\M$ gives rise to a~natural
linear action of this groups on the tangent sheaf $\mathcal T$ leaving invariant
the $\mathbb Z_2$-grading. As a~result,
we get a~linear representation of $G_0$ in the cohomology groups of $\mathcal T$
and, in particular, in the Lie superalgebra $\mathfrak v\M$. The corresponding
linear representation of the Lie algebra $\mathfrak g_0$ is given by the formula
$u\mapsto\ad_{u^*}$.

\ssbegin{Proposition}\label{P5.4}
If $r\ge 2$ or $s\ge 2$, then $\mathfrak v\M_{\bar 0}^{G_0} = 0$.
\end{Proposition}

\begin{Proof}
Any $\de\in\mathfrak v\M_{\bar 0}^{G_0}$ determines a~$P$-invariant even
derivation of the superalgebra $\mathcal O_o$ (we denote it by the same
character $\de$), Clearly, $\de$ preserves the maximal ideal~
$m_o$. Consider the vector subspace $\widetilde m_o\subset m_o$, spanned by the
germs of local coordinates at $o$. By Proposition \ref{P5.2}, $R$ preserves
the even and the odd parts of this subspace, inducing in each part an
irreducible representation, and the germs of local coordinates constitute a
weight basis of $\widetilde m_o$ with respect to $T$ with the weights
$-\la_i +\la_{r+\al}$, where $1\le i\le r$ and $1\le\al\le s$. Note
that the remaining weights of the representation of $R$ in the whole $m_o$
are certain sums of these weights, and hence we see that the weight
subspace of $m_o$ corresponding to any of these weights is
two-dimensional (and lies in $\widetilde m_o$). Since $\de$ is even and
$P$-invariant, the germs of local coordinates are eigenvectors for $\de$.
Moreover,  the Schur lemma implies that
$$
\de(x_{i\al}) = ax_{i\al},\; \de(\xi_{i\al}) = b\xi_{i\al},
$$
where $a,b\in\mathbb C$. We have $a = 0$. Indeed, consider the
vector field $\widetilde\de = \si_0(\de)\in\mathfrak v(M,\gr\mathcal O)_0$ (see
Subsection~\ref{ss2.1}). Clearly, $\widetilde\de$ is $G_0$-invariant, too, and hence
determines
the $G_0$-invariant vector field $\al(\widetilde\de)$ (see \eqref{(2.7)}. But it is
well known (see, e.g., \cite{23}) that the standard action of $\GL_n(\mathbb C)$
on $M$ is \textit{asystatic}, i.e.,  $M$ has no non-zero holomorphic
$G_0$-invariant vector fields (for the origin of the term \textit{asystatic}, see \cite{GP1, GP2} and interesting references therein). This implies that $\widetilde\de(x_{i\al} +\mathcal J) = 0$.
Therefore,  $\de(x_{i\al})\in\mathcal J^2$, whence $a = 0$. Now we prove
that $b = 0$, using the relation $[\de,v^*] = 0$ for all $v\in\mathfrak n_-$.
Proposition \ref{P5.1} implies that
$$
0 = [\de,E_{r+1,1}^*](x_{12}) = \de(E_{r+1,1}^*(x_{12})) =
\de(\xi_{11}\xi_{12}) = 2b\xi_{11}\xi_{12}.
$$
This implies our assertion whenever $s\ge 2$. To prove the assertion for
$r\ge 2$, one takes $x_{21}$ instead of $x_{12}$.
\end{Proof}

This result makes it possible to solve the splittness question
concerning the super-Grassmannians studied here.

\ssbegin[On splitness of $\Pi\Gr^{n|n}_{s|s}$]{Theorem}\label{T5.1}
The super-Grassmannian $\Pi\Gr^{n|n}_{s|s}$ is split if and only if $n = 2$ and $
s = 1$.
\end{Theorem}

\begin{Proof}
Consider the grading derivation $\ep$ of the $\mathbb Z$-graded sheaf
$\gr\mathcal O$ defined in Subsection~\ref{P2.2} and the natural homomorphism of Lie
superalgebras $\si_0: H^0(M,\mathcal T_{\bar 0})\tto H^0(M,\widetilde{\mathcal T}_0)$
defined in Subsection~\ref{ss2.1}. Proposition \ref{P5.4}
implies that $\ep\notin\Im\si$ whenever $s\ge 2$ or $r\ge 2$. Indeed, if
$\ep = \si_0(\de)$, where $\de\in H^0(M,\mathcal T_{\bar 0})$, then the complete
reducibility of the representation of $G_0$ in $H^0(M,\mathcal T_{\bar 0})$
implies that $\de$ can be chosen to be $G_0$-invariant. But then $\de = 0$,
whence $\ep = 0$, which gives a~contradiction. If $\M$ is split, then $\si$
is an isomorphism, but this is false whenever $s\ge 2$ or $r\ge 2$.
In the case $n = 2,\,s = 1$, we can see that the super-Grassmannian is split,
e.g., by calculating its transition functions.
\end{Proof}

An important property of $\Pi\Gr^{n|n}_{s|s}$ is the homogeneity, which we are
going to prove now.

\ssbegin[$\Pi\Gr^{n|n}_{s|s}$ is  homogeneous]{Proposition}\label{P5.5}
\begin{enumerate}
\item[$(1)$]
The canonical action of $\mathfrak{q}_{n}(\mathbb C)$ on the supermanifold
$\Pi\Gr^{n|n}_{s|s}$ is transitive.
\item[$(2)$]
The kernel of this action is $\langle I_{n|n}\rangle$.
\end{enumerate}
\end{Proposition}

\begin{Proof}
To prove (1), we have to calculate the vector fields $y^*$ corresponding to
certain odd elements of $\mathfrak{q}_{n}(\mathbb C)$. More precisely, take the matrix
$y = \begin{pmatrix} 0 & B \\ B & 0\end{pmatrix}\in\mathfrak{q}_{n}(\mathbb C)_{\bar 1}$, where
$$
B = \begin{pmatrix} 0 & Y \\ 0 & 0\end{pmatrix},
$$
$Y = (y_{i\al})$ being an $(r\x s)$-matrix. Denoting by $\ta$ an odd
parameter, we get
$$
\begin{pmatrix} I_n & \ta B \\ \ta B & I_n\end{pmatrix}\begin{pmatrix} X & \Xi\\
I_{s} & 0 \\ \Xi & X \\0 & I_s \end{pmatrix} =
\begin{pmatrix} X & \Xi + \ta Y\\ I_{s} & 0 \\ \Xi + \ta Y & X \\0 & I_s \end{pmatrix}.
$$
It follows that
$$
y^* = -\sum_{i,\al} y_{i\al}\pd{\xi_{i\al}}.
$$

Clearly, $\ev_o(y^*)$ span the vector space $T_o(\M)_{\bar 1}$. Since our
action is $\bar 0$-transitive, its transitivity follows from Proposition
\ref{P2.4}(2).

Let us denote by $\mathfrak q$ the kernel of our action. We see from
Proposition \ref{P5.1} that $I_{n|n}\in\mathfrak q$. Since $\mathfrak g_0 =
\mathfrak{gl}_n(\mathbb C)$ acts on $M$ in the standard way, it follows that
$\mathfrak q\cap\mathfrak g_0 = \langle I_{n|n}\rangle$. But it is known
(see, e.g., \cite{15}) that the only ideal of $\mathfrak{q}_{n}(\mathbb C)$ containing
$\langle I_{n|n}\rangle$ is
$$
\mathfrak {sq}_{n}(\mathbb C) := \left\{\begin{pmatrix} A & B\\ B & A\end{pmatrix}\,\vrule\,
\tr B = 0\right\}.
$$
As we have seen above, $\mathfrak q\ne\mathfrak{sq}_{n}(\mathbb C)$. Hence, 
$\mathfrak q = \langle I_{n|n}\rangle$.
\end{Proof}

Now we are able to prove our main result concerning $\Pi$-symmetric
super-Grassman\-nians.

\ssbegin{Theorem}\label{T5.2}
Let $\M = \Pi\Gr^{n|n}_{s|s}$ and $n\ge 3$.
\begin{enumerate}
\item[$1)$]
In the classification of non-split supermanifolds with retract $(\Gr^{n}_{,s},
\Om)$ given by Theorem $\ref{T4.5}$, $\M$ corresponds to the invariant $(2,1)$-form
$\et$.
\item[$2)$]
The natural action of the Lie superalgebra $\mathfrak{q}_{n}(\mathbb C)$ on $\M$
determines an isomorphism of Lie superalgebras
$$
\mathfrak v\M\simeq\fp\fq_{n}(\mathbb C):=\mathfrak{q}_{n}(\mathbb C)/\langle I_{n|n}\rangle.
$$
\item[$3)$]
If $2\ge s\ge n-2$, then
$$
\begin{aligned}
H^1(M,\mathcal T_{\bar 0})&\simeq\mathfrak{sl}_n(\mathbb C)\oplus\mathbb C,\\
H^1(M,\mathcal T_{\bar 1})&\simeq\mathfrak{sl}_n(\mathbb C).
\end{aligned}
$$
\item[$4)$]
If $s = 1$ or $n-1$, then
$$
H^1(M,\mathcal T) = H^1(M,\mathcal T_{\bar 1})\simeq\begin{cases} 0\text{ for } n\ge 4\\
\mathbb C\text{ for } n = 3\end{cases}.
$$
\end{enumerate}
\end{Theorem}

\begin{Proof}
By Corollary \ref{Cor3.3}, 
the supermanifold corresponding to the form
$a\th_2 + b\et$ in the case II cannot be homogeneous if $a\ne 0$.
By Proposition \ref{P5.5}(1), this implies (1).

By Proposition  \ref{P5.5}(2), the natural action of $\mathfrak{q}_{n}(\mathbb C)$ on $\M$
induces an injective homomorphism $\mathfrak{q}_{n}(\mathbb C)/\langle I_{n|n}\rangle\tto\mathfrak v\M$. Comparing this with Theorem \ref{T4.7}(4), we see that this
homomorphism is surjective. Thus, (2) is proved.

The assertions (2) and (3) follow from (1) and Theorem \ref{T4.7}.
\end{Proof}

Theorems \ref{T5.2}(1) and \ref{T4.5} imply

\ssbegin[A family of deformations of $\Pi\Gr^{n|n}_{s|s}$]{Corollary}
The supermanifold $\Pi\Gr^{n|n}_{s|s}$ for $2\le s\le n-2$ is included in  a 
$1$-parameter family of mutually non-isomorphic supermanifolds with the same
retract. In particular, it is not rigid.
\end{Corollary}

To conclude, we note that these properties of the $\Pi$-symmetric
super-Grassmannians contrast with
the rigidity of certain other series of super-Grassmannians (see Examples
\ref{E1.9}, \ref{E1.10}, \ref{E1.12}). Let us denote by $\M$ one of these super-Grassmannians,
by $(M,\mathcal O_{\gr})$ its retract and by $\mathcal T,\; \mathcal T_{\gr}$ the
corresponding tangent sheaves.

\ssbegin[Rigid super-Grassmannians]{Theorem}\label{T5.3}
Suppose that $\M$ is one of the following supermanifolds:
$$
\begin{aligned}
&\Gr^{n|m}_{k|l}\text{ with } 0 < k <m,\; 0 < l < m,\;\\
(k,l)\ne &(1,n-1),\;(m-1,1),\;(1,n-2),\;(m-2,1),\;(2,n-1),\;(m-1,2);\\
&\I\Gr^{2r|r}_{2s|s}\text{ with } r\ge 2,\; (r,s)\ne (2,1);\\
&\I_{\odd}\Gr^{n|n}_{s|n-s}\text{ with } 4\le s\le n-3.
\end{aligned}
$$
Then, $\M$ is the only non-split supermanifold with retract
$(M,\mathcal O_{\gr})$ and, moreover, $\M$ is rigid.
\end{Theorem}

\begin{Proof}
It is known that in all the cases listed above we have
$$
H^1(M,(\mathcal T_{\gr})_p) =\begin{cases}\mathbb C&\text{if } p = 2\\ 0&\text{otherwise}
\end{cases}
$$
(see \cite[Theorem 1]{22} for $\M = \Gr^{n|m}_{k|l}$, \cite[Theorem 1]{26} for $\M =
\I\Gr^{2r|r}_{2s|s}$, \cite[Theorem~ 1]{27} for $\M = \I_{\odd}\Gr^{n|n}_{s|n-s}$.
Moreover,
it was proved in these papers that the supermanifolds $\M$ are non-split. By
Proposition \ref{P3.3}, $\M$ is the only non-split supermanifold with retract
$(M,\mathcal O_{\gr})$, and the corresponding class $\la_2^*(\ga)$ is a~basic
element of $H^1(M,(\mathcal T_{\gr})_2)$. As in the proof of Theorem \ref{T4.7}, we
have $d_2(\ep) = - 2\la_2^*(\ga)$ in the spectral sequence $(E_r)$. 
Theorem \ref{T3.5} implies that $H^1(M,\mathcal T) = 0$, and hence $\M$ is rigid.
(The vanishing of $H^1(M,\mathcal T)$ was proved in the cited papers as well.)
\end{Proof}

Note that the super-Grassmannians listed in Theorem~\ref{T5.3}, together
with the $\Pi$-sym\-met\-ric super-Grassmannians, are just the supermanifolds
of flags that can be called symmetric superspaces (see \cite{32}).

\subsection*{Acknowledgements} The author wishes to thank the Laboratory of Mathematics of the
University of Poitiers (France) and the E. Schr\"odinger International
Institute for Mathematical Physics (Vienna, Austria) for the hospitality
during the spring semester 1996, when this paper was written. I am
grateful to P. Torasso, P. Michor and M. Eastwood for valuable discussions.

\end{paper}
\begin{references}

\bibitem[1]{1}
Akhiezer~D.~N.,
Lie Group Actions in Complex Analysis.
Vieweg,
Braunschweig/Wiesbaden,
1995.

\bibitem[2]{2}
Berezin~F.~A.,
\textit{Introduction to Superanalysis}.
Edited and with a foreword by A. A.~Kirillov. With an
appendix by V. I.~Ogievetsky. Translated from the Russian by
J.~Niederle and R.~ Koteck\'y. Translation edited by D.~Leites.
Mathematical Physics and Applied Mathematics, 9. D.~Reidel
Publishing Co., Dordrecht. xii+424 pp.

\bibitem[3*]{2b}
Berezin~F.~A.,
\textit{Introduction to Superanalysis}. 2nd edition revised and edited by D.~Leites and with
appendix by D.~Leites, V.~Shander, I.~Shchepochkina ``Seminar on
Supersymmetries v. 1$\frac12$"), MCCME, Moscow, 432 pp. (in Russian)


\bibitem[4]{3}
Bott~R.,
Homogeneous vector bundles.
Ann. Math.
V. 66.
(1957)
203--248.

\bibitem[5*]{BV}
Borovoi~M., Vishnyakova~E., Automorphisms and real structures for a $\Pi$-symmetric super-Grassmannian,  \url{https://arxiv.org/abs/2205.04380}.



\bibitem[6*]{BGLS}
Bouarroudj~S., Grozman~P., Leites~D., Shchepochkina~I..
Minkowski superspaces and superstrings as almost real-complex supermanifolds.  Theor. and Mathem. Physics, Vol. 173, no.~3, (2012) 1687--1708; \url{https://arxiv.org/abs/1010.4480}.

\bibitem[7]{4}
Bunegina~V.~A., Onishchik~A.~L.,
Homogeneous supermanifolds associated with the complex projective
line.
J. Math. Sc.
V. 82
(1996)
3503--3527.

\bibitem[8*]{Del}
Deligne~P., Etingof~P., Freed~D., Jeffrey~L., Kazhdan~D., Morgan~J., Morrison~D., Witten~E., (eds.). \textit{Quantum fields and strings: a~course for
mathematicians}. Vol. 1. Material from the Special Year on
Quantum Field Theory held at the Institute for Advanced Study,
Princeton, NJ, 1996--1997. American Mathematical Society,
Providence, RI; Institute for Advanced Study (IAS), Princeton, NJ,
1999. Vol. 1: xxii+723 pp.

\bibitem[9*]{DW}
Donagi,~R., Witten,~E., Super Atiyah classes and obstructions to splitting of supermoduli space. Pure Appl. Math. Q. 9 (2013), no. 4, 739--788.

\bibitem[10]{5}
Eastwood~M., LeBrun~C.,
Thickening and supersymmetric extensions of complex manifolds.
Amer. J. Math.
V. 108
(1986)
1177--1192.


\bibitem[11]{6}
Eastwood~M., LeBrun~C.,
Fattening complex manifolds: curvature and Kodaira--Spencer maps.
J. Geometry and Physics.
V. 8
(1992)
123--146.


\bibitem[12]{7}
Fr\"olicher~A., Nijenhuis~A.,
Theory of vector-valued differential
forms, P.1. Derivations in the graded ring of differential forms.
Proc. Kon. Ned. Akad. Wet. Amsterdam.
V. 59
(1956)
540--564.


\bibitem[13]{8}
Fr\"olicher~A., Nijenhuis~A.,
Some new cohomology invariants of complex manifolds.
Proc. Kon. Ned. Akad. Wet. Amsterdam.
V. 59
(1956)
540--564.

\bibitem[14*]{GM}
Gelfand~S., Manin~Yu., \textit{Methods of Homological Algebra} 2nd Edition (2003) 372pp.

\bibitem[15*]{GP1}
Gori~A.,  Podest\`a~F., Complex asystatic actions of compact Lie groups;
\url{https://arxiv.org/pdf/math/0411204.pdf}.

\bibitem[16*]{GP2}
Gori~A,  Podest\`a~F., Symplectically asystatic actions of compact Lie groups. Transform. Groups 11 (2006), no. 2, 177--184.

\bibitem[17]{9}
Green~P.,
On holomorphic graded manifolds.
Proc. Amer. Math. Soc.
V. 85
(1982)
 587--590.


\bibitem[18]{10}
Griffiths~P.~R., Harris~J.,
Principles of Algebraic Geometry.
J. Wiley \& Sons,
New York,
1978.


\bibitem[19]{11}
Grothendieck~A.,
Sur quelques points d'alg\`ebre homologique.
T\^ohoku Math. J.
V. 9
(1957)
 119--221.

\bibitem[20*]{Gz}
Grozman~P., Classification of bilinear invariant operators on tensor
fields. Functional Anal. Appl., v.~14 (1980), no.~2, 127--128; for
proofs, see \url{https://arxiv.org/abs/math/0509562} and this Special Volume, pp.

\bibitem[21]{12}
Helgason~S.,
Differential Geometry, Lie Groups, and Symmetric Spaces.
Academic Press,
New York e.a.,
1978.


\bibitem[22]{13}
Humphreys~J.~E.,
Introduction to Lie Algebras and Representation Theory.
Springer-Verlag,
New York e.a.,
1972.



\bibitem[23]{14}
Ise~M.,
Some properties of complex analytic vector bundles over compact
complex homogeneous spaces.
Osaka Math. J.
V. 12
(1960)
 217--252.


\bibitem[24]{15}
Kac~V.~G.,
Lie superalgebras.
Adv. Math.
V. 26
(1977)
 8--96.


\bibitem[25]{16}
Kol\'a\v r~I., Michor~P.~W., Slov\'ak~J.,
Natural Operations in Differential Geometry.
Springer-Verlag,
Berlin e.a.,
1993.


\bibitem[26]{17}
Kostant~B.,
Lie algebra cohomology and the generalized Borel--Weil theorem.
Ann. of Math.
V. 74
(1961)
 329--387.


\bibitem[27]{18}
Kostant~B.,
Graded manifolds, graded Lie theory, and prequantization. In:
Lecture Notes in Math.
Springer-Verlag,
Berlin e.a., 1977, 570 pp., 177--306.


\bibitem[28]{19}
Koszul~J.-L.,
Sur la forme hermitienne canonique des espaces homog\`enes complexes.
Canad. J. Math.
V. 7
(1955)
 562--576.

\bibitem[29*]{Lsos}
Leites~D. (ed.), Bernstein~J., Leites~D., Molotkov~V., Shander~V.. \textit{Seminar on supersymmetries. Vol. $1$: Algebra
and Calculus on supermanifolds}, MCCME, 
Moscow, 2011, 410~pp. (in Russian; a~draft in English is available) \url{https://staff.math.su.se/mleites/books/2011-sos1.pdf}.


\bibitem[30]{20}
Manin~Yu.~I.,
Gauge Field Theory and Complex Geometry.
Springer-Verlag,
Berlin e.a.,
1988. Second
edition. Springer-Verlag, Berlin, (1997) xii+346 pp.


\bibitem[31*]{MaAG}
Manin~Yu.~I.,
 \textit{Introduction to the theory of schemes}. Translated from the Russian and edited by D.~Leites, Springer, (2018) xvi+201 pp.
 
 
\bibitem[32*]{Meng}
Meng~L., Leray-Hirsch theorem and blow-up formula for Dolbeault cohomology; \url{https://arxiv.org/abs/1806.11435}.
 
\bibitem[33*]{Mont}
Montgomery~D., Simply connected homogeneous spaces. Proc. Amer. Math. Soc. 1 (1950), 467--469

\bibitem[34*]{Mo}
Molotkov~V., Infinite-dimensional and colored supermanifolds (In Russian; an expounded version of J. Nonlinear Math. Physics, 17: sup1 (2010), 375--446). \url{https://doi.org/10.1142/S140292511000088X}.

\bibitem[35*]{21}
Onishchik~A.~L.,
Transitive irreducible Lie superalgebras of vector fields. In: Leites D. (ed.) \textit{Seminar on Super \-ma\-ni\-folds}.
Reports Dep. Math. Univ. Stockholm.
V. 26 (1987) 1--21. See this Special Volume, pp.


\bibitem[36]{OV}
Onishchik~A., Vinberg~E., \textit{Lie groups and algebraic groups.}
Springer-Verlag, Berlin, 1990. xx+328 pp.

\bibitem[37]{22}
Onishchik~A.~L.,
On the rigidity of super-Grassmannians.
Ann. Global Analysis and Geom.
V. 11
(1993)
 361--372.


\bibitem[38]{23}
Onishchik~A.~L.,
Topology of Transitive Transformation Groups.
J.A. Barth Verlag,
Leipzig e.a.,
1994.


\bibitem[39]{24}
Onishchik~A.~L.,
A spectral sequence for the tangent sheaf cohomology of a
supermanifold.
SFB 288, Preprint No. 148,
Berlin 1994.


\bibitem[40]{25}
Onishchik~A.~L.,
About derivations and vector-valued differential forms.
Inst. Math. Ruhr-Universit\"at Bochum. Bericht N 181.
Bochum 1995.

\bibitem[41*]{OP1}
Onishchik~A.~L., Platonova~O.~V., Homogeneous supermanifolds associated with a~complex projective space. I. Sb. Math. 189 (1998), no. 1-2, 265--289.

\bibitem[42*]{OP2}
Onishchik~A.~L., Platonova~O.~V., Homogeneous supermanifolds associated with the complex projective space. II. Sb. Math. 189 (1998), no. 3-4, 421--441.



\bibitem[43]{26}
Onishchik~A.~L., Serov~A.~A.,
Vector fields and deformations of isotropic super-Grassmannians of
maximal type. In:
Lie Groups and Lie Algebras: E.B. Dynkin`s Seminar. Amer.
Math. Soc. Transl. Ser. 2. V. 169.
Amer. Math. Soc.,
Providence, R.I.,
1995.
 75--90.


\bibitem[44]{27}
Onishchik~A.~L., Serov~A.~A.,
On isotropic super-Grassmannians of maximal type associated with an
odd bilinear form.
E. Schr\"odinger Inst. for Math. Physics, Preprint No. 340.
Vienna 1996.


\bibitem[45]{28}
Palamodov~V.~P.,
Invariants of analytic $\mathbb Z_2$-manifolds.
Funct. Anal. Appl., 17:1 (1983), 68--69


\bibitem[46]{29}
Platonova~O.~V.,
Homogeneous supermanifolds associated with the projective space
(Thesis).
Yaroslavl Univ.
Yaroslavl 1995 (in Russian; see \cite{OP1} and \cite{OP2}).


\bibitem[47]{30}
Rothstein~M.~J.,
Deformations of complex supermanifolds.
Proc. Amer. Math. Soc.
V. 95
(1985)
 255--260.


\bibitem[48]{31}
Scheunert~M.,
The Theory of Lie Superalgebras. Lecture Notes in Math. 716,
Springer-Verlag,
Berlin e.a.,
1979. 217 pp.


\bibitem[49]{32}
Serganova~V.,
Classification of simple real Lie superalgebras and of symmetric
superspaces.
Funct. Anal. Appl., 17:3 (1983), 200--207.


\bibitem[50*]{33}
Vaintrob~A.~Yu., Deformations of complex structures on supermanifolds, Funct. Anal. Appl., 18:2 (1984), 135--136.

\bibitem[51*]{34}
Vaintrob~A.~Yu., Deformations of complex superspaces and of the coherent sheaves on them, J. Math. Sci. (N. Y.), 51:1 (1990), 2140--2188.

\bibitem[52*]{V6}
Vishnyakova~E.,  On complex Lie supergroups and homogeneous split supermanifolds, Transformation Groups, v. 16 (2011) 265--285 
\url{https://arxiv.org/abs/0908.1164}.

\bibitem[53*]{V5} 
Vishnyakova~E., Vector fields on $\Pi$-symmetric flag supermanifolds, S\~{a}o Paulo J. of  Math. Sciences, v. 10 (2016) 20--35; \url{https://arxiv.org/abs/1506.02295}.

\bibitem[54*]{V4}
Vishnyakova~E., Vector fields on $\fosp_{2m|2n}(\Cee)$- and $\pi\fsp_n(\Cee)$-flag supermanifolds, J. Algebra, v. 459 (2016) 1--28; 
\url{https://arxiv.org/abs/1706.00128} .

\bibitem[55*]{V3}
Vishnyakova~E., Vector fields on $\fosp_{2m-1|2n}(\Cee)$-flag supermanifolds, Commun. in Algebra, v. 47, (2019) 4247--4261; \url{https://arxiv.org/abs/1801.09239}.

\bibitem[56*]{V2}
Vishnyakova~E., Rigidity of flag supermanifolds,  Transformation groups,
\url{https://doi.org/10.1007/s00031-020-09629-6b}; \url{https://arxiv.org/abs/1908.11753}.
 

\end{references}
